%% file: problems.tex
\documentclass[oneside,11pt]{book}
\usepackage{xr,paralist,verbatim,xcite}

\usepackage{fancyhdr}
\usepackage{titlepic} 
\usepackage{amsmath,amsfonts,empheq,latexsym,epsfig,subfigure,fancybox,verbatim,paralist,mdframed,framed,amsthm,algorithm,isomath,multirow,makeidx,cancel,enumitem,algorithmic,minitoc,url,varioref,datetime,setspace,array,listings,enumitem}

\usepackage{pstricks-add,pst-sigsys,pst-tree,pst-bezier,pst-plot}
\input{myproblemcommands}

\usepackage{geometry}
\newcommand{\zt}{z_{t}^{(1)}}
\newcommand{\st}{\theta}
\newcommand{\ztp}{z_{t+1}^{(1)}}
\newcommand{\zht}{\hat{z}_{t}^{(1)}}

\newcommand{\ut}{r_l}
\renewcommand{\val}{\operatorname{val}}

\newtheorem*{theorem-non}{Theorem}

\newcommand{\act}[1]{u^{(#1)}}

\newcommand{\pol}[1]{\policy^{(#1)}}

\newcommand{\Jc}[1]{J^{(#1)}}
\newcommand{\corr}{\mathcal{C}}

\newcommand{\sr}[1]{S^{(#1)}}

\newcommand{\btheta}{R}

\newcommand{\lyap}{V}

\DeclareMathOperator{\dis}{dist}

\usepackage{lmodern}
\usepackage{tikz}
\usepackage{xcolor}
\usetikzlibrary{shapes,arrows}
\usetikzlibrary{decorations.pathmorphing, decorations.text}

\newcommand{\unitvec}{e}
\newcommand{\globbehav}{z}
\newcommand{\actprof}{\mathbf{u}}
\newcommand{\dtimee}{n}

\title{\bf Partially Observed Markov Decision Processes \\  {\Large From Filtering to Controlled Sensing} \\ \mbox{} \\ {\em Internet Supplement}}
\author{Vikram Krishnamurthy\\ University of British Columbia, Vancouver, Canada.\\ {This version: \today } \\ {\em This material is copyrighted.} \\
\copyright Vikram Krishnamurthy }
\titlepic{\includegraphics[width=1.12\textwidth]{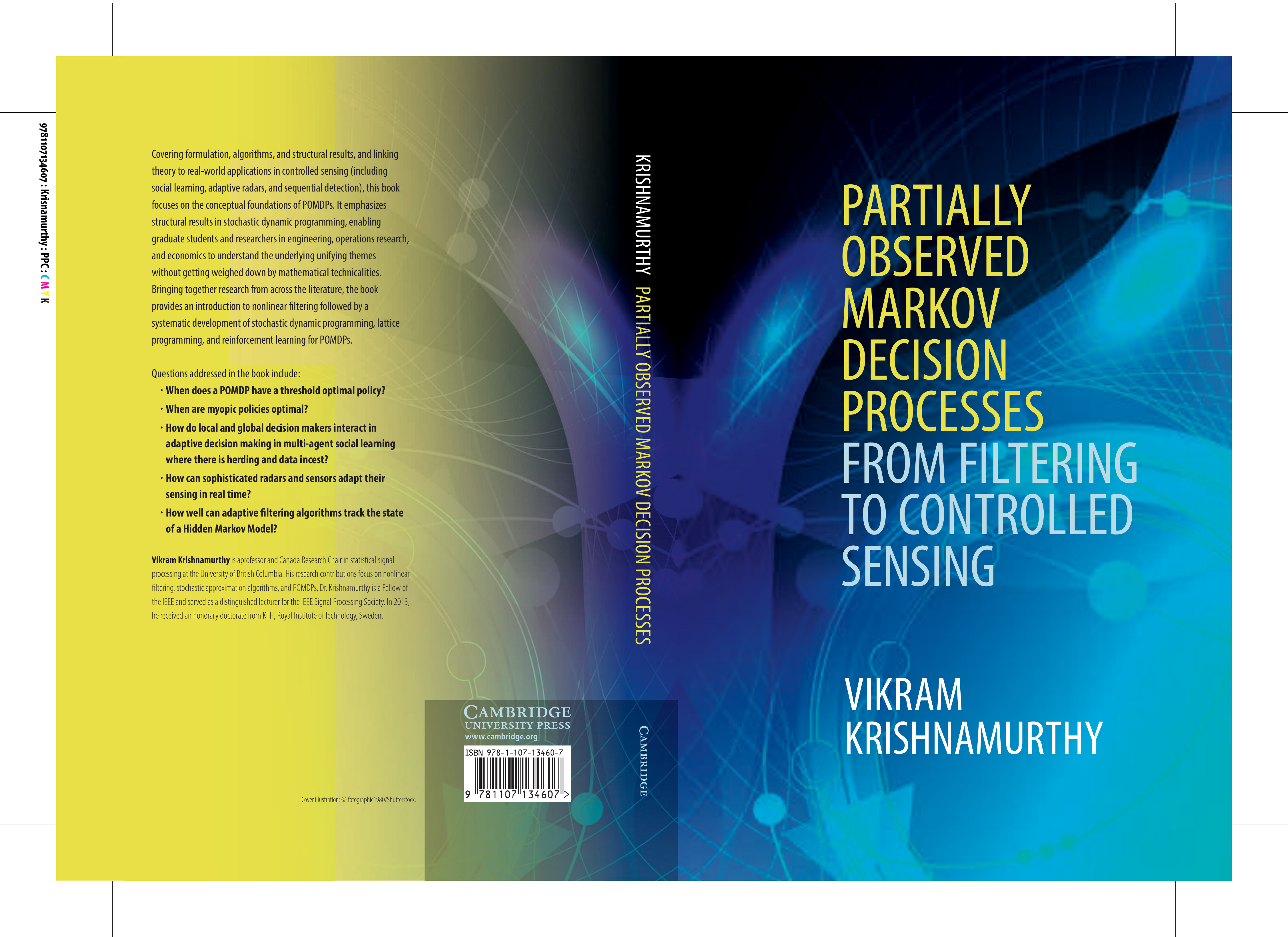}}
\date{}

\makeindex
\usepackage{chngcntr}
\counterwithout{equation}{chapter}

\begin{document}

\maketitle

\tableofcontents

\newpage  \thispagestyle{empty}
\chapter*{Preface to Internet Supplement}
This document is an 
internet supplement to my book ``Partially Observed Markov Decision Processes -- From Filtering to Controlled Sensing'' published
by Cambridge University Press in 2016.\footnote{Online ISBN:9781316471104
and 
Hardback ISBN:9781107134607}

This internet supplement  contains exercises,  examples and case studies.
The  material  appears in this internet supplement (instead of the book) so that  it can be  updated.
 This document will evolve over time  and further discussion and examples will be added.

The website \url{http://www.pomdp.org} contains downloadable  software for solving POMDPs and several examples of POMDPs.
I have found that by interfacing the POMDP solver with Matlab, one can solve several interesting types of POMDPs such as those 
with nonlinear costs (in terms of the information state) and bandit problems. 

I have given considerable thought into designing  the exercises and case studies in this internet supplement.
They are mainly mini-research type exercises rather than simplistic drill type exercises. Some of the problems
are extensions of the material in the book. As can be seen from the content list, this  document also contains some short (and in some cases, fairly incomplete) case studies which will be made more detailed over time.
These  case studies were put in this internet supplement in order to keep the size of the book manageable.
As time progresses, I hope to incorporate additional   case studies and other pedagogical notes  to this document to  assist in understanding some of the material in the book.  Time permitting, future plans include adding a detailed discussion on structural results for POMDP games; structural results for quasi-variational
inequalities, etc.

To avoid confusion in numbering, the equations in this internet supplement are numbered consecutively starting from (1) and not chapter wise.
In comparison, the equations in the book are numbered chapterwise.

 This internet supplement document is work in progress and will be updated periodically.
I welcome   constructive comments from readers of the book and this internet supplement.
Please email me at {\tt vikramk@ece.ubc.ca} \\
\\ \\

\hfill
\begin{minipage}{6cm}
Vikram Krishnamurthy,\\
 2016
 \end{minipage}

\newpage

\setcounter{chapter}{1}

\chapter{Stochastic State Space Models}   

\begin{compactenum}

\item Theorem \ref{thm:pf} dealt with the stationary distribution and eigenvalues of a stochastic matrix (transition  probability matrix of a Markov chain).
Parts of Theorem \ref{thm:pf} can be shown via elementary linear algebra.

{\bf Statement 2}:
 Define spectral radius  $\bar{\lambda}(\tp) = \max_i|\lambda_i|$ \\
{\em Lemma} : $\bar{\lambda}(\tp)  \leq \|\tp\|_\infty$ where $\|\tp\|_\infty =\max_i \sum_j \tp_{ij} $\\
{\em Proof}:  For all eigenvalues $\lambda$, $|\lambda | \|x\| = \|\lambda x\| = \|\tp  x \|  \leq \|\tp \| \|x \| \implies
 |\lambda| < \|\tp\|$.

For a stochastic matrix, $\|\tp\|_\infty =1 $  and $\tp$ has an eigenvalue at 1. So $\bar{\lambda}=1$.

{\bf Statement 3}: For non-negative matrix $\tp$, $\tp^\p \pi = \pi $ implies $\tp^\p|\pi| = |\pi|$ where $|\pi|$ denotes the vector with element-wise  absolute values. \\
Proof: $|\pi| = |\tp^\p \pi| \leq |\tp^\p| |\pi| = \tp^\p |\pi|$
So $\tp^\p|\pi| - |\pi| \geq 0$.\\
But $\tp^\p|\pi| - |\pi| > 0$ is impossible, since it implies $1^\p \tp^\p |\pi| > 1^\p |\pi|$,
i.e., $1^\p |\pi| > 1^\p |\pi|$.

\item Farkas' lemma  \index{Farkas' lemma}  is a widely used result in linear algebra. It states:  Let $M$ be an $m \times n$ matrix and $b$ an $m$-dimensional vector. Then only one of the following statements is true:
\begin{enumerate}
\item There exists a vector $x\in \reals^n$ such that $M x = b$ and $ x \geq 0$. 
\item There exists a vector $y\in \reals^m$ such that $M^\p y \geq 0$ and $b^\p y < 0$.
\end{enumerate}
Here $x\geq 0$ means that all components of the vector $x$ are non-negative.

Use Farkas lemma to prove that every transition matrix $\tp$ has a stationary distribution. That is, for any $\statedim \times\statedim$ stochastic matrix $\tp$, there
exists a probability vector $\belief$ such that $\tp^\p \belief = \belief$. (Recall a probability vector $\belief$ satisfies
$\belief(i) \geq 0$, $\sum_i \belief(i) = 1$).

Hint: Write alternative (a) of Farkas lemma as
$$  \begin{bmatrix} (\tp - I)^\p \\ \one^\p \end{bmatrix} \belief = \begin{bmatrix}  0_{\statedim} \\ 1 \end{bmatrix} , \quad \belief > 0$$
Show that this has a solution by demonstrating that alternative (b) does not have a solution.

\item
Using the maneuvering target model of Chapter \ref{chp:manuever},  \index{maneuvering target}
simulate the dynamics and measurement process of a target with the following specifications:\\
{\small 
\begin{tabular}{|l|ll|} \hline \hline
Sampling interval & $\samt = 7$ s & \\ \hline
Number of measurements & $\finaltime = 50$ &\\ \hline
Initial target position & $(-500 , -500)'$ m &\\ \hline
Initial target velocity & $(0.0, 5.0)'$ m/s &\\ \hline
Transition probability matrix & $\tp_{ij} = 
\left\{\begin{array}{cc} 0.9  \mbox{ if } i=j\\
                          0.05 \mbox{ otherwise }\end{array}\right.$ &\\
\hline 
Maneuver commands  (three)& $\inpm \mc = \left( \begin{array}{llll}
                                0 & 0  & 0 & 0
                                \end{array} \right)'$   &(straight)
                        \\
&                       $\inpm \mc = \left( \begin{array}{llll}
                                -1.225 &  -0.35  & 1.225 & 0.35
                                \end{array} \right)'$  &(left turn)
                        \\
&                       $\inpm \mc = \left( \begin{array}{llll}
                                1.225 &  0.35 & -1.225 & -0.35
                                \end{array} \right)'$ &(right turn) \\ \hline
Observation matrix & $\obsm=\identity_{4\times 4}$ &\\ \hline
Process noise      & $\snoisecov=(0.1)^2 \identity_{4\times 4}$ &\\ \hline
Measurement noise      & $\onoisecov= \mbox{diag}( 20.0^2, 1.0^2, 20.0^2, 1.0^2)$ &\\ \hline
Measurement volume $V$ & $[-1000, 1000]$ m in $x$ and $y$ position &\\ 
                   & $[-10.0, 10.0]$ m/s in $x$ and $y$  velocity &
\\ \hline \hline
\end{tabular} }

\item Simulate the optimal predictor via the composition method. The composition method is discussed in \secn \ref{sec:jmlspred}.

\item As should be apparent from an elementary linear systems course, the algebraic Lyapunov equation (\ref{eq:algebraic_lyapunov})
is intimately linked with the stability of a linear discrete time system. Prove that $A$ has all its eigenvalues strictly inside the unit circle
iff  for every positive definite matrix $Q$, there exists a positive definite matrix $\kalmancov_\infty$ such that (\ref{eq:algebraic_lyapunov}) holds.

\item Theorem \ref{thm:dobproperty}  states that 
$|\lambda_2| \leq \dob(\tp) $. That is, the Dobrushin coefficient  upper bounds the second largest eigenvalue modulus of a stochastic matrix $\tp$. 
Show that 
$$
\log |\lambda_2|  = \lim_{k\rightarrow \infty}  \frac{1}{k} \log  \dob(\tp^k) $$

\item 
Often for sparse transition matrices, $\dob(\tp)$ is typically equal to 1 and therefore not useful since it provides a trivial upper bound for $|\lambda_2|$.
For example, consider a random walk characterized by the tridiagonal transition matrix
$$ \tp = \begin{bmatrix}    
r_0 &  p_0   & 0  & 0 & \cdots &  0 \\
q_1 & r_1 & p_1 &  0 & \cdots & 0 \\
0 & q_2 & r_2 & p_2 & \cdots & 0 \\
\vdots & & \ddots & \ddots &  \ddots & \vdots \\
0 & \cdots & 0 & q_{\statedim-1} &  r_{\statedim-1}  & p_{\statedim-1} \\
0 & \cdots & 0  & 0  &  q_{\statedim} &   r_{\statedim}  \end{bmatrix}
$$
Then using  Property 3 of $\dob(\cdot)$ above, clearly
$\sum_l \min  \{ \tp_{il}, \tp_{jl}  \} = 0$, implying that
 $\dob(\tp) = 1$. So for this example, the Dobrushin coefficient does not say anything about the initial condition being forgotten geometrically  fast.
 
 For such cases, it is often useful to consider the Dobrushin coefficient of powers of $\tp$.
In the above example, clearly every state communicates with every other state in at least $\statedim$ time points. So
 $\tp^\statedim$ has strictly positive elements.  Therefore $\dob(\tp^\statedim)$ is strictly smaller than 1 and is a useful bound. Geometric ergodicity follows
 by consider blocks of length $\statedim$, i.e.,  
 $$ \dvar{  { \tp^\statedim}^\p \belief}   {{\tp^\statedim}^\p \bbelief } \leq \dob(\tp^\statedim)  \dvar{\belief}{\bbelief} $$

\item Show that the inhomogeneous Markov chain with transition matrix
$$ \tp(2n-1) = \begin{bmatrix} 0.5 & 0.5 \\ 1 & 0 \end{bmatrix} ,\quad   \tp(2n) = \begin{bmatrix} 0 & 1 \\ 1 & 0   \end{bmatrix}$$
is weakly ergodic.

\item {\bf Wasserstein distance.}  
\index{Wasserstein distance}
As mentioned in \secn \ref{sec:modelnotes} of the book, the Dobrushin coefficient is a special case of  a more general coefficient of ergodicity.
This general definition is in  terms of the 
Wasserstein metric which we now define:  Let $d$ be a metric on the state space $\statespace = \{e_1,e_2,\ldots,\}$ where the state space is possibly denumerable.
Consider the bivariate random vector  $(x,y) \in \statespace\times \statespace$ with marginals $\belief_x$ and $\belief_y$, respectively.

Define the Wasserstein distance as
$$ d(\belief_x,\belief_y) = \inf \E\{ d(x,y) \} $$
where the infimum is over the joint distribution of $(x,y)$.

\begin{compactenum}
\item Show that the variational distance is a special case of the  Wasserstein distance obtained by choosing $d(x,y)$ as the discrete metric
$$ d(x,y) = \begin{cases} 1 & x \neq y \\ 0  & x = y .\end{cases} $$
\item  Define the coefficient of ergodicity associated with the Wasserstein distance as 
$$ \dob(\tp) = \sup_{i \neq j} \frac{d( \tp^\p e_i, \tp^\p e_j)}{d(e_i,e_j)} $$
Show that the Dobrushin coefficient is a special case of the above coefficient of ergodicity corresponding to the discrete metric.

\item Show that the above coefficient of ergodicity satisfies properties 2, 4 and 5 of Theorem \ref{thm:dobproperty}.
\end{compactenum}

\item {\bf Ultrametric transition matrices.}  \index{ultrametric matrix}
It is trivial to verify that $\tp^n$ is a stochastic matrix for any integer $n \geq 0$.
Under what conditions is $\tp^{1/n}$ a stochastic matrix?
A symmetric  ultrametric stochastic matrix $\tp$ defined in \secn \ref{sec:blackwelldom} of the book satisfies this property.

\end{compactenum}

\chapter{Optimal Filtering}   

\section{Problems}
\begin{compactenum}
\item  Standard drill exercises include: 
\begin{compactenum}
\item Compare  via simulations the recursive least squares with the  Kalman filter
\item Compare via simulations the recursive least square and the least mean squares
 (LMS) algorithm with a HMM filter when tracking a slow Markov chain. Note that  Chapter \ref{sec:hmmsa} of the book gives performance bounds on how well a LMS algorithm can track a slow Markov chain.
 \item
Another standard exercise is to try out variations of the particle filter with different importance distributions and resampling strategies on different models.
Compare via simulations the cubature filter, unscented Kalman filter and a particle filter for a bearings only target tracking model.
\item
A classical result involving the Kalman filter is the so called innovations state space model representation and the associated spectral factorization
problem for the Riccati equation, see \cite{AM79}.
\item {\bf Posterior Cramer Rao bound.} The posterior Cramer Rao bound \cite{TMN98}  for filtering can be used to compute a lower bound to the mean square error.  \index{posterior Cramer Rao bound}
This requires twice differentiability of the logarithm of the joint density. For HMMs, one possibility is to consider the Weiss-Weinstein
bounds \index{Weiss-Weinstein bounds}, see \cite{RO05}.  Chapter~\ref{chp:filterstructure}
 of the book gives more useful sample path bounds on the HMM filter using
stochastic dominance.
\end{compactenum}

\item {\bf Bayes' rule interpretation of Lasso.}\cite{PC08}  \index{Lasso} Suppose that 
the state $\state \in \reals^\statedim$ is a random variable with prior pdf
$$ \pdf(\state) = \prod_{j=1}^\statedim \frac{\lambda}{2} \, \exp\left( -\lambda \state(j) \right) .$$
Suppose $\state$ is observed via the observation equation
$$ \obs = \statem \state + \onoise, \qquad \onoise \sim \normal(0,\sigma^2 I) $$
where $A$ is a known $n \times \statedim$ matrix.
The variance $\sigma^2$ is not known and has a prior pdf $\pdf(\sigma^2)$. Then show that the posterior of $(\state,\sigma^2)$ given 
the observation $\obs$  is of the form
$$ 
\pdf(\state,\sigma^2| \obs) \propto
\pdf(\sigma^2) \,   (\sigma^2)^{-\frac{n+1}{2}}  \exp\big( - \frac{1}{2\sigma^2}\, \text{Lasso}(\state,\obs,\mu)\,\big) $$
where $\mu = 2 \sigma^2 \lambda$ and 
$$ \text{Lasso}(\state,\obs,\mu) = \|\obs-\statem \state\|^2 + \mu \|\state\|_1 .$$
Therefore for fixed $\sigma^2$, computing the mode $\hat{x}$ of the posterior is equivalent to computing the minimizer $\hat{x}$ of $\text{Lasso}(\state,\obs,\mu) $.

The resulting Lasso (least absolute shrinkage and selection operator) estimator $\hat{x}$ was proposed in \cite{Tib96} which is one of the most influential papers 
 in statistics since the 1990s. Since $\text{Lasso}(\state,\obs,\mu) $ is convex in $\state$ it can be computed efficiently via convex optimization algorithms.

\item Show that if $X \leq Y$ (with probability 1), then $\E\{X|Z\} \leq \E\{Y|Z\}$ for any information $Z$.
\item  Show that for a linear Gaussian system (\ref{eq:filter_linears}), (\ref{eq:filter_linearobs}),
$$ \pdf(\obs_k|\obs_{1:k-1}) = \normal(\obs_k - \obs_{k|k-1},   \obsm_k  \kalmancov_{k|k-1} \obsm^\p_k + \onoisecov_k ) $$
where $\obs_{k|k-1}$ and $\kalmancov_{k|k-1}$ are defined in (\ref{eq:statepred}), (\ref{eq:covpred}), respectively.

\item Simulate in Matlab the HMM filter, and  fixed lag smoother. Study empirically how the error probability of the estimates
decreases with lag. (The filter is a fixed lag smoother with lag of zero). Please also refer to \cite{GSV05} for a very nice analysis
of error probabilities.

\item Consider a HMM where the Markov chain evolves slowly with transition matrix $\tp = I + \epsilon Q$ where $\epsilon $ is a small
positive constant and 
$Q$ is a generator matrix. That is $Q_{ii} < 0$, $Q_{ij} > 0$ and each row of $Q$ sums to zero.
Compare the performance of the HMM filter with the recursive least squares algorithm  (with an appropriate
forgetting factor chosen) for estimating the underlying state.

\item Consider the following Markov modulated auto-regressive time series model: 
$$ z_{k+1} = \statem(r_{k+1})\,z_k + \snoisem(r_{k+1})\,\snoise_{k+1} + \inpm(r_{k+1})\,\inp_{k+1} $$
where $\snoise_{k} \sim \normal(0,1)$, $\inp_k$ is a known exogenous input.
Assume the sequence $\{z_k\}$ is observed. 
 Derive an optimal filter for the underlying Markov chain $r_k$. 
 (In comparison to a jump Markov linear system,  $z_k$ is observed without noise in this problem. The optimal filter is very
 similar to the HMM filter).

 \item Consider a Markov chain $\state_k$ corrupted by iid zero mean Gaussian noise and a sinusoid:
  $$y_k = \state_k + \sin(k/100) + \onoise_k$$
Obtain a filtering algorithm for extracting $\state_k$ given the observations.
 
 \item {\bf Image Based Tracking.} The idea is to estimate the coordinates $z_k$ of the target by measuring its orientation $r_k$ in noise.
 For example an imager can determine which direction an aircraft's nose is pointing thereby giving useful information about which
 direction it can move.
 Assume that the target's orientation evolves according to a finite state Markov chain. (In other words, the imager quantizes the target
 orientation to one of a finite number of possibilities.) Then the model for the filtering problem is
 \[
 \begin{split}
 z_{k+1} &= \statem(r_{k+1})\,z_k + \snoisem(r_{k+1})\,\snoise_{k+1} \\
\obs_k &\sim \pdf(\obs | r_k) 
 \end{split} \]
 Derive the filtering expression for  $\E\{z_k | \obs_1,\ldots,\obs_k\}$. The papers \cite{SSD93,KE97a,EE99} consider image based filtering. \index{image-based tracking}
 
 \item Consider a jump Markov linear system.  Via computer simulations, compare the IMM algorithm, Unscented Kalman filter and  particle filter.

 
 \item {\bf Radar pulse train de-interleaving.}  \index{de-interleaving} 
 \index{jump Markov linear system! pulse de-interleaving} 
 In radar signal processing,  radar pulses are received from multiple periodic sources.
 It is of interest to estimate the periods of these sources. For example, suppose:
 \begin{compactitem} 
\item   source 1  pulses are received at times
$$ 2, 7, 12, 17, 22 , 27, 32, 37, 42,\ldots , \quad \text{    (period = 5, phase =2)} $$
\item source 2  pulses are received at times
$$ 4,  15, 26, 37, 48, 59, 70, 81,\ldots , \quad \text{  (period = 11, phase = 4)  }.$$\end{compactitem}  The  interleaved  signal consists of  pulses at times
$$2,4,7,12,15,22,26,27, 32, 37, 42,\ldots. $$  So the above interleaved signal contains time of arrival information.  
Note that at time 37, pulses are received from both sources; but it is assumed  that there is no amplitude information - so the received signal
is simply a time of arrival event at time 37.
At the receiver, interleaved signal (time of arrivals) is corrupted by jitter noise (modeled as iid noise).
 So the noisy received signal are,  for example, $$2.4, \; 4.1, \; 6.7, \; 11.4, \; 15.5, \; 21.9, \; 26.2, \; 27.5, \;  30.9, \; 38.2,\; 43.6, \ldots.$$
Given this noisy interleaved signal, the de-interleaving problem aims to determine which pulses came from which source.
This can be done by estimating the periods (namely, 5 and 11) and phases (namely, 2 and 4) of the 2 sources.

The de-interleaving problem can be formulated as a jump Markov linear system.
Define the state 
  $x_k' =  (T', \tau_k')$,
consists of the periods $T' = (T^{(1)}, \ldots, T^{(N)})$ of the $N$ sources and  
$\tau_k'=  (\tau_k^{(1)}, \ldots, \tau_k^{(N)})$, where $\tau_k^{(i)}$ denotes the last
time source $i$ was active up to and including the arrival of the $k$th pulse.
Let $\tau_1 =( \phi^{(1)}, \ldots, \phi^{(N)})$, be the phases of periodic pulse-train sources.
Then
\beq
 \tau_{k+1}^i = \left\{\begin{array}{cc} \tau_k^i + T^i & \mbox{if $(k+1)$th pulse is due to source $i$} \\
 \tau_k^i &  \mbox{otherwise} \end{array} \right.; \;\;\tau_1^i = \phi^{(i)}.
 \label{eq:basic}\eeq
 
Let $e_i, i=1,\ldots, N$, be the unit $N$-dimensional vectors  with $1$ in the $i$th
position.  Let  $r_k \in \{1,\ldots,N\}$ denote the active source at time $k$. Then one can express the time of arrivals as
the jump Markov linear system
\begin{align*}  x_{k+1} & = A(r_{k+1}) x_k + \snoise_{k}  \\
\obs_k &= C(r_k) x_k + \onoise_k \end{align*}
where 
$$
A(r_{k+1}) = \begin{bmatrix}
I_N & 0_{N\times N} \\
\mbox{diag}(e_{r_{k+1}}) &  I_N
\end{bmatrix},\;
C(r_k) = 
\begin{bmatrix}
 0_{1\times N} & e_{r_k}'
\end{bmatrix}
$$
Note that $r_k $ is a periodic process and so has transition probabilities  $$\tp_{i,i+1} = 1, \; \text{  for $i< N$, and  }\, \tp_{N,1} = 1. $$
$\onoise_k$ denotes the measurement  (jitter) noise; while $\snoise_k$ can be used to model  time varying periods.\\
{\em Remark}: Obviously, there are identifiability issues;  for example, if $\phi^{(1)} = \phi^{(2)}$ and $T^{(1}$ is a multiple of $T^{(2)}$ then
it is impossible to detect source 1.
  
  \item {\bf Narrowband Interference and JMLS.} \index{jump Markov linear system! narrowband interference}  Narrowband interference corrupting a Markov chain can be modeled as a jump Markov linear
  system. Narrowband interference can be modeled as an auto-regressive (AR) process with poles close to the unit circle: for example
$$i_k = a\,  i_{k-1} + w_k $$
where $a = 1- \epsilon$ and $\epsilon $ is a small positive number.
Consider the observation model
$$ y_k = \state_k + i_k + \onoise_k $$
where $\state_k$ is a finite state Markov chain, $i_k$ is narrowband interference and $\onoise_k$ is observation noise. Show that the above
model can be represented as a jump Markov linear system.
  
 \item {\bf Bayesian estimation of Stochastic context free grammar.}  \index{stochastic context free grammars}
 First some perspective: HMMs with a  finite observation space are also called  regular grammars.
 They are a subset of  a more general class of models called stochastic context free grammars as depicted by Chomsky's 
 hierarchy in Figure \ref{fig:chomsky}.
 \begin{figure}[h] \centering
\includegraphics[scale=0.8]{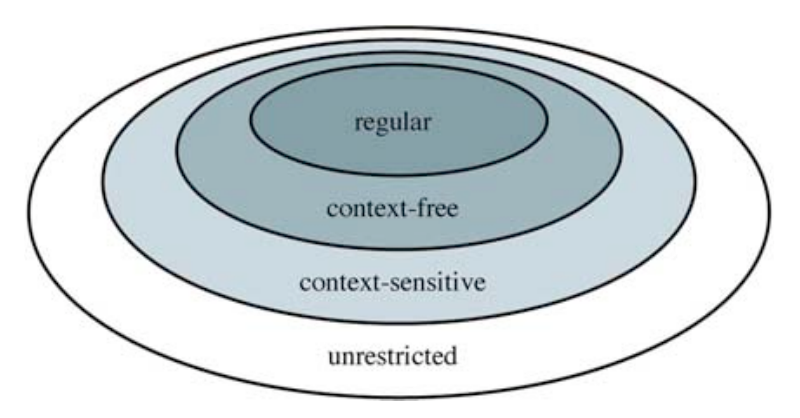}
\caption{The Chomsky hierarchy of languages} \label{fig:chomsky}
\end{figure} 

  Stochastic context free grammars (SCFGs) provide a powerful modeling tool for strings of alphabets and are used widely in
 natural language processing \cite{MS99a}. For example, consider the randomly generated string $a^n c^m b^n$ where $m,n$ are  non-negative  integer valued random variables. Here $a^n$ means the alphabet $a$ repeated $n$ times.
 The string $a^n c^m b^n$
 could model the trajectory of a  target that moves $n$ steps north and then  an arbitrary number of steps east or west and then
 $n$ steps south, implying that the target performs a U-turn.
A basic course in computer science would show (using a pumping lemma) 
that such strings cannot be generated exclusively using a Markov chain (since the memory $n$ is variable).

 If the string $a^n  c^m b^n$ was observed in noise, then Bayesian estimation
(stochastic parsing) algorithms can be used to estimate the underlying string.  Such  meta-level tracking algorithms have polynomial computational
cost (in the data length) and
are useful for estimating {\em trajectories} of targets (given noisy position and velocity measurements). They allow a human radar operator
to interpret tracks and can be viewed as middleware in the human-sensor interface.
Such stochastic context free grammars generalize HMMs and facilitate modeling complex spatial trajectories of targets.

Please refer to \cite{MS99a} for Bayesian signal processing algorithms and EM algorithms for stochastic context free grammars.
\cite{FK13,FK14} gives examples of meta-level target tracking using stochastic context free grammars.

 \item {\bf Kalman vs HMM filter}.  
 A Kalman filter is the optimal state estimator for the linear Gaussian state space model
 \begin{align*}  \state_{k+1} &= \statem \state_k + \snoise_k, \\
\obs_k  &  = \obsm^\p \state_k + \onoise_k  . \end{align*}
where $\snoise$ and $\onoise$ are mutually independent  iid Gaussian processes.

Recall from  (\ref{eq:intro_hmms}),  (\ref{eq:intro_hmmobs}) that for a Markov chain with state space 
$\statespace = \{e_1,\ldots, e_{\statedim}\}$ of unit vectors, an HMM can be expressed as
\begin{align*}  \state_{k+1} &= \tp^\p \state_k + \snoise_k, \\
\obs_k  &  = \obsm^\p \state_k + \onoise_k  . \end{align*}
A key difference  is that in  (\ref{eq:intro_hmms}), $\snoise$ is no longer i.i.d; instead it is 
 a  martingale difference process:
 $\E\{\snoise_k|\state_0,\state_1,\ldots,\state_{k}\} = 0$.
 
From \secn \ref{sec:linearminvar} of the book,  it follows that the Kalman filter is the minimum variance linear estimator for the above HMM.
Of course the optimal {\em linear} estimator (Kalman filter) can perform substantially worse than the optimal estimator (HMM filter).
Compare the performance of the HMM filter and Kalman filter numerically for the above example.

 \item {\bf Interpolation of a HMM.} \index{interpolation of HMM} Consider a   Markov chain $x_k$ with transition matrix $\tp$ where the discrete time clock ticks at intervals of 10 seconds.
 Assume noisy   measurements are obtained of at each time $k$.  Devise a smoothing algorithm  to estimate the state of the Markov chain at 5 second intervals.
 (Note: Obviously on the 5 second time scale, the transition matrix is  $\tp^{1/2}$. For this to be a valid stochastic matrix
 it is sufficient that $\tp$ is a symmetric ultrametric matrix or more generally $\tp^{-1}$ is an M-matrix \cite{HL11}; see also \secn \ref{sec:blackwelldom}
 of the book.)

 \end{compactenum}
 \section{Case Study. Sensitivity of HMM filter to transition matrix} Almost an identical proof to that of geometric ergodicity proof of the HMM filter in \secn \ref{sec:hmmforget} can be used to obtain expressions for the sensitivity of the HMM filter to the HMM parameters.
 \index{HMM filter! sensitivity bound to transition matrix}
 
{\bf  Aim}: We are interested in
a  recursion for $\|\belief_k - \lbelief_k\|_1$ when $\belief_k$ is updated with HMM filter using 
transition matrix $\tp$ and $\lbelief_k$ is updated with HMM filter using transition matrix $\ltp$. That is,
we want an expression
for 
\beq  \|\filter(\belief,\obs;\tp) - \filter(\lbelief,\obs; \ltp) \|_1 \text{  in terms of } \|\belief - \lbelief\|_1.  \label{eq:hmmsensc} \eeq
Such a bound if useful when the HMM filter is implemented with an incorrect transition matrix $\ltp$ instead of actual transition matrix $\tp$.
The idea is that when $\tp $ is close to $\ltp $ then $\filter(\belief,\obs;\tp) $ is close to $\filter(\belief,\obs;\tp) $.

A special case of   (\ref{eq:hmmsensc}) is  to obtain an expression for 
\beq   \|\filter(\belief,\obs;\tp) - \filter(\belief,\obs; \ltp) \|_1  \label{eq:hmmsensc2}\eeq that is when both HMM filters have the same initial belief $\belief$ but are updated
with different transition matrices, namely $\tp$ and $\ltp$.   

 The theorem below obtains expressions for both (\ref{eq:hmmsensc}) and  (\ref{eq:hmmsensc2}).
 
 \begin{theorem-non}  Consider a HMM with transition matrix $\tp$ and state levels $\levels$.
 Let $\epsilon > 0$ denote the  user defined parameter.
 Suppose $\|\ltp-\tp\|_1 \leq \epsilon$, where $\|\cdot\|_1$ denotes the induced 1-norm for matrices.\footnote{
The three  statements  $\| \tp^\p \belief - \ltp^\p \belief \|_1 \leq \epsilon$,  $\|\ltp-\tp\|_1 \leq \epsilon$ and $\sum_{i=1}^\statedim \| (\tp^\p - \ltp^\p)_{:,i} \|_1 \belief(i)  \leq \epsilon$
are all equivalent since $\|\belief\|_1 = 1$.}
Then
\begin{compactenum}
\item The expected absolute deviation between one step of filtering using $\tp$ versus $\ltp$ is upper bounded as:
\begin{multline} \E_{\obs}  \left| \levels^\p \left( \filter(\belief,\obs;\tp) - \filter(\belief,\obs;\ltp) \right) \right| \leq 
\epsilon \sum_\obs \max_{i,j}  \levels^\p (I - \filter(\belief,\obs;\ltp) \one^\p) \oprob_\obs (e_i - e_j)
\label{eq:ebound}
\end{multline}

\item The sample paths of the filtered posteriors and conditional means have the  following explicit bounds at each time $k$:
\begin{multline} \label{eq:samplepath1}
\| \belief_k - \lbelief_k\|_1 \leq \frac{\epsilon}{ \max\{ \statem(\lbelief_{k-1},\obs_k) - \epsilon, \, \mu(\obs_k) \}}    + \frac{ \dob(\ltp) \, \|\belief_{k-1} - \lbelief_{k-1} \|_1}{\statem(\lbelief_{k-1},\obs_k)}
%
\end{multline}
Here $\dob(\ltp)$ denotes the Dobrushin coefficient of the transition matrix $\ltp$ and   $\lbelief_k$ is  the posterior
computed using the HMM filter with  $\ltp$, and
\beq 
\statem(\lbelief,\obs) = \frac{\one^\p \oprob_{\obs} \ltp^\p \lbelief}{\max_i \oprob_{i,\obs}}, \quad
\mu(\obs) = \frac{ \min_i \oprob_{iy} } { \max_i \oprob_{iy} } .
\eeq
\end{compactenum}
\end{theorem-non}

 The above theorem gives explicit upper bounds between the filtered distributions using 
 transition matrices $\ltp$ and $\utp$. The $\E_\obs$ in (\ref{eq:ebound}) is with respect to the 
 measure $\filterd(\belief,\obs;\tp) = \one^\p \oprob_y \tp^\p \belief$  which corresponds to $\prob(\obs_k=\obs| \belief_{k-1} = \belief)$.

\begin{proof}
The triangle inequality for norms  yields
\begin{align}
& \dvar{\belief_{k+1} } { \lbelief_{k+1} }   = \dvar{ \filter(\belief_k,\obs_{k+1};\tp) }{ \filter(\lbelief_k,\obs_{k+1};\ltp) } \nonumber  \\
& \leq  
\dvar{ \filter(\belief_k,\obs_{k+1};\tp) } { \filter(\belief_k,\obs_{k+1};\ltp) } \nn \\ &  \hspace{1cm} + \dvar{ \filter(\belief_k,\obs_{k+1};\ltp) } {  \filter( \lbelief_k,\obs_{k+1}; \ltp)}. 
\label{eq:triangle}
\end{align}

{\bf Part 1}: Consider the first normed term in the right hand side of (\ref{eq:triangle}).
Applying (\ref{eq:rhshmmgeom}) with 
  $\belief = \tp^\p \belief_k$ and $\tbelief = \ltp^\p \belief_k$ yields
\begin{multline*}
  \levels^\p (\filter(\belief_k,\obs;\tp) - \filter(\belief_k,\obs;\ltp)) = \frac{1}{\filterd(\belief,\obs;\tp)}
  \levels^\p \left[ I - \filter(\belief,\obs,\ltp) \one^\p\right] \oprob_\obs (\tp - \ltp)^\p \belief  
\end{multline*}
where $\filterd(\belief,\obs;\tp) = \one^\p  \oprob_\obs \tp^\p \belief$.
  Then Lemma \ref{lem:cs}(i) yields 
\begin{multline*}  \levels^\p (\filter(\belief_k,\obs;\tp) - \filter(\belief_k,\obs;\ltp))\\  \leq 
\max_{i,j}  \frac{1}{\filterd(\belief,\obs;\tp)}  \levels^\p \left[ I - \filter(\belief,\obs,\ltp) \one^\p\right] \oprob_\obs (e_i - e_j) \dvar{\tp^\p \belief}{ \ltp^\p \belief} \end{multline*}
Since $\dvar{\tp^\p \belief}{ \ltp^\p \belief} \leq \epsilon$, taking expectations with respect to the measure $\filterd(\belief,\obs;\tp)$, completes the proof
of the first assertion.
  
{\bf   Part 2}:
Applying Theorem \ref{thm:bayesbackvar}(i) with the notation   $\belief = \tp^\p \belief_k$ and $\tbelief = \ltp^\p \belief_k$ yields
\begin{align}  &\dvar{\filter(\belief_k,\obs;\tp) } { \filter(\belief_k,\obs;\ltp) } \leq 
 \frac{ \max_i \oprob_{i,\obs}   \dvar{\tp^\p \belief_k}{\ltp^\p \belief_k} } { \one^\p  \oprob_{\obs} \ltp^\p \belief_k }  \nonumber\\
& \leq  \frac{\epsilon}{2}\,\frac{ \max_i \oprob_{i,\obs} } { \one^\p  \oprob_{\obs} \ltp^\p \belief_k } 
\leq   \frac{ \max_i \oprob_{i,\obs}  \, \epsilon/2} { \max\{ \one^\p  \oprob_{\obs} \ltp^\p \lbelief_k - \epsilon \max_i \oprob_{i\obs}, \min_i \oprob_{iy}  \}}.
\label{eq:term1}
\end{align}
The second  last inequality follows from the construction of $\ltp$ satisfying  (\ref{eq:con2})  (recall the variational norm is half the $l_1$ norm). The last inequality follows from 
Theorem \ref{thm:bayesbackvar}(ii).

Consider the second normed term in the right hand side of (\ref{eq:triangle}).
Applying Theorem  \ref{thm:bayesbackvar}(i) with notation $\belief = \ltp^\p \belief_k$ and $\tbelief = \ltp^\p \lbelief_k$ yields
\begin{multline} \dvar{ \filter(\belief_k,\obs;\ltp) } {  \filter( \lbelief_k,\obs; \ltp) } 
 \leq 
 \frac{ \max_i \oprob_{i,\obs}   \dvar{\ltp^\p \belief_k}{\ltp^\p \lbelief_k} } { \one^\p  \oprob_{\obs} \ltp^\p \lbelief_k } \\ \leq
\frac{ \max_i \oprob_{i,\obs} \, \dob(\ltp) \, \dvar{ \belief_k}{\lbelief_k} } { \one^\p  \oprob_{\obs} \ltp^\p \lbelief_k } \label{eq:term2}
  \end{multline}
 where the last inequality follows from the submultiplicative property of the Dobrushin coefficient.
Substituting   (\ref{eq:term1}) and  (\ref{eq:term2}) into the right hand side of the triangle inequality (\ref{eq:triangle}) proves the result.
\end{proof}

 \section{Case Study. Reference Probability Method for Filtering}
We describe here the so called {\em reference probability method} for deriving the un-normalized  filtering recursion (\ref{eq:unnormalized}).
The main idea is to start with the joint probability mass function of all observations and states until time $k$, namely, $\pdf(\state_{0:k}, \obs_{1:k})$.
Since this joint density contains all the information we need, it is not surprising that by suitable marginalization and integration, the filtering
recursion and hence the conditional
mean estimate can be computed.

Given the relatively straightforward derivations of the filtering recursions given in Chapter \ref{sec:filtermain} of the book, the reader might wonder why we present
yet another derivation. The reason is that in more complicated filtering problems, the reference probability method gives a systematic way of deriving
filtering expressions.
It is used extensively in \cite{EAM95} to derive filters in both discrete and continuous time.
In continuous time,
the reference probability measure is extremely useful -- it yields the so called Duncan-Mortenson-Zakai equations for nonlinear filtering.

\subsubsection{The Engineering Version}
Suppose the state and observation processes $\{\state_k\}$ and $\{\obs_k\}$  are in a probability space with probability measure $\prob$.
Since the state and observation noise processes are iid, under $\prob$, we have the following factorization:
\begin{align} & \pdf(\state_{0:k}, \obs_{1:k}) =  \prod_{n=1}^k  \pdf(\obs_n|\state_n)\, \pdf(\state_n| \state_{n-1})\,  \belief_0(\state_0) 
 \label{eq:factor}\\ 
 &\propto  \prod_{n=1}^k  \pdf_\onoise\left(\onoisem_n^{-1}(\state_n)\left[\obs_n - \obsm_n(\state_k)\right] \right)
 \pdf_\snoise\left(\snoisem_{n-1}^{-1}(\state_{n-1})\,\left[\state_{n}- \statem_{n-1}(\state_{n-1})\right]\right)  \,  \belief_0(\state_0)   \nonumber
\end{align}
Starting with $\pdf(\state_{0:k}, \obs_{1:k})$,  the 
conditional expectation of any function $\fun(\state_k)$  is
\beq
\E\{ \fun(\state_k) | \obs_{1:k}\} = \frac{ \int \fun(\state_k) \pdf(\state_{0:k}, \obs_{1:k}) d\state_{0:k}}{\int   \pdf(\state_{0:k}, \obs_{1:k}) d\state_{0:k}} 
=  \frac{ \int_\statespace  \fun(\state_k) \left[ \int  \pdf(\state_{0:k}, \obs_{1:k}) d\state_{0:k-1} \right] d\state_k}{\int   \pdf(\state_{0:k}, \obs_{1:k}) d\state_{0:k}} 
\label{eq:erp}
\eeq
The main idea then is to  define the term within the square brackets in the numerator as 
 the un-normalized density 
$ \ubelief_k(x_k) =  \int  \pdf(\state_{0:k}, \obs_{1:k}) d\state_{0:k-1}$.
(Of course then $\ubelief_k(x_k)  = \pdf(\state_k, \obs_{1:k})$).
We now derive the recursion  (\ref{eq:unnormalized})  for the un-normalized density $\ubelief_k$:
\begin{align*}
\int_\statespace \fun(\state_k) \ubelief_k(\state_k) d\state_k&= \int_\statespace \fun(\state_k)   \int  \pdf(\state_{0:k}, \obs_{1:k}) d\state_{0:k-1} d\state_k \\
 &\hspace{-1.5cm}= \int_\statespace   \int_\statespace  \fun(\state_k) \pdf(\obs_k| \state_k) \pdf(\state_k| \state_{k-1})  \left[\int \pdf(\state_{0:k-1},\obs_{1:k-1}) d\state_{0:k-2} \right]\, d\state_{k-1}   d\state_k
 \\
 &\hspace{-1.5cm}= \int_\statespace  \int_\statespace   \fun(\state_k)  \pdf(\obs_k| \state_k) \pdf(\state_k| \state_{k-1})  \ubelief_{k-1}(\state_{k-1}) d\state_{k-1} d\state_k
 \end{align*}
where the second equality follows from (\ref{eq:factor}). 
Since the above holds for any test function $\fun$, it follows that the integrands within the outside integral are equal, thereby yielding the
un-normalized  filtering
recursion  (\ref{eq:unnormalized}).

\subsubsection{Interpretation as Change of Measure}
We now interpret the above derivation as the engineering version of the reference probability method.\footnote{In continuous time, the change of measure of a random process involves Girsanov's theorem, see \cite{EAM95}.  Indeed the Zakai form of the continuous time filters  in the appendix of the book can be derived in a fairly
straightforward manner using Girsanov's theorem.}
Define a new probability measure $\bar{\prob}$ as having associated density
$$\pref(\state_{0:k}, \obs_{1:k}) = \prod_{n=1}^k  \pdf_\onoise(\obs_n) \,
 \pdf_\snoise(\state_n)  \,  \belief_0(\state_0)  .  $$
 The above equation is tantamount to saying that under this new measure $ \bar{\prob}$, the processes $\{\state_k\}$ and $\{\obs_k\}$ are iid sequences
with density functions $\pdf_\snoise$ and $\pdf_\onoise$, respectively. $\bar{\prob}$ will be called the {\em reference probability measure} - under this measure, due to the iid nature of  $\{\state_k\}$ and $\{\obs_k\}$, the filtering
recursion can be derived conveniently, as we now describe.

Let
 $\bE$ denote expectation associated with measure $\bar{P}$, so that for any function $\fun(\state_k)$, the conditional expectation is
 $$ \bE\{ \fun(\state_k) | \obs_{1:k}\} = \int \fun(\state_k) \pref(\state_{0:k}, \obs_{1:k})  d\state_{0:k} $$
 Obviously, to obtain the expectation $\E\{\fun(\state_k)| \obs_{1:k}\}$ under the probability
measure $\prob$, it follows from (\ref{eq:erp}) that
 \begin{align}
  \E\{ \fun(\state_k) | \obs_{1:k}\} &= \frac{\int \fun(\state_k)  \Lambda_k   \pref(\state_{0:k}, \obs_{1:k})  d\state_{0:k}}
 {\int \Lambda_k  \pref(\state_{0:k}, \obs_{1:k})  d\state_{0:k}}  ,
  \quad \text{ where }   \Lambda_k = \frac{\pdf(\state_{0:k}, \obs_{1:k}) }{ \pref(\state_{0:k}, \obs_{1:k}) }  \\
  &= \frac{\bE\{  \Lambda_k\fun(\state_k) | \obs_{1:k} \} }{ \bE\{ \Lambda_k| \obs_{1:k}\}}
  \nonumber
  \end{align}

The derivation then proceeds as follows. \begin{align*}
& \int_\statespace  \ubelief_k(\state) \fun(\state) d\state  =
\bE\{\Lambda_k \fun(\state_k) | \mathcal{Y}_k\} \qquad  \text{ (definition of $\ubelief_k$) } \\ 
&= \int \frac{\pdf(\state_{0:k},\obs_{1:k})}{\pref(\state_{0:k},\obs_{1:k})} \fun(\state_k)  \pref(\state_{0:k},\obs_{1:k}) d\state_{0:k} \\
&=  \int \frac{\pdf(\state_{0:k-1},\obs_{1:k-1})}{\pref(\state_{0:k-1},\obs_{1:k-1})} 
\frac{\pdf(\obs_k | \state_k) \pdf(\state_k| \state_{k-1})}{\cancel{\pdf_\onoise(\obs_k) } \cancel{ \pdf_\snoise(\state_k)} }
\fun(\state_k)  \cancel{\pdf_\onoise(\obs_k)}  \cancel{\pdf_\snoise(\state_k) } \pref(\state_{0:k-1},\obs_{1:k-1})   d\state_{0:k} \\
&=  \int \frac{\pdf(\state_{0:k-1},\obs_{1:k-1})}{\pref(\state_{0:k-1},\obs_{1:k-1})} 
\left[\int_\statespace {\pdf(\obs_k | \state_k) \pdf(\state_k| \state_{k-1})} \fun(\state_k) d\state_k\right]  \pref(\state_{0:k-1},\obs_{1:k-1})   d\state_{0:k-1} 
\\
&= \bE\{ \Lambda_{k-1} \left[\int_\statespace {\pdf(\obs_k | \state_k) \pdf(\state_k| \state_{k-1})} \fun(\state_k) d\state_k\right]  | \obs_{1:k-1} \} \\
&= \int \ubelief_{k-1}(\state_{k-1})  \left[\int_\statespace {\pdf(\obs_k | \state_k) \pdf(\state_k| \state_{k-1})} \fun(\state_k) d\state_k\right]  d\state_{k-1}
\end{align*}
where the last equality follows from the definition of $\ubelief$ in the first equality.

Since this holds for any test function $\phi(\state)$, we have that the material inside the integral in the left and right hand side are equal.
So 
$$ \belief_k(\state_k )  =  \pdf(\obs_k| \state_k)   \int_\statespace \ubelief_{k-1}(\state_{k-1}) \pdf(\state_k| \state_{k-1}) d\state_{k-1}. $$


\chapter{Algorithms for Maximum Likelihood Parameter Estimation}   

\begin{compactenum}

\item A standard drill exercise involves deriving  the Cram\'er-Rao bound in terms of the Fisher information matrix; see  wikipedia or any book
in statistical signal processing  for an elementary
description.

\item {\bf Minorization Maximization Algorithm (MM Algorithm).} The EM algorithm is a special case of the  MM algorithm\footnote{MM can also be used equivalently to denote majorization minimization}; see \cite{HL04} for a nice tutorial on MM 
\index{Minorization Maximization algorithm} 
algorithms.  MM algorithms constitute a general purpose method for optimization and are not restricted just to maximum likelihood estimation.

The main idea behind the MM algorithm  is as follows: Suppose we wish to compute the maximizer $\th^*$ of a function $\fun(\model)$.
The idea is to construct a minorizing  function $g(\model, \model^{(m)})$  such that
\beq \begin{split}
g(\model, \model^{(m)})  & \leq \fun(\model)  \quad \text{ for all } \model \\
  g(\model^{(m)} , \model^{(m)}  ) &=  \fun(\model^{(m)}) . \label{eq:minorf}
\end{split} \eeq
That is, the minorizing function $g(\model, \model^{(m)}) $  lies above $ \fun(\model)$ and is a tangent to it at the point $\model^{(m)}$.
Here  $$\model^{(m)}  = \argmax_\th g(\th^{(m-1)}, \th) $$
 denotes the estimate of the maximizer at iteration $m$ of MM algorithm.

The property (\ref{eq:minorf})
implies that successive iterations of the MM algorithm yield $$\fun(\model^{(m+1)}) \geq \fun(\model^{(m)}). $$
In words, successive iterations of the MM algorithm yield increasing values of the objective function which is a very useful property for a general
purpose numerical optimization algorithm.
This is shown straightforwardly as follows:
\begin{align*}  \fun(\th^{(m+1}) &=  \fun(\th^{(m+1}) - g(\model^{(m+1)}, \model^{(m)}) + g(\model^{(m+1)}, \model^{(m)}) \\
& \stackrel{a}{\geq}  \fun(\th^{(m+1}) - g(\model^{(m+1)}, \model^{(m)})  + g(\model^{(m)}, \model^{(m)})  \\
& \stackrel{b}{\geq} \fun(\th^{(m}) - \cancel{g(\model^{(m)}, \model^{(m)}) } + \cancel{g(\model^{(m)}, \model^{(m)})}
 \end{align*}
Inequality (a) follows since $g(\model^{(m+1)}, \model^{(m)}) \geq g(\model^{(m)}, \model^{(m)})$ by definition since $\model^{(m+1)} =
\argmax_\th g(\model, \model^{(m)}) $. Inequality (b) follows from (\ref{eq:minorf}).

The EM algorithm is a special case of the MM algorithm where  \index{EM algorithm} \index{Minorization Maximization algorithm! EM algorithm}
$$ g(\th,\model^{(m)} = Q(\th,\model^{(m)}) - Q(\model^{(m)},\model^{(m)}), \quad
\fun(\th) = \logl(\th) -  \logl(\th^{(m)} ) $$
Here $\logl(\model) = \log \pdf(\obs_{1:N} | \model)$ is the log likelihood which we want to maximize to compute the MLE
and $Q(\th,\model^{(m)}) $ is the auxiliary log likelihood defined
in (\ref{eq:auxl0}) which is maximized in the M step of the  EM algorithm.

Indeed the minorization property  (\ref{eq:minorf}) was established for the EM algorithm in  Lemma \vref{lem:em} of the book by using Jensen's inequality.

\item {\bf EM algorithm in more elegant (abstract) notation.}
Let $\{P_\theta\,,\,\theta\in \Theta\}$ be a family of probability
measures on a measurable space $(\Omega, {\cal F})$ all absolutely
continuous with respect to a fixed probability measure $P_0$,
and let ${\cal Y} \subset {\cal F}$. The likelihood function for
computing an estimate of the parameter $\theta$ based on the
information available in ${\cal Y}$ is
\begin{displaymath}                                                                                                                                              
   L(\theta) = \E_0[ \frac{dP_\theta}{dP_0}\mid {\cal Y}]\ ,                                                                                                     
\end{displaymath}
and the MLE estimate is defined by
\begin{displaymath}                                                                                                                                              
   \widehat{\theta}\in \argmax_{\theta\in \Theta} L(\theta)\ .                                                                                                   
\end{displaymath}
In general, the MLE is difficult to compute directly,
and the EM algorithm provides an iterative approximation method~:
\begin{description}
   \item{{\it Step 1.}} Set $p=0$ and choose $\widehat{\theta}_0$.

   \item{{\it Step 2.}} (E--step) Set $\theta' = \widehat{\theta}_p$
and compute $Q(\cdot,\theta')$, where
\begin{displaymath}                                                                                                                                              
   Q(\theta,\theta')                                                                                                                                             
   = \E_{\theta'}[ \log\frac{dP_\theta}{dP_{\theta'}}\mid {\cal Y}]\ .                                                                                           
\end{displaymath}

   \item{{\it Step 3.}} (M--step) Find
\begin{displaymath}                                                                                                                                              
   \widehat{\theta}_{p+1}\in \argmax_{\theta \in \Theta} Q(\theta,\theta')\ .                                                                                    
\end{displaymath}

   \item{{\it Step 4.}} Replace $p$ by $p+1$ and repeat beginning
with Step~2, until a stopping criterion is  satisfied.
\end{description}
The sequence generated $\{\widehat{\theta}_p\,,\,p\geq 0\}$
gives non--decreasing values of the likelihood function~:
indeed, it follows from Jensen's inequality that
\begin{displaymath}                                                                                                                                              
   \log L(\widehat{\theta}_{p+1}) - \log L(\widehat{\theta}_p)                                                                                                   
   \geq Q(\widehat{\theta}_{p+1},\widehat{\theta}_p)                                                                                                             
   \geq Q(\widehat{\theta}_p,\widehat{\theta}_p) = 0\ ,                                                                                                          
\end{displaymath}
with equality if and only if $\widehat{\theta}_{p+1} = \widehat{\theta}_p$.

\item  {\bf Forward-only EM algorithm for Linear Gaussian Model.}  In \secn 4.4 of the book, we described  a forward-only EM algorithm for ML parameter  estimation of the a HMM.
Forward-only EM algorithms can also be constructed for maximum likelihood estimation of the parameters of a  linear Gaussian state space model
\cite{EK99}. These involve computing filters for functionals of the state and use Kalman filter estimates.

\item  {\bf Sinusoid  in HMM.} Consider a sinusoid with amplitude $A$ and phase $\phi$. It is observed  as
$$y_k = \state_k + A\, \sin(k/100 + \phi) + \onoise_k$$
where $v_k$ is an iid Gaussian noise process.
Use the EM algorithm to estimate $A, \phi$ and the parameters of the Markov chain and noise variance.

\item In the forward-only  EM algorithm of  \secn 4.4, the filters for the number of jumps involves $O(\statedim^4)$ computations at each time while filters for the duration time involve
$O(\statedim^3)$ at each time. Is it possible to reduce the computational cost by approximating some of these estimates? 

\item  Using computer simulations, compare the methods of  moments estimator for a HMM in \secn 4.5  with the maximum likelihood estimator in terms of efficiency. That is generate several $N$ point trajectories of an HMM with a fixed set of parameters, then compute the variance
of the estimates. (Of course,  instances where the MLE
the algorithm converges to local maxima should be eliminated from the computation).

\item Non-asymptotic statistical inference using concentration of measure if very popular today.
Assuming the likelihood is a Lipschitz function of the observations, and the observations are Markovian, show that
the likelihood function concentrates to the Kullback Leibler function.

\item {\bf EM Algorithm for State Estimation.} The EM algorithm was used in Chapter \ref{chp:mle} as a numerical algorithm
for maximum likelihood {\em parameter} estimation.  It turns out that the EM algorithm can be used for  {\em state} estimation, particularly for
a jump Markov linear system (JMLS). Recall from \secn \ref{sec:particle}  that a JMLS has model
\begin{align*}
z_{k+1} &= \statem(\mc_{k+1})\,z_k + \snoisem(\mc_{k+1})\,\snoise_{k+1} + \inpm (\mc_{k+1})\,\inp_{k+1} \\
\obs_k &= \obsm(\mc_k)\,z_k + \onoisem(\mc_k)\,\onoise_k + \oinpm(\mc_k)\inp_k. 
\end{align*}
As described in  \secn \ref{sec:particle}, 
 the optimal filter for a JMLS is  computationally intractable.
  In comparison for a JMLS, the EM algorithm can be used to estimate the MAP (maximum aposteriori state estimate).
system (assuming the parameters of the JMLS are known).  Show how 
one can compute this MAP state estimate $\max_{z_{1:k},r_{1:k}} P(y_{1:k} | z_{1:k},r_{1:k})$ using the EM algorithm.
In \cite{LK99} is shown that the resulting EM algorithm involves the cross coupling of a Kalman and HMM smoother.
A data augmentation algorithm in similar spirit appears in \cite{DLK00}. \index{jump Markov linear system! EM algorithm for state estimation}

\item {\bf Quadratic Convergence of Newton Algorithm.}

\index{Newton algorithm! quadratic convergence}
We start with some  definitions:
Given a sequence $\{\th^{(n)}\}$ generated by an optimization algorithm, the {\em order}  of convergence is $p$ if 
\beq  \beta =  \limsup_{n \rightarrow \infty} \frac{\|\th^{(n+1) } - \th^*\|} {\|\th^{(n) } - \th^*\|^p}   \text{ exists } 
\label{eq:convergenceorder}
\eeq
Also if $p=1$ and $\beta < 1$, the sequence is said to converge linearly to $\th^*$ with {\em convergence ratio (rate)} $\beta$.
Moreover, the case  $p=1$ and $\beta = 0$ is referred to as superlinear convergence.

\begin{compactenum}
\item
Recall that the Newton Raphson algorithm computes the MLE iteratively as
$$ \th^{(n+1)} = \th^{(n)} +   \bigl({\nabla^2} \logl (\th^{(n)}) \big)^{-1} \nabla \logl(\th^{(n)})  $$
 The Newton Raphson algorithm has quadratic order of convergence in the following sense.
Suppose  the log   likelihood $\logl(\th)$ is twice continuous differentiable and that at a local maximum $\th^*$, the Hessian
$\nabla_\th^2\logl$ is positive definite. Then if started sufficient close to $\th^*$,  Newton Raphson converges to $\th^*$ at a quadratic rate.
that the model estimates satisfy $\th^{(n)}$ satisfy
$$\| \th^{(n+1)} - \th^* \|  \leq \beta \| \th^{(n)} - \th^*\|^2 $$
for some constant $\beta$.

This is shown straightforwardly (see any optimization textbook) as follows:
\beq
\begin{split}  \|\th^{(n+1)} - \th^* \|  & = \|   \th^{(n)} - \th^* + \bigl({\nabla^2} \logl (\th^{(n)}) \big)^{-1} \nabla \logl(\th^{(n)}) \|
\\ 
&= \|  \bigl({\nabla^2} \logl (\th^{(n)}) \big)^{-1}  \biggl(  \nabla \logl(\th^{(n)}) - \nabla \logl(\th^*)    - \nabla^2 \logl (\th^{(n)}) \big(\th^{(n) }- \th^*  \big)
\biggr)
\end{split}
\eeq
For $\|\th^{(n)} - \th^*\| < \rho$, it is clear from a Taylor series expansion that
$$  \|\nabla \logl(\th^*)  -  \nabla \logl(\th^{(n)})  - \nabla^2 \logl (\th^{(n)}) \big(\th^* - \th^{(n)} \big)  \| \leq  \beta_1\| \th^{(n)} - \th^*\|^2$$
for some positive constant $\beta_1$.
Also,  $\|  \bigl({\nabla^2} \logl (\th^{(n)}) \big)^{-1}  \| \leq \beta_2$.

\item  The convergence order and rate
of the EM algorithm has been studied in great detail since the early 1980s; there are numerous papers in the area; see \cite{XJ96} and the
references therein.
The EM algorithm has linear convergence order, i.e., $p=1$ in (\ref{eq:convergenceorder}). 
Please see
\cite{MXJ00} and the references therein for examples where EM exhibits superlinear convergence.

\end{compactenum}
\end{compactenum}

\chapter{Multi-agent Sensing: Social Learning and Data Incest}  

\section{Problems}
\begin{enumerate}

\item A substantial amount of insight can be gleaned by actually simulating the setup (in Matlab) of the social learning
filter for both the random variable and Markov chain case. Also simulate the risk-averse social learning filter discussed in  \secn \ref{sec:classicalsocial}
of the book.


\item {\bf CVaR Social Learning Filter.} Consider the risk averse social learning  \index{CVaR social learning filter}
discussed in    \secn \ref{sec:classicalsocial}.
Suppose
agents  choose their  actions $a_{k}$ to minimize the CVaR   risk averse measure 
$$
a_{k} =   {\underset{a \in \mathcal{A}}{\text{argmin}}} \{ {\underset{z \in \mathbb{R}}{\text{min}}} ~ \{ z + \frac{1}{\alpha} \mathbb{E}_{y_{k}}[{\max} \{ (c(x_{k},a)-z),0 \rbrace] \} \} 
$$
Here $\alpha \in (0,1]$ reflects the degree of risk-aversion for the agent (the smaller $\alpha$ is, the more risk-averse the agent is). 
Show that the structural result Theorem \ref{thm:monotone} continues to hold for the CVaR social learning filter.
Also show that for sufficiently risk-averse agents (namely, $\alpha$ close to zero), social learning ceases and agents always herd.

Generalize the above result to any coherent risk measure.

\item The necessary and sufficient condition given in Theorem \ref{thm:sufficient} for exact data incest removal requires
 that 
$$ A_n(j,n)=0   \implies w_n(j)= 0,  \text{ where }  \quad w_n =  T_{n-1}^{-1}  t_n,
$$
and $T_n =  \text{sgn}((\mathbf{I}_n-A_n)^{-1}) = \begin{bmatrix} T_{n-1} &  t_n \\ 0_{1 \times n-1} & 1 \end{bmatrix} $ is the transitive closure matrix.
Thus the condition depends purely on the adjacency matrix.  Discuss what types of matrices satisfy the above condition.

\item  Theorem \ref{thm:sufficient}  also applies to data incest where the prior and likelihood are Gaussian. The posterior is then
evaluated by a Kalman filter. Compare the performance of exact data incest removal with the covariance intersection algorithm
in \cite{CAM02} which assumes no knowledge of the correlation structure (and hence of the network).

\item Consensus algorithms \cite{SM04} have been  extremely popular during the last decade and there are numerous papers in the area.
 They are non-Bayesian and seek to compute,
for example, the average over measurements observed at a number of nodes in a graph.  It is worthwhile comparing the performance
of the optimal Bayesian incest removal algorithms with consensus algorithms.

\item The data incest removal algorithm in \secn \ref{sec:incest} of the book arises assumes that agents do not send additional information
apart from their incest free estimates. Suppose agents are allowed to send a fixed number of  labels of previous agents from whom
they have received information. What is the minimum about of additional labels the agents need to send in order to completely remove data incest.

\item Quantify the bias introduced by data incest as a function of the adjacency matrix.

\item Prospect theory \index{prospect theory} (pioneered by the psychologist Kahneman \cite{KT79} who won the 2003 Nobel prize in economics) is a behavioral economic theory
that seek to model how humans make decisions amongst probabilistic alternatives. (It is an alternative to expected utility theory considered
in the social learning models of this chapter.)  The main features are:
\begin{compactenum}
\item Preference is an S-shaped curve with reference point $x=0$
\item The investor maximizes the expected  value $V(x)$ where $V$ is a preference and $x$ is the change in wealth.

\item Decision maker employ decision weight $w(p)$ rather than objective probability $p$, where the weight function $w(\cdf)$ has a reverse
S shape where $\cdf$ is the cumulative probability.
\end{compactenum}
Construct a social learning filter where the utility function satisfies the above assumptions. Under what conditions do information cascades
occur?

\item {\bf Rational Inattention.} Another powerful way  for modeling the behavior of (human) decision makers is in terms of rational inattention. 
See the seminal work of \cite{Sim03} where essentially the ability of the human to absorb information is modeled via the 
information theoretic capacity of a communication channel.\index{rational inattention}

\item There are several real life experiments that seek to understand how humans interact in decision making.  See for example \cite{BOL10}
and  \cite{KH15}. In \cite{BOL10}, four  models are considered. How can these models be linked to social learning?
\end{enumerate}

\section{Social Learning with limited memory} \index{social learning! limited memory}
Here we briefly describe a variation of the  vanilla social learning protocol.  In order to mitigate herding, assume that agents randomly
sample only a fixed number of previous actions. The aim below is to describe the resulting setup; see \cite{SS97} for a detailed discussion.

Let the variable $\st \in \{1,2\}$ denote the states. Let $a \in \{1,2\} $ denote the action alphabet and $y \in \{1,2\} $ denote the observation alphabet. In this model of social learning with limited memory, it is assumed that each agent (at time $t \geq N+1$) observes \textit{only} $N$ randomly selected actions from the history $h_{t}= \{a_{1},a_{2},\hdots,a_{t-1} \}$. In the periods $t \leq N$, each agent acts according to his private belief. This phase is termed as the seed phase in the model. \\

Let $\zt$ denote the number of times action $1$ is chosen until time $t$, i.e, $$\zt = \sum_{j=1}^{t} I(a_{j}=1).$$ Let $\zht$ denote the number of times action $1$ is chosen in a sample of $N$ randomly observed actions in the past, i.e, $\zht = \sum_{j=1}^{N} I(a_{j}=1)$. \\

The social learning protocol with limited memory is as follows:
\begin{itemize}
\item[1.)] \textit{Private belief update}: Agent $t$ makes two observations at each instant $t(>N)$. These observations correspond to a noisy private signal $y_{t}$ and a sample of $N$ past actions from the history $h_{t}$ sampled uniformly randomly.  Let $B_{y_{t}}$ and $D_{z_{t}=k}$ denote the probability of observing $y_{t}$ and $(\zht=k)$ respectively. The private belief is updated as follows.

For each draw from the past, the probability of observing action 1 is $\zt/ (t-1)$. So the
probability that at time $t$, action $1$ occurs $k$ times in a random sample of $N$ observed actions  is
\begin{equation*}
\mathbb{P}(\zht = k|\zt) = \frac{N!}{(N-k)!k!} \left( \frac{\zt}{t}\right) ^{k} \left(1-\frac{\zt}{t}\right)^{N-k}
\end{equation*}
Therefore,
the number of times action `$1$' is chosen in the sample, $\zht$, has a distribution that depends on $\st$ according to:
\begin{equation*}
\mathbb{P}(\zht = k|\st) = \sum_{\zt=1}^{t} \mathbb{P}(\zht = k|\zt) \mathbb{P}(\zt|\st)
\end{equation*}

After obtaining a private noisy signal $y_{t}$, and having observed ($\zht=k$), the belief  $\belief_t = [\belief_t(1), \belief_t(2)]^\p$ where
$\belief_t(i) = \prob(\th=i| \zht,y_t)$
is updated by agent $t$ as:
\begin{equation*}
{\hspace{1.5cm}}\pi_{t} = \frac{B_{y_{t}}D_{z_{t}=k}}{\textbf{1}'B_{y_{t}}D_{z_{t}=k}}.
\end{equation*}
Here $B$ and $D$ are the observation likelihoods of $y_t$ and $\zht$ given the state:
\[   ~~B_{y_{t}}=\text{diag}(\mathbb{P}(y_{t}|\st=i),i\in \{1,2\}), \quad 
D_{z_{t}=k} = \begin{bmatrix}
\mathbb{P}(\zht=k|\st=1) \\
\mathbb{P}(\zht=k|\st=2)
\end{bmatrix}
\]

\item[2.)] \textit{Agent's decision}: With the private belief $\pi_{t}$, the agent $t$ makes a decision as:
\begin{equation*}
a_{t} = {\underset{a\in \{1,2\}}{\text{argmin}}}~ c_{a}^{T} \pi_{t}
\end{equation*}
where $c_{a}$ denotes the cost vector.

\item[3.)] \textit{Action distribution}: The distribution of actions $\mathbb{P}(\zt|\st)$ in the two states $\st = 1,2$ is assumed to be common knowledge at time $t$. It is updated after the decision of agent $t$ as follows. \\

The probability of $(a_{t}=1)$ in period $t$ depends on the actual number of `$1$' actions $\zt$ and on the state according to:
\begin{equation*}
{\hspace{-0.75cm}}\mathbb{P}(a_{t}=1|\zt=n,\st) = \sum_{k=0}^{N} \sum_{i=1}^{2} \mathbb{P}(a_{t}=1|y=i,\zht=k,\zt=n,\st)~ \mathbb{P}(y=i|\st)~\mathbb{P}(\zht = k|\zt=n)
\end{equation*}

where,
\[ \mathbb{P}(a_{t}=1|y=i,\zht=k,\zt=n,\st) = \left\{ \begin{array}{ll}
         1 & \mbox{if $c_{1}^{T}B_{y=i}D_{z_{t}=k} < c_{2}^{T}B_{y=i}D_{z_{t}=k}$};\\
        0 & \mbox{otherwise}.\end{array} \right. \]

After agent $t$ takes an action, the distribution is updated as:
\begin{equation}{\label{eq:dis}}
{\hspace{-0.75cm}}\mathbb{P}(\ztp=n|\st) = \mathbb{P}(\zt=n|\st)(1-\mathbb{P}(a_{t}=1|\zt=n,\st)) + \mathbb{P}(\zt=n-1|\st)\mathbb{P}(a_{t}=1|\zt=n,\st)
\end{equation}
\end{itemize}

According to equation~\eqref{eq:dis}, the sufficient statistic $\mathbb{P}(\zt|\st)$ is growing with time $t$. It is noted that this has $(t-2)$ numbers at time $t$ and hence grows with time. $\mathbb{P}(\ztp=n|\st)$ in equation~\eqref{eq:dis} is used to compute $D_{z_{t+1}}$.

With the above model, consider the following questions:
\begin{compactenum}
\item Show that there is asymptotic herding when $N =1$.
\item  Show that for $N=2A$, reduction in the historical information will improve social learning. Also, comment on whether there is herding when $N=2$.
\item  Show that as $N$  increases, the convergence to the true state is slower.
Hint: Even though more observations are chosen, greater weight on the history precludes the use of private information.
\end{compactenum}

\chapter{Fully Observed Markov Decision Processes}  

\section{Problems}
\begin{compactenum}

\item  The following nice example from \cite{KV86} gives a useful motivation for feedback control in stochastic systems. It shows that for
stochastic systems, using feedback control can result in behavior that cannot be obtained by an open loop system.
\index{why feedback control?}
\begin{compactenum}
\item First, recall from undergraduate control courses that  for a deterministic linear time invariant system with forward transfer function $G(z^{-1})$ and negative feedback $H(z^{-1})$, the equivalent transfer function
is $\frac{G(z^{-1})}{1+G(z^{-1}) H(z^{-1}) }$. So an open loop system with this equivalent transfer function is identical to a feedback system.

\item More generally, consider the deterministic system $$\state_{k+1} = \phi(\state_k, \action_k),  \obs_k = \psi(\state_k,\action_k) $$ Suppose the actions are given
by a  policy of the form $$u_k = \policy(\state_{0:k},\obs_{1:k})$$
Then clearly, the open loop system, $$\state_{k+1} = \phi(\state_k, \policy(\state_{0:k},\obs_{1:k})),  \obs_k = \psi(\state_k,\policy(\state_{0:k},\obs_{1:k})) $$ generates the same
state and observation sequences. 

So for a deterministic system (with fully specified model), open and closed loop behavior are identical.

\item Now consider a fully observed stochastic system with feedback:
\beq
\begin{split}  \state_{k+1} &= \state_k + \action_k + \snoise_k , \\
 u_k &= - \state_k  \end{split} \eeq
 where $\snoise_k$ is iid with zero mean and variance $\sigma^2$ (as usual we assume $x_0$ is independent of $\{\snoise_k\}$.)
 Then $\state_{k+1} = \snoise_k$ and so $u_k = -\snoise_{k-1}$ for $k=1,2,\dots$.
 Therefore $\E\{x_k\} = 0$ and $\var\{x_k^2\} = \sigma^2$.
 
\item  Finally, consider an open loop stochastic system where $u_k$ is a deterministic sequence:
$$
  \state_{k+1}= \state_k + \action_k + \snoise_k  $$
 Then $\E\{x_k\} = \E\{\state_0\} + \sum_{n=0}^{k-1} u_k $ and $\var\{x_k^2\} = \E\{x_0^2\} + k \sigma^2$.
 Clearly, it is impossible to construct a deterministic input sequence that yields a zero mean state with variance $\sigma^2$.
 \end{compactenum}

\item {\bf Trading of call options.} An investor buys a call option at a price $p$. He has $\finaltime$ days to exercise this option. If the investor exercises the option when the stock price is $x$, he
gets $x-p$  dollars.  The investor can also decide not the exercise the option at all.

Assume the stock price evolves as $x_k = x_0 + \sum_{n=1}^k w_n $ where $\{w_n \}$ is in iid process.
Let $\tau$ denote the day the investor decides to exercise the option.  Determine the optimal investment strategy to maximize
$$\E\{ (x_\tau - p) I(\tau \leq T) \}. $$
This is an example of a fully observed stopping time problem. Chapter \ref{ch:pomdpstop} considers more general stopping time POMDPs.

Note: Define $s_k \in \{0,1\}$ where $s_k = 0$ means that the option has not been exercised until time $k$.
$s_k = 1$ means that the option has been exercised before time $k$. Define the state $z_k = (x_k,s_k)$.

Denote the action $u_k = 1$  to exercise option and $u_k = 0$ means do not exercise option.
Then the dynamics are
$$ s_{k+1} = \max\{ s_k,u_k\} ,  \quad  x_{k+1} = x_k + w_k $$
The reward at each time $k$ is $r(z_k,u_k,k) = (1-s_k) u_k (x_k-p)$ and the problem can be formulated as
$$\max_\policy  \E\{ \sum_{k=1}^\finaltime r(z_k,u_k,k)\} $$

\item 
Discounted cost problems  can also be motivated as stopping time problems (with a random termination time).
Suppose at each time $k$, the MDP can terminate with probability $1 - \discount$ or continue with probability $\discount$.
Let $\tau$ denote the random variable for the termination time. Consider the undiscounted cost MDP
\begin{align*}
 \E_{\bpolicy}  \left\{  \sum_{k=0}^{\tau} \cost(\state_k,\action_k)  \mid \state_0 = i\right\} &= 
\E_{\bpolicy}  \left\{  \sum_{k=0}^{\infty} I(k \leq \tau) \, \cost(\state_k,\action_k)  \mid \state_0 = i\right\} \\
&= \E_{\bpolicy}  \left\{  \sum_{k=0}^{\infty} \discount^k  \cost(\state_k,\action_k)  \mid \state_0 = i\right\}.
\end{align*}
The last equality follows since $\prob(k \leq \tau) = \discount^k$.

\item We discussed risk averse utilities and dynamic risk measures briefly in \secn \ref{sec:riskaverse}. Also \secn \ref{subsec:afriat} discussed revealed
preferences for constructing a utility function from a dataset.  
Given a utility  function  $U(x)$,
a widely used measure for the degree of risk aversion is the Arrow-Pratt risk aversion coefficient which is defined as 
$$ a(x) =- \frac{d^2 U/dx^2}{ dU/dx} .$$
This is often termed as an absolute risk aversion measure, while $x\, a(x)$ is termed a relative risk aversion measure. 
Can  this risk averse coefficient be used for mean semi-deviation risk, conditional value at risk( CVaR) and exponential risk?

\item A classical result involving utility functions is the following \cite[pp.42]{HS84}:
A rational decision maker who compares random variables only according to their means and variances must have preferences consistent
with a quadratic utility function. Prove this result.

\end{compactenum}

\section{Case study.  Non-cooperative Discounted Cost Markov games}  \index{Markov game} \label{sec:markovgame}
\secn \ref{sec:mdpdiscount} of the book dealt with infinite horizon discounted MDPs.
Below we introduce briefly some elementary ideas in non-cooperative infinite horizon discounted Markov games. There are several excellent 
books in the area
\cite{KV12,BO91}. 

Markov games can be viewed as a multi-agent  decentralized extension of MDPs. They arise in a variety of applications including  dynamic spectrum allocation, financial models and smart grids. Our aim here is to consider  some simple cases where the Nash
equilibrium can be obtained by solving a linear programming problem.\footnote{The reader should be cautious with decentralized stochastic
control. The famous Witsenhausen's counterexample formulated in the 1960s shows that even a deceptively simple toy problem in decentralized
stochastic control can be very difficult to solve,  see \url{https://en.wikipedia.org/wiki/Witsenhausen\%27s_counterexample} \index{Witsenhausen's counterexample}}

Consider the following infinite horizon discounted cost two-payer Markovian game.
There are two decision makers (players) indexed by $l=1,2$.
\begin{compactitem}
\item Let $\act1_k \in \actionspace $ and $\act2_k\in \actionspace $  denote the action of player 1 and player 2, respectively, at time $k$. For convenience we assume
the same action space for both players.
\item The cost incurred
by player $l \in \{1,2\}$  for state $\state$, actions $\act1,\act2$ is $c_l(\state,\act1,\act2)$.
\item The transition probabilities  of the Markov process $x$ depends on the  actions of both players:
$$ \tp_{ij}(\act1,\act2) =   \prob(\state_{k+1} = j | \state_k = i, \act1_k = \act1, \act2_k = \act2) $$
\item 
Define the policies for the stationary (randomized) Markovian policies for two players as  $\pol1$, $\pol2$, respectively.
So $\act1_k $ is chosen from probability distribution $\pol1(\state_k)$ and $\act2_k $ is chosen from probability distribution $\pol2(\state_k)$.
For convenience denote the class of stationary Markovian policies as $\Policy_S$.
\item The cumulative cost incurred by each player $l \in \{1,2\}$ is 
\beq  \Jc{l}_{\pol1,\pol2}(x) =  \E\big\{  \sum_{k=0}^\infty \discount^k   c_l(\state_k,\act1_k,\act2_k)  \vert x_0 = x\big\} 
\label{eq:cumcostgame} \eeq
where as usual $\discount \in (0,1)$ is the discount factor.
\end{compactitem}

The non-cooperative assumption in game theory is that the players are interested in minimizing their individual cumulative costs only; they do not
collude.

\subsection{Nash equilibrium of general sum Markov game}  \index{Markov game! Nash equilibrium}
 \index{game theory! Markov game}
Assume that each player has complete knowledge of the other player's cost function. Then the  policies
${\pol1}^*, {\pol2}^*$ of the non-cooperative infinite horizon Markov game constitute a Nash equilibrium if 
\beq \label{eq:nasheq}
\begin{split}
 \Jc{1}_{{\pol1}^*,{\pol2}^*}(x) &\leq  \Jc{1}_{\pol1,{\pol2}^*}(x),  \quad \text{ for all }   \pol1 \in \Policy_S \\
  \Jc{1}_{{\pol1}^*,{\pol2}^*}(x) &\leq  \Jc{1}_{ {\pol1}^*,{\pol2}}(x),  \quad \text{ for all  } \pol2 \in \Policy_S.
  \end{split} \eeq
This means that unilateral deviations from ${\pol1}^*, {\pol2}^*$ result in either player being worse off (incurring a larger cost).
Since in a non-cooperative game collusion is not allowed, there is no rational reason for players to deviate from the Nash equilibrium
(\ref{eq:nasheq}).

In game theory, two important issues are:  
\begin{compactenum}
\item {\em Does a Nash equilibrium exist? } 
For the above discounted cost game with finite action and state space, the answer is "yes".
\begin{theorem}
A discounted Markov game has at least one Nash equilibrium within the class of Markovian stationary (randomized) policies.
\end{theorem}
The proof is in \cite{FV12} and  involves Kakutani's fixed point theorem.\footnote{Existence proofs for equilibria  involve using either Kakutani's  fixed point theorem  (which generalizes Brouwer's fixed point theorem to set valued correspondences) or Tarski's fixed point theorem (which applies to supermodular games). Please see \cite{MWG95} for a nice intuitive  visual illustration
of these fixed point theorems.}

\item {\em How can the Nash equilibria be computed?} Define the randomized policy of player 1  (corresponding to $\pol1$)  and player 2  
(corresponding to $\pol2$)  as 
$$p(i,\act1) = \prob(\act1_k = \act1 | \state_k = i),  \quad q(i,\act2) = \prob(\act2_k = \act2 | \state_k = i) $$
Then for an infinite horizon discounted cost Markov game, the Nash equilibria $(p^*,q^*)$ are global optima  of the following  non-convex optimization problem:
\beq \begin{split} &  \text{ Compute }   \max \sum_{l=1}^2 \sum_{i=1}^\statedim \alpha_i \bigg(  \uV^{(l)}(i)
-  \sum_{\act1,\act2}  \cost_l(i,\act1,\act2) p(i,\act1) q(i,\act2)  \\ & -  
\discount  \sum_{j \in \statespace} \sum_{\act1,\act2}  \tp_{ij}(\act1,\act2) p(i,\act1) q(i,\act2) \uV^{(l)}(j) \bigg) \\ &\text{ with respect to $(\uV^{(1)},\uV^{(2)},p, q)$}    \\
 \text{ subject to } &\uV^{(1)}(i) \leq \sum_{\act2}  \cost(i,\act1,\act2) q(i,\act2) +  
\discount  \sum_{j \in \statespace} \sum_{\act2}  \tp_{ij}(\act1,\act2) q(i,\act2) \uV^{(1)}(j) ,\\
&\uV^{(2)}(i) \leq \sum_{\act2}  \cost(i,\act1,\act2) p(i,\act1) +  
\discount  \sum_{j \in \statespace} \sum_{\act1}  \tp_{ij}(\act1,\act2) p(i,\act1) \uV^{(2)}(j) ,\\
& q(i,\act2) \geq 0,  \quad \sum_{\act2} q(i,\act2) = 1,  \quad i = 1,2,\ldots, \statedim ,\;  \act2 = 1,\ldots,\actiondim
\\
& p(i,\act1) \geq 0,  \quad \sum_{\act1} p(i,\act1) = 1,  \quad i = 1,2,\ldots, \statedim ,\;  \act1 = 1,\ldots,\actiondim.
\end{split} \label{eq:genmarkovgame}  \eeq
 In general, solving the non-convex optimization problem (\ref{eq:genmarkovgame}) is difficult; there
 can be multiple  global optima (each corresponding to a Nash equilibrium) and multiple local optima. In fact there is a fascinating property that if all
 the parameters (transition probabilities, costs) are rational numbers,  the Nash equilibrium policy can involve irrational numbers. This points to the 
 fact that in general one can only approximately compute the Nash equilibrium.

\index{Markov game! general sum}
\end{compactenum}

{\bf Proof}.
First write  (\ref{eq:genmarkovgame}) in more abstract but intuitive notation in terms of the randomized policies $p,q$ as
\beq \begin{split} &  
 \max \sum_{l=1}^2  \alpha^\p \bigg(  \uV^{(l)} - c_l(p,q) - \rho \tp(p,q) \uV^{(l)}  \biggr) \\
\text{ subject to } \; &  \uV^{(1)} \leq c_1(\act1,q) + \rho \tp(\act1,q) \uV^{(1)} , \quad \act1 = 1,\ldots,\actiondim \\
 &  \uV^{(2)} \leq c_2(p,\act2) + \rho \tp(p,\act2) \uV^{(2)} , \quad \act2 = 1,\ldots,\actiondim \\
&  p , q \text{ valid pmfs } 
 \end{split} \label{eq:genmarkovgame2}  \eeq
 It is clear from the constraints that the objective function is always $\leq 0$. In fact the maximum is attained when the objective function
 is zero, in which case the constraints hold with equality.  When the constraints hold at equality, they  satisfy
  $$  V_*^{(l)}  = \big(I - \discount \tp(p^*,q^*) \big)^{-1}  c_l(p^*,q^*) ,  \quad l = 1,2. $$
  This serves as definition of $V_*^{(l)}$ and 
is equivalent to saying\footnote{This holds since from (\ref{eq:cumcostgame}), $  \Jc{l}_{p^*,q^*}(x) = c_l(p^*,q^*) + 
 \discount \tp  c_l(p^*,q^*)  + \discount^2 \tp^2 c_l(p^*,q^*) + \cdots + $. Indeed a  similar expression holds for discounted cost MDPs.}
 that $V_*^{(l)}$ is the  infinite horizon cost attained by the policies $(p^*,q^*)$. That is,
\beq  V_*^{(l)}  =  c_l(p^*,q^*) + \rho \tp(p^*,q^*) V_*^{(l)}  \implies 
 \Jc{l}_{p^*,q^*}(x)  = V_*^{(l)}.  \label{eq:neeqa} \eeq

Also  setting $\uV^{(l)} =V_*^{(l)}$,  the constraints in (\ref{eq:genmarkovgame2})  satisfy
$$  V_*^{(1)} \leq c_1(p,q^*) + \rho \tp(p,q^*)V_*^{(1)} ,  \quad    V_*^{(2)} \leq c_2(p^*,q) + \rho \tp(p^*,q) V_*^{(2)}  $$
implying that
\beq \Jc{1}_{p,q^*}(x)  \geq  V_*^{(1)},  \quad   \Jc{2}_{p^*,q}(x)  \geq  V_*^{(2)}.   \label{eq:neeqb}  \eeq
(\ref{eq:neeqa}) and (\ref{eq:neeqb}) imply that $(p^*,q^*)$ constitute a Nash equilibrium. \qed
 
{\em  Remark}.
The reader should compare the above proof  with the linear programming formulation for a discounted cost MDP. In that derivation we started with
a similar  constraint
\beq  \uV \leq c_1(\act1) + \rho \tp(\act1) \uV . \label{eq:edmdpc} \eeq
 This implies that $\uV < V$  where $V$ denotes the unique value function of Bellman's equation. Therefore the objective was to find $\max \alpha^\p \uV$
subject to (\ref{eq:edmdpc}). 
So in  MDP case we obtain a linear program.
In the dynamic game case, in general, there is no value function to clamp (upper bound) $\uV$.


\subsection{Zero-sum discounted Markov game} \index{Markov game! zero sum}
With the above brief introduction,
the main aim below is to give  special cases of {\em zero-sum} Markov games where the Nash equilibrium can be 
computed via linear programming. (Recall \secn \ref{sec:vipilp} of the book shows how a discounted cost MDP can be solved via linear programming.)

A discounted Markovian game is said to be zero sum\footnote{A constant sum game 
 $c_1(\state,\act1,\act2) +  c_2(\state,\act1,\act2) = K$ for constant $K$ is equivalent to a zero sum game.
 Define $\bar{c}_l(\state,\act1,\act2) = c_l(\state,\act1,\act2)+K/2 $, $l=1,2$,   resulting in a zero sum game in terms of $\bar{c}_l$.}
 if 
$$ c_1(\state,\act1,\act2) +  c_2(\state,\act1,\act2) = 0.$$ 
That is, 
$$c(\state,\act1,\act2)  \ole c_1(\state,\act1,\act2) = - c_2(\state,\act1,\act2) .$$
For a zero sum game, the Nash equilibrium (\ref{eq:nasheq}) becomes  a saddle point:
$$J_{{\pol1}^*,{\pol2}}(x) \leq J_{{\pol1}^*,{\pol2}^*}(x)  \leq J_{{\pol1},{\pol2}^*}(x),
$$
that is, it is a minimum in the $\pol1$ direction and a maximum in the $\pol2$ direction.

A well known result from the 1950s due to Shapley is: \index{Markov game! Shapley's theorem}
\begin{theorem}[Shapley] A zero sum  infinite horizon discounted cost Markov game has a unique value function, even though
there could be multiple Nash equilibria (saddle points).  Thus all the Nash equilibria are equivalent.
\end{theorem}
The  value function of the zero-sum game is 
$$ J_{{\pol1}^*,{\pol2}^*}(i) = V(i)$$ where
$V$ satisfies an  equation that resembles  dynamic programming:
\beq  V(i) = \val \big[ (1-\discount) c(i,\act1,\act2) + \discount \sum_{j} \tp_{ij}  (\act1,\act2) V(j) \big]_{\act1,\act2} \label{eq:dpzg} \eeq
Here $\val[M]_{\act1,\act2}$ denotes the value of the matrix\footnote{A  zero sum matrix game is of the form: Given a $m\times n$ 
matrix $M$, determine the Nash equilibrium
$$ (x^*, y^*) =  \argmax_x \argmin_y y^\p M x, \quad  \text{ where } x, y \text{ are probability vectors }$$
The value of this matrix game is $\val[M] = {y^*}^\p M x^*$ and is  computed as the solution of a linear programming (LP) problem as follows:
Clearly $ \max_ x \min_y y^\p M x = \max_x \min_{i} e_i^\p M x $ where $e_i$, $i=1,2,\ldots,m$ denotes the unit $m$-dimensional vector with 1 in the $i$-th position. This follows since a linear function is minimized at its extreme points. So the minimization over continuum has been reduced
to one over a finite set.  Denoting $z = \min_{i} e_i^\p M x $, the value of the game is the solution of the following LP:
\beq  \val[M] = \begin{cases}
  \text{ Compute }    \max z  \\
   z < e_i^\p M x  ,  \quad i =1,2,\ldots, m, \\
      \one^\p x = 1 ,  \quad x_j \geq 0,  j=1,2\ldots,n \end{cases}
\eeq
}
  game with elements  $M(\act1,\act2)$. 
Even though for a specific vector $V$, the $\val[\cdot]$ in the right hand side of (\ref{eq:dpzg}) can be evaluated by solving an LP,  it is not useful  for the Markov zero sum
game, since we have a functional equation in the variable $V$. So solving a zero sum Markov game is difficult in general.

\subsubsection*{Nash Equilibrium as a Non-convex Bilinear Program}  \index{Markov game! Nash equilibrium as bilinear program}
To give more insight, as we did in the discounted cost MDP case, let us formulate computing the Nash equilibrium (saddle point) of  the zero sum Markov game as an optimization problem. In the MDP case we obtained a LP; for the Markov  game (as shown below) we obtain a non-convex
bilinear optimization problem.

Define the randomized policy of player 1 (minimizer) and player 2 (maximizer) as 
$$p(i,\act1) = \prob(\act1_k = \act1 | \state_k = i),  \quad q(i,\act2) = \prob(\act2_k = \act2 | \state_k = i) $$
In complete analogy to the discounted MDP  case in (\ref{eq:discountprimal}), player 2 optimal strategy $q^*$ is the solution of the bilinear program
\beq \begin{split}  &  \max \sum_i \alpha_i  \uV(i)  \;\text{ with respect to $(\uV, q)$}   \\
& \text{ subject to } \uV(i) \leq \sum_{\act2}  \cost(i,\act1,\act2) q(i,\act2) +  
\discount  \sum_{j \in \statespace} \sum_{\act2}  \tp_{ij}(\act1,\act2) q(i,\act2) \uV(j) ,\\
& q(i,\act2) \geq 0,  \quad \sum_{\act2} q(i,\act2) = 1,  \quad i = 1,2,\ldots, \statedim ,\;  \act2 = 1,2,\ldots,\actiondim.
\end{split} \label{eq:gdiscountprimal}  \eeq
By symmetry, player 1 optimal  strategy $p^*$ is the solution of the bilinear program
\beq \begin{split}  &  \min \sum_i \alpha_i  \uV(i)  \;\text{ with respect to $(\uV, p)$}   \\
& \text{ subject to } \uV(i) \geq \sum_{\act2}  \cost(i,\act1,\act2) p(i,\act1) +  
\discount  \sum_{j \in \statespace} \sum_{\act1}  \tp_{ij}(\act1,\act2) p(i,\act1) \uV(j) ,\\
& p(i,\act1) \geq 0,  \quad \sum_{\act1} p(i,\act1) = 1,  \quad i = 1,2,\ldots, \statedim ,\;  \act1 = 1,2,\ldots,\actiondim.
\end{split} \label{eq:gdiscountprimal2}  \eeq
The key difference between the above discounted Markov game problem and the discounted MDP (\ref{eq:discountprimal}) is that the above equations
are no longer LPs. Indeed the constraints are {\em bilinear} in $ (\uV,q)$ and $(\uV,p)$. So the constraint set for a zero-sum Markov game is non-convex. Despite (\ref{eq:gdiscountprimal}) and (\ref{eq:gdiscountprimal2}) being nonconvex, in light of Shapley's theorem all local minima are global minima.

Finally (\ref{eq:gdiscountprimal}) and (\ref{eq:gdiscountprimal2}) can be combined into a single optimization problem. To summarize, the (randomized)  Nash equilibrium $p^*,q^*$ of 
a zero-sum Markov game is the solution of the following bilinear (noconvex) optimization problem:
\beq \begin{split} &    \max \sum_i \alpha_i \big(  \uV^{(1)}(i)  -\uV^{(2)}(i) \big) \;\text{ with respect to $(\uV^{(1)},\uV^{(2)},p, q)$}    \\
 \text{ subject to } &\uV^{(1)}(i) \leq \sum_{\act2}  \cost(i,\act1,\act2) q(i,\act2) +  
\discount  \sum_{j \in \statespace} \sum_{\act2}  \tp_{ij}(\act1,\act2) q(i,\act2) \uV^{(1)}(j) ,\\
&\uV^{(2)}(i) \geq \sum_{\act2}  \cost(i,\act1,\act2) p(i,\act1) +  
\discount  \sum_{j \in \statespace} \sum_{\act1}  \tp_{ij}(\act1,\act2) p(i,\act1) \uV^{(2)}(j) ,\\
& q(i,\act2) \geq 0,  \quad \sum_{\act2} q(i,\act2) = 1,  \quad i = 1,2,\ldots, \statedim ,\;  \act2 = 1,2,\ldots,\actiondim
\\
& p(i,\act1) \geq 0,  \quad \sum_{\act1} p(i,\act1) = 1,  \quad i = 1,2,\ldots, \statedim ,\;  \act1 = 1,2,\ldots,\actiondim.
\end{split} \label{eq:cgdiscountprimal2}  \eeq

\subsubsection*{Special cases where computing Nash Equilibrium is an LP} \index{Markov game! linear programming}
We now give two special examples of zero-sum Markov games that can be solved as a linear programming problem (LP); single controller
games and switched controller games. In both cases the bilinear terms in (\ref{eq:cgdiscountprimal2}) vanish and the computing the Nash
equilibrium reduces to solving linear programs.

\subsection{Example 1. Single Controller zero-sum Markov Game} \index{Markov game! single controller}
In a single controller Markov game, the transition probabilities are controlled by one  player only;  we assume that this is player 1. So 
$$ \tp_{ij}(\act1,\act2) =   \tp_{ij}(\act1) = \prob(\state_{k+1} = j | \state_k = i, \act1_k = \act1) $$
 Due to this assumption, the bilinear constraint in (\ref{eq:gdiscountprimal}) becomes {\em linear}, namely 
$$ \uV(i) \leq \sum_{\act2}  \cost(i,\act1,\act2) q(i,\act2) +  
\discount  \sum_{j \in \statespace}  \tp_{ij}(\act1)  \uV(j) $$
since $\sum_{\act2} q(i,\act2) = 1$. Therefore  (\ref{eq:gdiscountprimal}) is now an LP which can be solved for $q^*$, namely:
 \beq \begin{split}  &  \max_{\uV} \sum_i \alpha_i  \uV(i)  \;\text{ with respect to $(\uV, q)$}   \\
& \text{ subject to }  \uV(i) \leq \sum_{\act2}  \cost(i,\act1,\act2) q(i,\act2) +  
\discount  \sum_{j \in \statespace}  \tp_{ij}(\act1)  \uV(j) ,   \\
& q(i,\act2) \geq 0,  \quad \sum_{\act2} q(i,\act2) = 1,  \quad i = 1,2,\ldots, \statedim ,\; \act2 = 1,2,\ldots,\actiondim.
\end{split} \label{eq:zdiscountprimal}  \eeq
Solving the above LP yields the Nash equilibrium policy $\pol2$ for player 2.

The dual problem to (\ref{eq:zdiscountprimal}) is the linear program 
\beq \begin{split}
 \text{ Minimize } &  \sum_{i\in \statespace} 
 z(i) 
 \;\text{ with respect to $(z, p)$} \nonumber \\
\text{ subject to }  & p({i,\act1}) \geq 0, \quad i \in \statespace, \action  \in \actionspace \\
&  \sum_{\act1} p({j,\act1} )= \discount\, \sum_{i} \sum_{\act1}  \tp_{ij}(\act1)\, p({i,\act1}) + \alpha_j ,\;  j \in \statespace. \\
& z(i) \geq  \sum_{\act1} p(i,\act1)\, c(i,\act1,\act2)
\end{split}
\label{eq:gdiscountdual}
\eeq
 The above dual gives the randomized Nash equilibrium policy $p^*$  for player 1.
 
 \subsection{Example 2. Switching Controller Markov Game} \index{Markov game! switched controller}
 This is a special case of a zero sum Markov game where the state space $\statespace$ is partitioned into disjoint sets $\sr1,\sr2$ such
 that $\sr1 \cup \sr2 = \statespace$ and 
 $$ \tp_{ij}(\act1,\act2) = \begin{cases}  \tp_{ij}(\act1),  & i \in \sr1 \\
   							  \tp_{ij}(\act2), & i \in \sr2 \end{cases} $$
So for states in $\sr1$, controller 1 controls that transition matrix, while for states in $\sr2$, controller 2 controls the transition matrix.

Obviously for $ i \in\sr2$, (\ref{eq:gdiscountprimal}) becomes  an linear program while for $i \in \sr1$, (\ref{eq:gdiscountprimal2}) becomes a
linear program.
As discussed in \cite{FV12}, the Nash equilibrium can be computed by solving  a finite sequence of linear programming problems.

\chapter{Partially Observed Markov Decision Processes (POMDPs)}  

Several well studied instances of POMDPs and their parameter files can be found at \url{http://www.pomdp.org/examples/}
\\ \\

\begin{compactenum}
\item Much insight can be gained by simulating the dynamic programming recursion for a 3-state POMDP.  The belief
state needs to be quantized to a finite grid.  We also strongly recommend using the exact POMDP solver in \cite{POMDP} to gain
insight into the piecewise linear concave nature of the value function.

\item Implement Lovejoy's suboptimal algorithm and compare its performance with the optimal policy.

\item {\bf Tiger problem}:   \index{POMDP tiger problem} This is a colorful name given to the following POMDP problem.

A tiger resides behind one of two doors, a left door $(l)$ and a right door $(r)$.
The state $\state \in \{l,r\}$ denotes the position of a tiger. The action $\action \in \{l,r,h\}$ denotes a human either opening the left door $(l)$,
opening the right door $(r)$, or simply hearing $(h)$  the growls of the tiger. If the human opens a door, he gets a perfect measurement
of the position of the tiger (if the tiger is not behind the door he opens, then it must be behind the other door). If the human chooses action
$h$ then he hears the growls of the tiger which gives noisy information about the tiger's position. Denote the probabilities
$\oprob_{ll}(h) = p$, $\oprob_{rr}(h) = q$.

Every time the human chooses the action to open a door, the problem resets and the tiger is put with equal probability behind one of the doors.
(So the transition probabilities for the actions $l$ and $r$ are $0.5$). 

The cost of  opening the door behind where the tiger is hiding is $\alpha$, possibly reflecting injury from the tiger. The cost
of opening the other door is $-\beta$ indicating a reward. Finally the cost of hearing and not opening a door is $\gamma$.

The aim is to minimize the cost (maximize the reward) over a finite or infinite horizon.
To summarize, the POMDP parameters of the tiger problem are:
\begin{align*} \statespace &= \{l,r\},  \obspace= \{l,r\},  \actionspace = \{l, r, h\},  \\
\oprob(l) &= \oprob(r) = I_{2 \times 2}, \oprob(h) = \begin{bmatrix} p & 1-p \\ 1-q & q \end{bmatrix} \\
\tp(l) &= \tp(r) = \begin{bmatrix} 0.5 & 0.5 \\ 0.5 & 0.5 \end{bmatrix} ,
\tp(h) = I_{2 \times 2}, \\
\cost_l &= (\alpha, -\beta)^\p, \cost_r = (-\beta, \alpha)^\p, \cost_h = (\gamma, \gamma)^\p
\end{align*}


\item {\bf Open Loop Feedback Control.} As described in \secn \ref{sec:olfc}, open loop feedback control is a useful suboptimal scheme for solving
POMDPs. Is it  possible to exploit knowledge that the value function of a POMDP is piecewise linear and concave in
the design of an open loop feedback controller?

\item Finitely transient policies were discussed in \secn 7.6. For a 2-state,  2-action, 2-observation POMDP,  give an example of POMDP parameters that yield a finitely transient policy with $n^* = 2$.

\item {\bf Uniform sampling from Belief space.} Recall that the belief space $\Belief$ is the unit $X-1$ dimensional simplex.
Show that a convenient way of sampling uniformly from $\Belief$ is to use 
the Dirichlet distribution $$\belief_0(i) = \frac{x_i}{\sum_{j=1}^X  x_j}, \quad \text{  where }  x_i \sim \text{  unit exponential distribution. }  $$ \index{Dirichlet distribution}

\item {\bf Adaptive Control of a fully observed MDP formulated as a POMDP problem.}
\label{prob:adaptive control}  \index{adaptive control of MDP as a POMDP}
Consider a fully observed MDP with transition matrix $\tp(u)$ and cost $c(i,u)$, where $u\in \{1,2,\ldots,\actiondim\}$ denotes the action. Suppose the true transition matrices 
$\tp(u)$ are not known.
However, it is known apriori that they  belong to a known finite set of matrices $\tp(u,\th)$ where   $\th \in \{1,2, ,\ldots, L\}$. As data accumulates, the controller
must simultaneously control the Markov chain and also estimate the transition matrices.

The above problem can be formulated straightforwardly as a POMDP.
Let $\th_k$ denote the parameter  process. Since the parameter  $\th_k = \th$ does not evolve with time,  it has identity transition matrix.
Note that $\th$ is not known; it is partially observed since we only see the sample path realization of the Markov chain $\state$ with transition matrix
$\tp(\action,\th)$.\\
{\bf Aim}: Compute  the optimal policy
$$ \mu^* =  \argmin_\mu  J_\mu(\belief_0) =  \E\{ \sum_{k=0}^{\finaltime-1}  c\big(\state_k, u_k 
 \big) | \belief_0 \} $$
where $\belief_0$ is the prior pmf of $\th$. The key point here is that as in a POMDP (and unlike an  MDP),
the action $u_k$ will now depend on the history of past actions and the trajectory of the Markov chain as we will now describe.

{\bf Formulation}:
Define the augmented state $(\state_k, \th_k)$.  Since $\th_k =\th$ does not evolve, clearly the augmented state  has transition probabilities
$$ \prob(\state_{k+1} = j, \th_{k+1} = m | \state_{k} = i, \th_k = l, \action_k = u)  = \tp_{ij}(\action,l)\,   \delta(l-m), \quad 
m=1,\ldots, L. $$
At time $k$, denote the history as
$\history_k = \{\state_0,\ldots,\state_k, \action_1,\ldots,\action_{k-1}\}$.  Then  define the belief state which is the posterior pmf of the model
parameter estimate:
$$ \belief_k(l) = \prob( \th_k = l| \history_k) , \qquad l = 1,2.\ldots, L.$$
\begin{compactenum}
\item Show that  the posterior is updated  via  Bayes' formula as
\begin{equation}
\begin{split}
 \belief_{k+1}(l) &=  T(\belief_k, \state_k,\state_{k+1},\action_k)(l) \ole \frac{\tp_{\state_k,\state_{k+1}}(\action_k, l) \,  \belief_k( l)}
{\sigma(\pi_k,x_k,x_{k+1})}, \;   l = 1,2.\ldots, L \\ \text{ where } \; &
\sigma(\pi_k,x_k,x_{k+1},u_k) =  \sum_m \tp_{\state_k,\state_{k+1}}(\action_k, m)\,   \belief_k( m).  \end{split}
\eeq
Note that $\belief_k$ lives in the $L-1$ dimensional unit simplex.

Define the belief state as $(\state_k, \belief_k)$. The actions are then chosen as
$$\action_{k} = \mu_k(\state_k,\belief_k) $$

Then the optimal policy $\mu_k^*(i,\belief)$ satisfies Bellman's equation
\beq \begin{split}
J_k(i,\belief) &= \min_u  Q_k( i, u, \belief) ,  \quad \mu^*_k(i,\belief) = \argmin_u  Q_k( i, u, \belief) \\
Q_k(i,u,\belief) &= 
c(i,u) +  \sum_{j}   J_{k+1}\big(j, T(\belief,i,j,u) \big) \,  
\sigma(\pi,i,j,u)
 \end{split}
\eeq 
initialized with the terminal cost $J_N(i,\belief) =   c_N(i)$.

\item Show that the value  function $J_k$ is piecewise linear and concave in $\belief$. Also show how the exact POMDP solution algorithms 
in  Chapter \ref{ch:pomdpbasic}
can be used to compute
the optimal policy. 
\end{compactenum}
The above problem is related to the concept of {\em dual control} which dates back to the 1960s \cite{Fel65}; see also
\cite{Lov93} for the use of Lovejoy's suboptimal algorithm to this problem.  \index{dual control}
Dual control relates to the tradeoff between estimation and control:  if the controller is uncertain about the model
 parameter, it needs to control   
the system more aggressively in order to probe the system to estimate it;  if the controller is more certain about  the model parameter, 
it can deploy a less aggressive control.
In other words, initially the controller explores and as the controller becomes more certain it exploits. Multi-armed bandit problems optimize the
tradeoff between exploration and exploitation. \index{dual control}

\item{\bf Optimal Search and Dynamic (Active)  hypothesis testing.}  \index{dynamic hypothesis testing}  \index{active hypothesis testing}   In \secn 7.7.4 of the book,  we considered the classical optimal  search problem where the objective was to search for a non-moving target amongst a finite number of cells.
A crucial assumption was that there  are no false alarms;  if an object is not present in a cell and the cell is searched, the observation recorded
is $\bar{F}$ (not found). 

A generalization of this problem is studied in \cite{Cas95}.  Assume there are $\actionspace = \{1,2,\ldots,U\}$ cells. When cell $u$ is searched
\begin{compactitem}
\item If the target is in cell $u$ then an observation $y$ is generated with pdf or pmf $\phi(y)$ if 
the target is in cell $u$
\item If the target is {\bf not} in cell $u$, then an observation $y$ is generated with pdf or pmf $\bar{\phi}(y)$. 
(Recall in classical search $\bar{\phi}(y)$ is  dirac measure on the observation symbol $\bar{F}$.)
\end{compactitem}
The aim is to determine the optimal search policy $\bpolicy$ over a time
horizon $\finaltime$
to maximize
$$ J_\bpolicy = \E_{\bpolicy} \max_{u \in \{1,\ldots, U\}}  \belief_\finaltime (u) \} $$
at the final time $\finaltime$.

Assume  the pdf or pmf $\bar{\phi}(y)$ is symmetric in $y$, that is  $\bar{\phi}(y) =  \bar{\phi}(b-y)$  for some real constant $b$. Then \cite[Proposition 3] {Cas95} shows the nice
result that the optimal policy is to search either of the two most likely locations given the belief $\belief_k$.

The above problem can be viewed as an active  hypothesis testing problem, which is an instance of a  controlled
sensing problem.
The decision  maker seeks to adaptively select the most informative sensing action for making a decision in a hypothesis testing problem.  Active hypothesis testing goes  all the way
back to the 1959 paper by 
Chernoff \cite{Che59}.
For a more general and recent take of active hypothesis testing please see~\cite{NJ13}. 
\end{compactenum}

\chapter{POMDPs in Controlled Sensing and Sensor Scheduling}  

\begin{compactenum}

\item {\bf Optimal Observer Trajectory for Estimating a Markovian Target.} This problem is identical to the search problem described in
\secn \ref{chp:optimal search}. A target moves in space according to a Markov chain. (For convenience assume $\statedim$-cells in two dimensional
space.
A moving  observer (sensor) measures the target's state (position) in noise. Assume that the noise depends on the relative distance between the target
and the observer.  \index{optimal observer trajectory}
How should the observer move amongst the $\statedim$-cells in order to locate where the target is?  
One metric that has been used in the literature \cite{LI99}  is the stochastic observability (which is related to the mutual information) of the target; see also \secn \ref{sec:radarkf}.  The aim of the observer is to move so as to maximize
the stochastic observability of the target. As described in \secn \ref{chp:optimal search}, the problem is equivalent  to a POMDP.

A more fancy version of the setup involves multiple observers (sensors) that move within the state space and collaboratively seek to locate the 
target. Assume that the observers exchange information about their observations and actions. The problem can again be formulated as a POMDP with a larger action and observation space.

Suppose the exchange of information between the observers occurs over a noisy communication channel where the error probabilities
evolve  according to a Markov chain as in \secn \ref{sec:minh}. Formulate the problem as a POMDP.

\item {\bf Risk averse sensor scheduling.}  As described in \secn \ref{sec:nonlinearpomdpmotivation}, in controlled sensing applications, one is interested in incorporating
the uncertainty in the state estimate into the instantaneous cost. This cannot be modeled  using a linear cost since the uncertainty is minimized
at each vertex of the simplex $\Belief$. In \secn  \ref{sec:nonlinearpomdpmotivation}, quadratic functions of the belief were used to model the conditional
variance. A more principled  alternative is to use  dynamic coherent risk measures; recall three examples of such risk measures were discussed
in \secn \ref{sec:riskaverse}.  

Discuss how open loop feedback control can be used for a POMDP with dynamic coherent risk measure.

 \item  {\bf Sensor Usage Constraints.} The aim here is to how the POMDP formulation  of a controlled sensing problem can be modified 
 straightforwardly to incorporate sensing 
 constraints 
 on the total usage of particular sensors.
Such constraints  are often used in 
sensor resource management.

\begin{compactenum}
\item Consider a $N$ horizon problem  \index{sensor management! usage constraints}
where  sensor 1 can be used at most $L$ times
where $L \leq N$. For notational simplicity, assume that there are two sensors, so $\actionspace = \{1,2\}$.
Assume that there are no constraints on the usage of the other sensors.

For notational convenience we consider rewards denoted as $R(\belief,u) = \sum_{i=1}^X R(i,u) \belief(i) $ instead of costs $C(\pi,u)$ expressed in terms of the belief state
$\belief$.
Show that Bellman's equation is given by
\begin{multline*}
 V_{n+1}(\belief,l) = \max\{ R(\belief,1) + \sum_\obs V_{n}(\filter(\belief,\obs,1), l-1) \filterd(\belief,\obs,1), \\
R(\belief,2) + \sum_\obs V_{n}(\filter(\belief,\obs,2), l) \filterd(\belief,\obs,2)\}
\end{multline*}
with boundary  condition $V_n(\belief,0) = 0$, $n=0,1,\ldots, N$.

\item If the constraint is that sensor 1 needs to be used exactly $L$ times, then show that the following additional  boundary condition needs to be included:
$$V_n(\belief,n) = R(\belief,1) + \sum_\obs V_{n-1}(\filter(\belief,\obs,1), n-1) \filterd(\belief,\obs,1) , \; \text{ for } n=1,\ldots,L .$$

\item In terms of the POMDP solver software, the  constraint for using sensor~1 at most $L$ times is  easily incorporated by augmenting the state space. Define the controlled finite state process $r_k \in \{0,2,\ldots,L\}$ with $(L+1 )\times (L+1)$ transition matrices
$$ Q(1) = \begin{bmatrix} 0 & 1 & 0 & \cdots & 0\\
                0 & 0 & 1 & \cdots & 0 \\
        \vdots & \vdots& \vdots & \ddots & 1 \\
        0 & 0 & 0 & \cdots & 1 \end{bmatrix}, \quad Q(2) = I .$$
Then define the POMDP with:
\begin{compactitem}
\item transition matrices $\tp(1) \otimes Q(1)$ and $\tp(2) \otimes Q(2)$,
\item   observation probabilities
$\pdf(\obs | \state, r,\action) = \pdf(\obs| \state, \action)$, 
\item rewards $R(x,r,u) = R(x,u)$ for $r> 0$ and   $R(x,r=0,u) = 0$.
\end{compactitem}

In the problems for Chapter \ref{ch:pomdpstop}, we consider a simpler version of the above problem for optimal measurement selection of a HMM.
In that simpler case, one can develop structural results for the optimal policy.
\end{compactenum}

\item As described in \secn \ref{sec:nonlinearpomdpmotivation},  in controlled sensing it makes sense to choose a cost that is nonlinear in the belief state  $\belief$ in order to penalize uncertainty in the state estimate.
One choice of a nonlinear cost that has zero cost at the vertices of the belief space is 
$$C(\belief,\action) = \min_{ i   \in \{1,\ldots,\statedim\}} \belief(i) . $$
This cost $C(\belief,\action)$  is  piecewise linear and concave in $\belief \in \Belief$ where $\Belief$ denotes the belief space.

Since $C(\belief,\action)$ is positively homogeneous, show that the value function is piecewise linear and concave for any finite horizon $\finaltime$.
Hence the optimal POMDP solvers of Chapter  \ref{ch:pomdpbasic} can be used to solve this nonlinear cost POMDP exactly and therefore compute the optimal
 policy.

\end{compactenum}

\chapter{Structural Results for Markov Decision Processes}  

\begin{compactenum}

\item {\bf Supermodularity, Single Crossing Condition \& Interval Dominance Order.} \index{supermodular} \index{single crossing}  \index{interval dominance order}
A key step in establishing structural results for  MDPs is to  give sufficient conditions for  $u^*(x) = \argmax_u \fun(x,u)$ to be increasing in $x$.
In \secn \ref{chp:fullsupermod} of Chapter \ref{chp:monotonemdp} we gave two conditions, namely supermodularity and the single crossing condition (which
is a more general condition than supermodularity).
More recently, the interval dominance order has been introduced in \cite{QS09} as an even more general condition.
All three conditions boil down to the following statement:
\beq \fun(x+1,u+1) - \fun(x+1,u)  \geq  \rho(u) \,  \big(  \fun(x,u+1) - \fun(x,u)  \big)    \label{eq:intdominance} \eeq
where $\rho(u)$ is a strictly positive function of $u$.  In particular,
\begin{compactitem}
\item Choosing $\rho(u)=1$ in (\ref{eq:intdominance}) yields the supermodularity condition.
\item  If there exists a fixed positive constant  $\rho(u)$ such that  (\ref{eq:intdominance}) holds, then  the single crossing condition holds.
\item If there exists a positive function $\rho(u)$ that is increasing\footnote{Recall that in the book we use increasing in the weak sense to mean non-decreasing} in  $u$, then  (\ref{eq:intdominance}) yields the interval dominance order condition (actually this is a sufficient condition for interval dominance, see \cite{QS09} for details).
\end{compactitem}
Note that single crossing and interval dominance are ordinal properties in the sense that they are preserved by monotone transformations.

The sum of supermodular functions is supermodular.  Unfortunately, in general, the um of single crossing functions is not single crossing; however, see
\cite{QS12} for some results.
Discuss if the interval dominance order holds for sums of functions. Can it be used to develop structural results for an MDP?

\item Clearly, in general,  the sum of single crossing functions is not single crossing.  Even a constant plus a single crossing function
is not necessarily single crossing.  Sketch the curve of a single crossing function which wiggles close to zero. Then adding a positive constant
implies that the curve will cross zero more than once.
Also the sum of a supermodular plus single crossing is not single crossing.  In terms of  $\phi(x) = f(x,u+1) - f(x,u)$,  supermodular implies
$\phi(x) $ is increasing in $x$. Clearly the sum of an increasing function and a single crossing is not single crossing in general.

\item {\bf Invariance of optimal  policy to costs.} Recall that Theorem \ref{thm:mdpmonotone} require that the MDP costs satisfy assumptions
(A1) and (A3) for the optimal policy to be monotone.
Show that  for a discounted cost infinite horizon MDP, assumption (A1) and (A3) can be  relaxed as follows:

There exists a single vector $\fun \in \reals^\statedim$ such that for every action $\action \in \actionspace$, 
\begin{compactenum}
\item[(A1')]  $(I - \discount \tp(\action)) \fun  $ is a vector with increasing elements.
(Recall $\discount$ is the discount factor.)
\item[(A3')] $  (\tp(\action+1) - \tp(\action)) \fun$ is a vector with decreasing elements.
\end{compactenum}
In other words the structure of the transition matrix is enough to ensure a monotone policy and no assumptions are required on the cost
(of course the costs are assumed to be bounded)

{\em Hint}: Define the new value function
$\Vb(i) = V(i) - \fun(i)$ . Clearly the optimal policy remains unchanged and  $\Vb$ satisfies Bellman's equation
$$ \Vb(i) = \min_\action \{ c(i,u) - \fun(i) + \discount \sum_j \fun(j) \tp_{ij}(u)   + \discount  \sum_j  \Vb(j) \tp_{ij}(u)\} $$
where $\discount \in (0,1)$ denotes the discount factor.

\item {\bf Myopic lower bound to optimal policy.}
Recall that supermodularity of the transition matrix (A4) was a key requirement for the optimal policy to be monotone.
In particular,  Theorem  \ref{thm:mdpmonotone} shows that $Q(i,u)$ is submodular, i.e., 
$Q(i,u+1) - Q(i,u)$ is decreasing in $i$.
Sometimes supermodularity of the transition matrix is too much to ask for. 
Consider instead of (A4) the relaxed condition
\begin{compactitem}
\item[(A4')]$\tp_i(u+1) \gs \tp_i(u)$ for each row $i$.
\end{compactitem}
Show that (A4') together with (A1), (A2) implies that 
$$ \sum_j \tp_{ij}(u+1) V(j) \leq  \sum_j \tp_{ij}(u) V(j) $$
Define the myopic policy $\policyl(i) = \argmin_u c(i,u)$.
Show that 
under (A1), (A2), (A4'),
$\mu^*(i) \geq \underline{\mu}(i)$. In other words, the myopic policy $\policyl$ forms a lower bound to the optimal policy $\mu^*$.

\item {\bf Monotone policy iteration algorithm.} Suppose an MDP has a monotone policy.
 If the MDP parameters are known,
then the policy iteration algorithm of \secn \ref{sec:vipilp} can be used. If the policy $\policy_{n-1}$ at iteration $n-1$ is  monotone  then show that under the assumptions of (A1), (A2) of Theorem \ref{thm:mdpmonotone}, the policy evaluation step yields $J_{\mu_{n-1}}$ as a decreasing vector. Also show that  under (A1)-(A4),  (a similar  proof to Theorem \ref{thm:mdpmonotone}) implies that the policy improvement
 step yields  $\policy_n$ that is monotone. So the  policy iteration algorithm will automatically be confined to  monotone policies if initialized
 by a monotone policy.

\item {\bf Stochastic knapsack problem.} Consider the following version of the stochastic knapsack problem;\footnote{The classical NP hard knapsack problem deals with
$\actiondim$ items with costs $\cost(1),\cost(2),\ldots,\cost(\actiondim)$ and lifetimes $t_1,t_2,\ldots t_\actiondim$. The aim is to compute 
the minimum cost subset of these items whose total lifetime is at most~$T$.}
see \cite{Ros83} and also \cite{CR14a}.  \index{stochastic knapsack problem}
A machine  must operate for $T$ time points. Suppose that one specific component  of the machine fails intermittently.
This component is replaced when it fails. There are $\actiondim$-possible brands one can choose to replace this component when it fails.
Brand $\action\in \{1,2\ldots,\actiondim\}$  costs $\cost_\action$ and has an operating lifetime that is exponentially distributed with rate $\lambda_\action$. 
The aim is to minimize the expected total cost incurred by replacing the failed component so that the machine operates for $T$ time points.

Suppose a component has just failed. Let $t$ denote the remaining time left to operate the machine.
The optimal policy for deciding  which of the $\actiondim$ possible brands to choose the replacement  satisfies Bellman's equation
\begin{align*}  Q(t,u) &=  \cost(\action)  + \int_0^t V(t-\tau) \, \lambda_\action e^{-\lambda_\action \tau} d\tau , \quad Q(0,u) = 0, \\
V(t) &= \min_{u \in \{1,2,\ldots,\actiondim\}} Q(t,u) ,  \quad \optpolicy(t) = \argmin_{u \in \{1,2,\ldots,\actiondim\}} Q(t,u)
\end{align*}
Show that if $\lambda_\action \cost(\action)$ is decreasing with $\action$, then $Q(t,u)$ is submodular.
In particular, show that 
$$\frac{d}{dt}  Q(t,u) =  \lambda_\action \cost(\action) $$
Therefore, the optimal policy $\optpolicy(t)$  has the following structure:
Use brand 1 when the time remaining is small, then switch to brand 2 when the time increases, then brand 3, etc.

Generalize the above result to the case when time $k$ is discrete and the brand $u$ has life time pmf $\pdf(k,\action)$, $k=0,1\ldots$. 
Then Bellman's equation reads
\begin{align*} Q(n,u) &= \cost(\action) + \sum_{k=0}^n V(n-k) \, \pdf(k,\action) \\
 V(n) &= \min_{\action \in  \{1,2,\ldots,\actiondim\}}  Q(n,u), \quad  \optpolicy(n) = \argmin_{u \in \{1,2,\ldots,\actiondim\}} Q(n,u)
 \end{align*}
 What are sufficient conditions in terms of submodularity of  the lifetime pmf $\pdf(k,\action)$ for the optimal policy to be monotone?

 \item {\bf Monotonicity of optimal policy with respect to horizon.}  \index{effect of planning horizon} Show that the following result holds for a 
 finite horizon MDP.
 If $Q_n(i,u)$ is supermodular in $(i,u,n)$ then  $V_n(i) = \max_u Q_n(i,u)$ is supermodular in $i,u$.
 Note that  checking supermodularity with respect to  $(i,u,n)$ is pairwise: so it suffices to check
 supermodularity with respect to $(i,u)$, $(i,n)$ and $(u,n)$.
 
With the above result, consider  a finite horizon  MDP satisfies the assumptions (A1)-(A4) of \secn \ref{sec:monotonecond}.
Under what further conditions is 
  $\mu_n^*(i)$ is increasing in $n$ for fixed $i$?  What does this mean intuitively?

  \item {\bf Monotone Discounted Cost Markov Games.}    \index{structural result! Markov game}    \index{Markov game! structural result}  
    \index{Markov game! Nash equilibrium! structural result}  
  In \secn \vref{sec:markovgame} of this internet supplement we briefly described the formulation of infinite horizon discounted cost Markov games.
  Below we comment briefly on structural results for the Nash equilibrium of such games.
  
  Consider the infinite horizon discounted cumulative cost of (\ref{eq:cumcostgame}).
  The structural results developed in this chapter for MDPs extend straightforwardly to 
  infinite horizon discounted cost Markov games.
  The assumptions (A1) to (A4) of  \secn \ref{sec:monotonecond} of the book need to be extended as follows:
  
  \begin{description} 
\item[(A1)]  Costs $\cost(\state,\action,\action^-)$ are  decreasing in $\state$ and $\action^-$.  Here $\action^-$ denotes the actions of others players. 
\item[(A2)] $\tp_i(\action,\action^-) \ls \tp_{i+1}(\action,\action^-)$ for each $i$ and fixed $\action,\action^-$. 
Here $\tp_i(\action,\action^-)$ denotes the $i$-th row of the transition matrix for action $\action,\action^-$.
\item[(A3)] $\cost(\state,\action,\action^-)$ is submodular in $(\state,\action)$ and $(\action,\action^-)$
\item[(A4)] $\tp_{ij}(\action,\action^-)$ is tail-sum supermodular in $(i,\action,\action^-)$. That is,
$$\sum_{j\geq l} \big(\tp_{i j}(\action+1,\action^-) - \tp_{i j}(\action,\action^-) \big) \text{  is increasing in }  i .$$
\end{description}
\begin{theorem}
Under  conditions (A1)-(A4), there exists a pure Nash equilibrium $({\pol1}^*, {\pol2}^*)$ such that the pure policies
${\pol1}^*$ and ${\pol2}^*$ are increasing in state $i$. 
\end{theorem}
Contrast this with the case of a general Markov game (\secn \ref{sec:markovgame} of this internet supplement) where one can only
guarantee the existence of a randomized  Nash equilibrium
in general.
 
 The proof of the above theorem is as follows. First for any  increasing fixed  policy $\pol2$ for player 2, one can show via an identical proof to Theorem
 \ref{thm:mdpmonotone},   the optimal policy ${\pol1}^*(\state,\pol2(\state))$ is increasing in $\state$. Similarly, for any  increasing fixed  policy $\pol1$ for player 1, 
 ${\pol2}^*(\state,\pol1(\state))$ is increasing in $\state$.  These are obtained as the solution of Bellman's equation. In game theory, these are called best response strategies.
 Therefore the vector function $[{\pol1}^*(\state), {\pol2}^*(\state)]$ is increasing in   $\state$.
It then follows from Tarski's fixed point theorem\footnote{Let $X$ denote a compact lattice and $f: X \rightarrow X$ denote an increasing function.
Then there exists a fixed point $x^* \in X$ such that $f(x^*) = x^*$} that such a function has a fixed point. Clearly this  fixed point is a Nash equilibrium since any unilateral deviation
makes one of the players worse off.

Actually for submodular games a lot more holds.  The smallest and largest Nash equilibria are pure (non-randomized) and satisfy the monotone property
of the above theorem. These can be obtained
via a best response algorithm the simply iterates the best responses ${\pol1}^*(\state,\pol2(\state))$ and ${\pol2}^*(\state,\pol1(\state))$ until convergence.
There are numerous papers and books in the area.
 
\end{compactenum}

\chapter{Structural Results for Optimal Filters}  

\begin{compactenum}

\item In the  structural results presented in the book, we have only considered first order stochastic dominance and monotone likelihood ratio dominance (MLR)  since they are sufficient
for our purposes. Naturally there are many other concepts of stochastic dominance \cite{MS02}. Show that 
$$ \text{ MLR } \implies \text{Hazard rate order} \implies \text{first order}  \implies \text{second order} $$
Even though second order stochastic dominance is useful for concave decreasing functions (such as the value function of a POMDP),
just like first order dominance, it cannot cope with conditioning (Bayes' rule).

\item Consider a reversible Markov chain with transition  matrix $\tp$, initial distribution $\belief_0$ and stationary distribution $\belief_\infty$.
Suppose $\belief_0 \lr \belief_\infty$. Show that if $\tp$ has rows that are first order increasing then  $\belief_n \lr \belief_\infty$.

\item {\bf TPn matrix.}  A key assumption \ref{A3} in the structural results is that the transition matrix $\tp$ is TP2.  More generally, suppose $n =2,3,\ldots$.
Then a $\statedim\times \statedim$ 
matrix $\tp$ is said to be totally positive of order $n$ (denoted as TPn)  if for each $k \leq n$, all the $k \times k$ minors of $\tp$ are non-negative.

\item {\bf TP2 matrix properties.\footnote{Note that a TP2 matrix does not need to be a square matrix; we consider $\tp$ to be square here since it is a transition probability matrix.}} 
\secn \ref{sec:assumpdiscussion} gave some useful properties of TP2 matrices.

Suppose the $\statedim \times \statedim$ stochastic matrix $\tp$ is TP2. 
\begin{compactenum}
\item  Show that this implies that the elements satisfy
\begin{align*}  & \tp_{11} \geq \tp_{21} \geq \cdots \geq \tp_{X1} \\
 & \tp_{1X} \geq \tp_{2X} \geq \cdots \geq \tp_{XX}   \end{align*}
 \item
 Suppose $\tp$ has no null columns. Show that if $\tp_{ij} = 0$, then either $\tp_{kl} = 0$ for $k \leq i$ and $l \geq j$, or
 $\tp_{kl} = 0$ for $k \geq i$  and $l \leq j$.
 
 \item Show that 
 $$ e_1^\p  ({\tp^{n}})^\p   e_1  \downarrow n, \qquad    e_X^\p  ({\tp^{n}})^\p  e_1 \uparrow n.
 $$
Also show that  for each  $n$,
$$ e_1^\p  ({\tp^{n}})^\p e_i  \downarrow i,  \qquad e_X^\p  ({\tp^{n}})^\p e_i  \uparrow i
$$

 \end{compactenum}
 Please see \cite{KK77} for several other interesting properties of TP2 matrices.

\item MLR dominance is intimately linked with the TP2 property. Show that 
$$\belief_1 \lr \belief_2  \iff 
\begin{bmatrix} \belief_1^\p \\ \belief_2^\p \end{bmatrix} \text{ is TP2 .} $$

\item{\bf Properties of MLR dominance.} Suppose $X$ and $Y$ are random variables and recall that $\gr$ denotes MLR dominance.\footnote{Stochastic dominance is a property of the distribution of a random variable and has nothing to do with the random variable itself. Therefore in the book, we defined stochastic dominance in terms of the pdf or pmf. Here to simplify notation we use the random variable instead of its distribution.}
\begin{compactenum}
\item Show that $X\gr Y$ is equivalent to 
$$ \{ X | X \in A\} \gs \{Y | Y \in A\} $$
for all events $A$ with $P(X \in A) > 0$ and $P(Y\in A) > 0$ where $\gs$ denotes first order dominance.  This property is due to \cite{Whi82}.
\item Show that $X\gr Y$ implies that $g(X) \gr g(Y)$ for any increasing function $ g$.
\item Show that $X \gr Y$ implies that $\max\{X,c\} \gr \max\{Y,c\}$ for any positive constant $c$.

\item Under what conditions does $X \gr Y$ imply that $-X \lr -Y$?
\end{compactenum}
Do the above two properties hold for first order dominance?

\item {\bf MLR monotone optimal predictor.} Consider the HMM predictor given by the Chapman Kolmogorov equation
$ \belief_k = \tp^\p \belief_{k-1}$.
Show that  if $\tp$ is a TP2 matrix and $\belief_0 \lr \belief_1$, then 
$\belief_0 \lr \belief_1 \lr \belief_2 \lr \ldots$.

\item {\bf MLR constrained importance sampling.}
One of the main results of this chapter was to construct reduced complexity HMM filters that provably form lower and upper bounds to the optimal HMM filter
in the MLR sense.  
In this regard, consider the following problem. Suppose it is known that 
$ \ltp^\p \pi \lr  \tp^\p \pi$. Then given the reduced complexity computation of $\ltp^\p \pi$,
how can this be exploited to compute $\tp^\p \pi$?

It is helpful to think of the following toy example:
Suppose it is known that $x^\p p \leq 1$ for a positive vector $x$ and probability vector $p$. How can this constraint be exploited to actually compute the inner 
product 
$x^\p p$?  Obviously from a deterministic point of view there is little one can do to exploit this constraint. \index{constrained importance sampling}
But one can use constrained important sampling:  one simple  estimator is as follows:
$$ \frac{1}{N} \sum_{i=1}^N x_i I(x_i \leq 1)$$ 
where index $i$ is simulated iid from probability vector $p$.
In \cite{KR14} a more sophisticated constrained importance sampling approach is used to estimate $\tp^\p \pi$ by exploiting the constraint
$ \ltp^\p \pi \lr  \tp^\p \pi$.


\item {\bf Posterior Cramer Rao bound.} The posterior Cramer Rao bound \cite{TMN98}  for filtering can be used to compute a lower bound to the mean square error.  \index{posterior Cramer Rao bound}
This requires twice differentiability of the logarithm of the joint density. For HMMs, one possibility is to consider the Weiss-Weinstein
bounds \index{Weiss-Weinstein bounds}, see \cite{RO05}. Alternatively, the analysis of \cite{GSV05} can be used. Compare these 
with the sample path bounds for the HMM filter obtained in this chapter.

 \item  \index{shifted likelihood ratio order} The shifted likelihood ratio order is a stronger order than the MLR order. Indeed,  $p> q$ in the shifted likelihood ratio order sense
if $p_i/q_{i+j}$ is increasing in $i$ for any $j$. (If $j=0$ it coincides with the standard MLR order.)  What additional assumptions are required
to preserve the shifted likelihood ratio order under Bayes' rule?  Show that 
the shifted likelihood ratio order is closed under convolution.  How can this property be exploited to bound an optimal filter?

\item In deriving  sample path bounds for the optimal filter, we did not exploit the fact that $T(\pi,y)$ increases with $y$.
How can this fact be used in bounding the sample path of an optimal filter?

\item {\bf Neyman-Pearson Detector} \index{Neyman-Pearson detector}
Here we briefly review elementary  Neyman-Pearson detection theory and show the classical result that MLR dominance results in a threshold optimal detector.

Given the observation $x$ of a random variable, we wish to decide if $x$  is from pdf $f$ or $g$. To do this, we construct a decision policy  $\phi(x)$. The detector
decides
\beq  \begin{split}  f \quad  \text{ if }  \phi(x) & =  0  \\
        g \quad  \text{ if } \phi(x) & = 1  \end{split} \label{eq:detectornp} \eeq
The performance of the decision policy $\phi$  in (\ref{eq:detectornp}) is determined in terms of two metrics:
\begin{compactenum}
\item $\mathcal{P} =  \prob( \text{ reject } f | f \text{ is true }) $
\item $\mathcal{Q} =  \prob( \text{ reject } f | f \text{ is false }) $
\end{compactenum}
Clearly for the decision policy  $\phi(\cdot)$ in (\ref{eq:detectornp}),
$$ \mathcal{P} = \int_\reals f(x) \phi(x) dx,  \quad  \mathcal{Q} = \int_\reals g(x) \phi(x) dx. $$
The  well known Neyman-Pearson detector seeks to determine the  optimal decision policy  $\phi^*$ that maximizes $\mathcal{Q}$ subject to
the constraint  $\mathcal{P} \leq  \alpha$ for some user specified $\alpha \in (0,1]$.
The main result is
\begin{theorem-non}[Neyman-Pearson lemma]
Amongst all decision rules $\phi$ such that $\mathcal{P} \leq \alpha$,  the decision rule $\phi^*$  which maximizes $\mathcal{Q}$ is given by
$$ \phi^*(x) =  \begin{cases} 0   & \frac{f(x)}{g(x)} \geq c \\
		1 &  \frac{f(x)}{g(x)} < c  \end{cases} $$
where $c$ is chosen so that  $\mathcal{P}  = \alpha$.
\end{theorem-non}
\begin{proof}
Clearly for any $x \in \reals$,
$$ \big(\phi^*(x) - \phi(x) \big)  \big( c g(x) - f(x) \big) \geq 0. $$
Please verify the above inequality by showing that if $\phi^*(x) =1$ then both the terms in the above product are nonnegative; while if $\phi^*(x) = 0$, then
both the terms are nonpositive.
Therefore,
$$
c \bigg( \int  \phi^*(x) g(x) dx - \int \phi(x) g(x) dx\bigg)  \geq \int \phi^*(x) f(x) dx - \int \phi(x) f(x) dx 
$$
The right hand side is non-negative since by construction $ \int  \phi^*(x) f(x) dx = \alpha$ , while $ \int  \phi(x) f(x) dx \leq  \alpha$.
\end{proof}

{\bf Threshold structure of optimal detector.}  \index{Neyman-Pearson detector! optimal threshold structure}
Let us now give conditions so that the optimal Neyman-Pearson decision policy is a threshold policy:
Suppose now that $f$ MLR dominates $g$, that is $f(x)/g(x) \uparrow x$. 
Then clearly
\beq   \phi^*(x) = \begin{cases} 0  & x \geq x^*\\
					1  & x < x^*  \end{cases}   \label{eq:npt}\eeq
					where  threshold $x^*$ satisfies
					$$ \int_{-\infty}^{x^*} f(x) dx = \alpha $$
Thus if $f \gr g$, then  the optimal detector  (in the Neyman-Pearson sense) is  the threshold detector  (\ref{eq:npt}).

\end{compactenum}

\chapter{Monotonicity of Value Function for POMDPs}  

\begin{compactenum}

\item Theorem \ref{thm:pomdpmonotoneval} is the main result of the chapter and it gives conditions under which the value function
of a POMDP is MLR decreasing. Condition \ref{A1} was the main assumption on the possibly non-linear cost.
Give sufficient conditions for a quadratic cost $1-\pi^\p \pi + c_u^\p \pi$ to satisfy \ref{A1}. Under  what conditions does the entropy
$-\sum_i \pi(i) \log \pi(i) +  c_u^\p \pi$  satisfy \ref{A1}.

\item The shifted likelihood ratio order is a stronger order than the MLR order. Indeed,  $p> q$ in the shifted likelihood ratio order sense
if $p_i/q_{i+j}$ is increasing in $i$ for any $j$. If $j=0$ it coincides with the standard MLR order. 
(Recall also the problem in the previous chapter which says that  the shifted likelihood ratio order is closed under convolution.)
By using the shifted likelihood  ratio order, what further results on the value function
$V(\pi)$ can one get by using Theorem \ref{thm:pomdpmonotoneval}. 

\item  Theorem \ref{thm:pomdp2state} gives sufficient conditions for a 2-state POMDP to have a threshold policy. We have assumed that 
the observation probabilities are not action dependent.  How should the assumptions and proof be modified to allow for action
dependent observation probabilities?

\item How can  Theorem \ref{thm:pomdp2state} be modified if dynamic risk measures of \secn \ref{sec:riskaverse} are considered?
(see also  \secn \ref{sec:risk}).

\item {\bf Finite dimensional characterization of Gittins index for POMDP bandit} \cite{KW09}: \secn \ref{sec:POMDPbandit} dealt with
POMDP multi-armed bandit problem.
Consider a POMDP bandit where the Gittins index (\ref{eq:gittindef}) is characterized as the solution of Bellman's equation
(\ref{eq:bellmanbandit}).
 Since the value function of a POMDP is piecewise linear and concave (and therefore a finite dimensional characterization), it follows that 
a value iteration algorithm for (\ref{eq:bellmanbandit}) that characterizes the Gittins index also has a finite dimensional characterization.  Obtain an expression for 
this finite dimensional characterization for the Gittins index  (\ref{eq:gittindef})  for a horizon $N$ value iteration algorithm.

\item \secn \ref{sec:POMDPbandit} of the book deals with structural results for POMDP bandits. Consider the problem where several searchers are looking
for a stationary target.  Only one searcher can operate  at a given time and the searchers cannot receive state estimate information from other searchers or a base-station. The base station simply sends a  0 or 1 signal to each searcher telling them when to operate and when to shut down.
When it operates, the searcher obtains moves according to a Markov chain and obtains noisy information about the target.
Show how the problem can be formulated as a POMDP multi-armed bandit.

Show how a radar seeking to hide its emissions (low probability of intercept radar) can be formulated approximately as a POMDP bandit.

\item How does the structural result for the Gittins index for a POMDP bandit specialize to that of a full observed  Markov decision
process bandit problem?

\item Consider Problem \vref{prob:adaptive control} of Chapter \ref{ch:pomdpbasic} where optimal adaptive control of a fully observed MDP was 
formulated as a POMDP.  Give conditions that ensure that the value function $J_k(i,\belief)$  is MLR decreasing in $\belief$ and also monotone in $i$. What are the implications of this monotonicity in terms of dual control (i.e., exploration vs exploitation)?

\item {\bf Optimality of Threshold Policy for 2-state POMDP} \index{POMDP! optimality of threshold policy}
Recall that Theorem \ref{thm:pomdp2state} in the book gave sufficient conditions for the optimal policy of a 2-state POMDP to be a threshold.
Consider the proof of  Theorem \ref{thm:pomdp2state} in Appendix 11.A of the book. The last step  involved going from (\ref{eq:subex}) 
to a simpler expression via tedious but elementary steps. Here we specify what these steps are.

Start with (\ref{eq:subex})  in the book:
\begin{align}
\begin{aligned}\label{eq:i3}
I_3 &= \left[\sigma(\bar{\pi},y,2) + \sigma(\bar{\pi},y,1)\cfrac{T(\bar{\pi},y,1) - T(\pi,y,2)}{T(\pi,y,2) - T(\bar{\pi},y,2)} + \sigma(\pi,y,1)\cfrac{T(\pi,y,2) - T(\pi,y,1)}{T(\pi,y,2) - T(\bar{\pi},y,2)}\right]\\
&=\cfrac{I_{31} + I_{32} + I_{33}}{\sigma(\pi,y,2)\left(T(\pi, y, 2) - T(\bar{\pi}, y, 2)\right)}\\
I_{31} &= \sigma(\pi,y,2)\sigma(\bar{\pi},y,1)\left(T(\bar{\pi},y,1) - T(\pi,y,2)\right)\\ 
I_{32} &= \sigma(\pi,y,2)\sigma(\bar{\pi},y,2)\left(T(\pi,y,2) - T(\bar{\pi},y,2)\right)\\ 
I_{33} &= \sigma(\pi,y,2)\sigma(\pi,y,1)\left(T(\pi,y,2) - T(\pi,y,1)\right) 
\end{aligned}
\end{align}

The second element of HMM predictors $P(a)^\prime\pi$ and $(P(a)^\prime\bar{\pi})$ are denoted by $b_{a2}$, $b_{a1}, a = 1, 2$ respectively. Here $b_{a2}$ is defined as follows
\begin{align}\label{eq:sim_fil}
b_{a2} = (1-\pi(2))P_{12}(a) + \pi(2)P_{22}(a). 
\end{align}

Consider the following simplification of the term $I_{31}$ by using $b_{a2}$ and $b_{a1}$.
\begin{align}\label{eq:i31}
\begin{aligned}
I_{31} = &(B_{1y}(1-b_{22}) + B_{2y}b_{22})B_{2y}b_{11} - (B_{1y}(1-b_{11}) + B_{2y}b_{11})B_{2y}b_{22}\\
       = &B_{1y}B_{2y}(b_{11} - b_{22}) 
\end{aligned}
\end{align}
Similarly, $I_{32}$ and $I_{33}$ are simplified as follows
\begin{align}\label{eq:i3233}
\begin{aligned}
I_{32} = B_{1y}B_{2y}(b_{22} - b_{21}), I_{33} = B_{1y}B_{2y}(b_{22} - b_{12}) 
\end{aligned}
\end{align}
Substituting \eqref{eq:i31}, \eqref{eq:i3233} in \eqref{eq:i3} yields the following
\begin{align}\label{eq:i3f}
\begin{aligned}
I_3 &=  B_{1y}B_{2y}\cfrac{b_{11} + b_{22} - b_{21} - b_{12}}{\sigma(\pi,y,2)\left(T(\pi, y, 2) - T(\bar{\pi}, y, 2)\right)}
\end{aligned}
\end{align}
Substituting \eqref{eq:sim_fil} for $b_{ij}$ and some trivial algebraic manipulations yield the following
\begin{align}\label{eq:i3f}
\begin{aligned}
I_3 &=  B_{1y}B_{2y}(\pi(2) - \bar{\pi}(2))\cfrac{P_{22}(2) - P_{12}(2) - (P_{22}(1) - P_{12}(1))}{\sigma(\pi,y,2)\left(T(\pi, y, 2) - T(\bar{\pi}, y, 2)\right)}.
\end{aligned}
\end{align}

\item  Consider the following special case of a POMDP. Suppose the prior belief  $\belief_0 \in \Belief $ is known. From time 1 onwards,
the state is fully observed. How can the structural results in this chapter be used to characterize the optimal policy?
\end{compactenum}

\chapter{Structural Results for Stopping Time POMDPs}  

\section{Problems}
Most results in stopping time POMDPs in the literature  use the fact that the stopping set is convex
(namely, Theorem  \ref{thm:pomdpconvex}). 
Recall that the only requirements of 
Theorem  \ref{thm:pomdpconvex}  
are that the value function is convex and the stopping cost is linear.
Another important result for finite horizon POMDP stopping time  problems is the nested 
stopping set property  $\stopset_0 \subseteq \stopset_1 \subseteq \stopset_2 \ldots$. The following exercises discuss both these aspects.

\begin{compactenum}
\item {\bf Nested stopping set structure.}  \index{stopping time POMDP! nested stopping sets}
Consider the stopping time POMDP dynamic programming equation
$$ V(\belief) = \min\{c_1^\p \belief,  c_2^\p \belief + \sum_y V( \filter(\belief,y,u) ) \filterd(\belief,y,u) \}.$$
Define the stopping set as
$$ \stopset = \{ \belief: c_1^\p \belief \leq   c_2^\p \belief + \sum_y V( \filter(\belief,y,u) ) \filterd(\belief,y,u) \} = \{\belief: \policy^*(\belief) = 1 \text{ (stop) } \}$$
Recall the value iteration algorithm is
$$ V_{n+1}(\belief) = \min\{c_1^\p \belief,  c_2^\p \belief + \sum_y V_n( \filter(\belief,y,u) ) \filterd(\belief,y,u) \}, \quad V_0(\belief) = 0. $$
Define the stopping sets $\stopset_n = \{\belief: c_1^\p \belief \leq c_2^\p \belief + \sum_y V_n( \filter(\belief,y,u) ) \filterd(\belief,y,u) \}$.

Show that  the stopping sets satisfy $\stopset_0 \subseteq \stopset_1 \subseteq \stopset_2 \ldots$ implying that $$\stopset = \cup_{n} \stopset_n $$

\item {\bf Explicit characterization of stopping set.} Theorem \ref{thm:pomdpconvex} showed that for a stopping time POMDP, the stopping set $\stopset$ is convex.
By imposing further conditions, the set $\stopset$ can be determined explicitly.
Consider the following set of belief states
\beq \stopset^o= \{ \belief:  c_1^\p \belief \leq c_2^\p \belief + c_1^\p \tp^\p \belief \}  \label{eq:stopsetoo}  \eeq
Suppose the transition matrix $\tp$ and observation probabilities $\oprob$ of the stopping time POMDP satisfy the following property:
\beq \belief \in \stopset^o \implies \filter(\belief,\obs) \in \stopset^o, \quad  \forall  \obs \in \obspace . \label{eq:onesteppp}
\eeq

\begin{compactenum}
\item Prove that $\stopset^o = \stopset$. Therefore, the hyperplane $c_1^\p \belief = c_2^\p \belief + c_1^\p \tp^\p \belief $ determines the stopping set $\stopset$.

The proof proceeds in two steps: First prove by induction on the value iteration algorithm that for $\belief \in \stopset^o$, $V_n(\belief) = c_1^\p \belief$, for $n=1,2\ldots$.   \index{stopping time POMDP! characterization of stopping set}

Second, consider  a belief $\belief $ such that the optimal policy  goes one step and then stops.  This  implies that  the value function is $V(\belief) = c_2^\p \belief + c_1^\p \tp^\p \belief $. Therefore clearly   $c_2^\p \belief + c_1^\p \tp^\p \belief 
< c_1^\p \belief$. This implies that $\belief \notin \stopset^o$. So for any belief $\belief $ such that  $\optpolicy(\belief) $ goes one step and stops, then $\belief \notin \stopset^o$.       Therefore, for any belief $\belief$ such that $\optpolicy(\belief) $ goes more than one step and stops, then $\belief \notin \stopset^o$.

The two steps imply that  $\stopset^o = \stopset$.Therefore  that the stopping set is explicitly given by the polytope in (\ref{eq:stopsetoo}).

\item Give sufficient conditions on $\tp$ and $\oprob$ so that condition  (\ref{eq:onesteppp}) holds for a stopping time POMDP.
\end{compactenum}

\item 

Show that an identical, proof to Theorem  \ref{thm:pomdpconvex}
implies that the stopping sets $\stopset_n$, $n=1,2,\ldots$  are convex for a finite horizon problem.

\item {\bf Choosing a single sample from a HMM.} \index{multiple stopping problem}
Suppose a Markov chain $x_k$ is observed in noise sequentially over time as $\obs_k \sim \oprob_{\state_k,y}$, $k=1,2\ldots,N$.
 Over a horizon of length $N$, I need to choose a single  observation $y_k$ to maximize
$\E\{y_k\}$, $k \in 1,\ldots,N$. If at time $k$ I decide to choose observation $y_k$, then I get reward 
$\E\{y_k\}$ and the problem stops.
If I decide not to choose  observation $y_k$, then I can use it to update my estimate of the state and proceed to the next time
instant.  However, I am not allowed to choose $\obs_k$  at a later time.

\begin{compactenum}
\item Which single observation should I choose?

 Show that Bellman's equation becomes
$$ V_{n+1}(\belief) = \max_{u\in \{1,2\} }  \{ r^\p \belief ,  \sum_y V_n( \filter(\belief,  y)) \filterd(\belief,y) \} $$
where the elements of $r$ are $r(i) = \sum_{y} y \oprob_{iy}$, $i=1,\ldots,\statedim$.
Here $u=1$ denotes choose an observation, while $u=2$ denotes do not choose an observation.

\item Show using an identical proof to Theorem  \ref{thm:pomdpconvex} that the region of the belief space $\stopset_n = \{\belief: \mu^*(\policy) = 1\}$  is convex.
Moreover if (\ref{A2},\ref{A3}) hold, show that $e_1$ belongs to $\stopset_n$.  
Also show that $\stopset_0 \subseteq \stopset_1 \subseteq \stopset_2 \ldots$.

\item {\bf Optimal Channel sensing.}
Another interpretation of the above problem is as follows: The quality $x_k$ of a communication channel is observed in noise.  I need to
transmit a  packet using this channel.  If the channel is in state $x$, I incur a cost $c(x)$ for transmission. Given $N$ slots, when should
I transmit?  \index{optimal channel sensing}
\end{compactenum}

\item {\bf Optimal measurement selection for  a Hidden Markov Model (Multiple stopping problem)}. \index{multiple stopping problem}
 The following problem generalizes the previous problem as follows. \index{multiple stopping problem}
 I need to choose the best $L$ observations of a Hidden Markov model in a horizon of length $N$
where $L \leq N$? If I select observation $k$ then I get a reward $\E\{y_k\}$, if I reject the observation then I get no reward.
In either case, I use the observation $y_k$ to update my belief state. (This problem is also called the multiple stopping problem in 
\cite{Nak95}.)
Show that Bellman's dynamic programming recursion reads:
\begin{multline*}  V_{n+1}(\belief,l) = \max\{ r^\p \pi + \sum_\obs V_{n}(T(\belief,y),l-1) \filterd(\belief,y) , \\   \sum_\obs V_{n}(T(\belief,y),l) \filterd(\belief,y) \} , \quad n=1,\ldots,N \end{multline*}
with initial condition $V_n(\belief,0) = 0$, $n=0,1,\ldots$ and boundary conditions
$$V_{n}(\belief,n) = r^\p \pi + \sum_y  V_{n-1} ( T(\pi,y), n-1)  \sigma(\pi,Y)  , \quad n=1,\ldots,L. $$
The boundary condition says that if I  have only $n$ time points left to make $n$ observations, then I need to make an observation
at each of these $n$ time points. Obtain a structural result for the optimal measurement selection policy. (Notice that the actions
do no affect the evolution of the belief state $\belief$, they only affect $l$, so  the problem is simpler than a full blown POMDP.)

\item {\bf Separable POMDPs.}  \index{separable POMDPs} 
Recall that the action space is denoted as $\actionspace = \{1,2,\ldots,U\}$.
In analogy to \cite[Chapter 7.4]{HS84}, define a POMDP to be separable if:
the exists a subset $\bar{\actionspace}= \{1,2\ldots,\bar{U}\}$ of the action space $\actionspace$ such that for $\action \in \bar{\actionspace}$
\begin{compactenum}
\item The cost is additively separable: $c(x,u) = \phi(u) + g(x)$ for some scalars $\phi(u)$ and $g(x)$.
\item The transition matrix $\tp_{ij}(u)$ depends only on $j$. That is the process evolves independently of the previous state.
\end{compactenum}
Assuming that the actions $u \in \bar{\actionspace}$ are ordered so that $\phi(1) < \phi(2) < \ldots < \phi(\bar{U})$, clearly it is never
optimal to pick actions $2,\ldots,\bar{U}$. So solving the POMDP involves choosing between actions $
\{1, \bar{U}+1,\ldots, U\}$.
So from Theorem \ref{thm:pomdpconvex}, the set of beliefs where the optimal policy 
$\mu^*(\belief) = 1$ is convex.

Solving for the optimal policy for which the actions $\{\bar{U}+1,\ldots, U\}$ arise is still as complex as a solving a standard POMDP. However, the 
bounds proposed in Chapter \ref{chp:myopicul} can be used.  

Consider the special case of the above model  where $\bar{\actionspace} = \actionspace $ and instead of (a), $c(x,u)$ are arbitrary costs.  Then
show that the optimal policy is a  linear threshold policy.
\end{compactenum}

\section{Case Study: Bayesian Nash equilibrium  of one-shot global game for coordinated sensing} \index{coordinated sensing}
\index{Bayesian global game}
\index{Bayesian Nash equilibrium (BNE)}
 \index{game theory! global game}
This section gives a short description of Bayesian global games.
The ideas involve MLR dominance of posterior distributions and supermodularity and serves as a useful illustration of the 
structural results developed in the chapter.

We start with some perspective:
 Recall that in the classical Bayesian  social learning, agents act sequentially in time. The global games model that has been studied
in economics during the last two decades, considers multiple agents that act simultaneously by predicting the behavior
of other agents. 
The theory of global games was first introduced in \cite{CD93} as a tool for refining equilibria in economic game theory;
see \cite{MS00} for an excellent exposition.
Global games
represent a useful method for decentralized coordination amongst agents;   they 
have  been used to model speculative currency attacks and regime change in social systems,  see
\cite{MS00,KLM07,AHP07}. Applications in sensor networks and cognitive radio appear in \cite{Kri08,Kri09}.

\subsection{Global Game Model}
Consider a  continuum of agents  in which each agent $i$ obtains noisy measurements $Y^{(i)}$
of  an underlying state of nature $X$.  Here 
$$ Y^{(i)} = X + W^{(i)} , \quad X \sim \pi,\;
W^{(i)} \sim p_{W}(\cdot) $$

Assume all agents have the same noise distribution $p_W$.
 Based on its observation $y^{(i})$, each agent 
takes an action $u^i \in \{1, 2\}$ to optimize its expected reward \begin{equation} \label{eq:utilityi}
 R(X,\alpha,u=2) =  X + 
 f(\alpha), \quad
 R(X,u=1) =   0 
 \end{equation}
Here $\alpha\in [0,1]$ denotes the fraction of agents that choose action 2 and $f(\alpha)$ is a user specified function. We will call
$f$ the congestion function for reasons explained below.

As an illustrative
example,  suppose $x$  (state of nature) denotes the quality of a social group and  $y^{(i})$ denotes the measurement of this quality by agent $i$.  The action $u^i = 1$  means that agent $i$ decides not to  join the social group, while  $u^i = 2$ means that agent $i$ joins the group.
The utility function $R(u^i=2,\alpha)$ for joining the social group depends on $\alpha$, where $\alpha$ is the fraction of people who decide to join the  group. 
If $\alpha \approx 1$, i.e.,  too many people join the group, then  the utility to each agent is small since the group is too congested and agents do not receive sufficient individual service.
On the other hand, if $\alpha \approx 0$, i.e.,  too few people join the group, then the utility  is also small since there is not enough social interaction.
In this case the congestion function $f(\alpha)$ would be chosen as a quasi-concave function of $\alpha$ (that increases with $\alpha$ up to a certain value of $\alpha$ and then decreases with $\alpha$).

Since each agent is  rational, it  uses its  observation $y^{(i)}$ to predict $\alpha$, i.e., the fraction of  other agents 
that choose action 2. The main question is: {\em What is the optimal strategy for each agent $i$ to maximize its expected reward?}

\subsection{Bayesian Nash Equilibrium}
Let us now  formulate this problem:
Each agent
chooses its action $u \in \{1 , 2 \}$  based on a (possibly randomized) strategy $\mu^{(i)}$ that maps the current
observation $Y^{(i)}$ to the action $u$.  In a global game we are interested
in {\em symmetric strategies}, i.e.,  where all   choose the same
strategy denoted as  $\mu$. 
That is, each agent $i$  deploys the strategy $$ \mu:Y^{(i)} \rightarrow\{1  ,2 \} . $$
(Of course,  the action $\mu(Y^{(i)})$ picked by individual agents $i$ 
depend on their random observation $Y^{(i)}$. So the actions picked
are not necessarily identical even though the strategies are identical).

Let $\alpha(x)$ denote the fraction of agents that
select action $u=2$ (go)  given the quality of music $X = x$.
Since we are considering an infinite number of agents  that  behave independently, 
$\alpha(x)$ is also (with probability 1) the 
conditional probability that an agent  receives signal $Y^{(i)}$
and decides to pick $u=2$, given
$X$. So
\beq \label{eq:alfagdef}
\alpha(x) = P(\mu(Y) =2  | X=x) .
\eeq

We can now define the Bayesian Nash equilibrium (BNE)  of the global game.
For each agent $i$ 
given its observation $Y^{(i)}$, the goal  is to  choose  a strategy to optimize its local reward. That is, 
agent $i$  seeks to compute strategy $\mu^{(i),*}$ such that
\begin{equation} \label{eq:localpolicy}
\mu^{(i),*}(Y^{(i)})\in \{1 \text{ (stay) } ,2 \text{ (go) }\} \text{ maximizes }
\E [R(X,\alpha(X),\mu^{(i)}(Y^{(i)}) )|Y^{(i)} ].
 \end{equation}
Here $ R(X,\alpha(X),u )$ is defined as in (\ref{eq:utilityi}) with $\alpha(X)$
defined in (\ref{eq:alfagdef}). 

If such
a strategy $\mu^{(i),*}$ in (\ref{eq:localpolicy})
 exists and is the same for all agents  $i$,
  then they constitute a {\em symmetric} BNE  for the global game.
We will use the notation $\mu^*(Y)$ to denote this  symmetric BNE.

\noindent {\em Remark}: Since we are dealing with an incomplete information game, players
use randomized strategies. If  a BNE exists, then a pure (non-randomized) version exists straightforwardly
(see Proposition 8E.1, pp.225 in \cite{MWG95}). Indeed, with $y^{(i)}$ denoting realization
of random variable $Y^{(i)}$, 
$$ \E[R\big(X,\alpha(X),\mu(Y^{(i)}) \big)| Y^{(i)} = y^{(i)}] = \sum_{u=1}^2 \E[ R(X,\alpha(X),u)|Y^{(i)}= y^{(i)}] P(u|Y^{(i)}= y^{(i)}) .$$
Since a linear combination is maximized at its extreme values,
 the optimal (BNE) strategy is to choose $P(u^*|Y^{(i)} = y^{(i)}) = 1$ where 
\beq\label{eq:bang}  u^* = \mu^*(y^{(i)}) = \argmax_{u \in \{1,2\}} \E[R(X,\alpha(X),u)|Y^{(i)}= y^{(i)}] .\eeq

For notational convenience denote
$$ R(y,u) = \E[R(X,\alpha(X),u)|Y^{(i)}= y^{(i)}] $$

\subsection{Main Result. Monotone BNE} \index{monotone Bayesian Nash equilibrium}
With the above description, we will now give sufficient conditions for the BNE $\mu^*(y)$  to be monotone increasing in $y$ (denoted
$\mu^*(y) \uparrow y$).
This implies that the BNE is a threshold policy of the form:
$$ \mu^*(y) = \begin{cases} 1  & y \leq  y^* \\
					2  & y  > y^* \end{cases} $$
Before proving this monotone structure, first note that $\mu^*(y) \uparrow y$ implies that $\alpha(x)$ in (\ref{eq:alfagdef}) becomes
$$ \alpha(x) = P( y > y^* | X = x) = P(x + w > y^*) = P( w > y^* - x) = 1- \cdf_W(y^*-x)$$

Clearly from (\ref{eq:bang}), a sufficient condition for $\mu^*(y) \uparrow y$  is that
$$R(y,u ) = \int R(x, \alpha(x) ,u ) \,\pdf( x | y) dx $$
is supermodular in $(y,u)$ that is 
$$ R(y,u+1) - R(y,u) \uparrow y.$$
Since $R(X,u=0)$ it follows that $R(y,1) = 0$.  So it suffices that  $R(y,2) \uparrow y $.

\begin{compactenum}
\item What are sufficient conditions on the noise pdf $p_W(\cdot)$, and congestion function $f(\cdot)$ in
(\ref{eq:utilityi}) so that  $R(y,2) \uparrow y$ and so BNE $\mu^*(y) \uparrow y$?

Clearly  sufficient conditions for $R(y,2) \uparrow y$ are:
\begin{compactenum}
\item  $p(x|y)$ is MLR increasing in $y$,
\item $\R(x,\alpha(x),2)$ is increasing in $x$.
\end{compactenum}

But we know that $p(x|y)$ is MLR increasing in $y$ if the noise distribution is such  that $\pdf_W(y - x)$ is TP2 in $x,y$

Also $R(x,\alpha(x),2)$  is increasing in $x$ if its derivative wrt $x$ is positive. That is,
$$ \frac{d}{dx} R(x,\alpha(x),2) =  1 + \frac{df}{d\alpha} \frac{d \alpha} {dx}  = 1 +  \frac{df}{d\alpha}  p_W(y^*-x) >  0 $$

To summarize:  The BNE $\mu^*(y) \uparrow y$ if the following two conditions hold:
\begin{compactenum}
\item  $p(y|x) = p_W(y-x)$ is TP2 in $(x,y)$
\item  $$ \frac{df}{d\alpha}  > -\frac{1}{p_W(y^*-x)  } $$
\end{compactenum}
Note that a sufficient condition for the  second  condition is that 
$$  \frac{df}{d\alpha}  > -\frac{1}{\max_w p_W(w) } $$

\item Suppose $W$ is uniformly distributed in $[-1.1]$. Then  using the above conditions show that  a sufficient condition on the congestion function $f(\alpha)$ 
for the BNE to be monotone is that $df/d\alpha > - 2$.

\item Suppose  $W$ is zero mean Gaussian noise with variance $\sigma^2$. Then using the above conditions show that  a sufficient condition on the congestion function $f(\alpha)$ 
for the BNE to be monotone is that $df/d\alpha > - \sqrt{2 \pi} \sigma$.

\subsection{One-shot HMM Global Game}  \index{HMM global game}
Suppose that $X_0 \sim \belief_0$, and given $X_0$, $X_1$ is obtained by simulating from transition matrix $\tp$.  The observation for agent
$i$ is  obtained as the HMM observation
$$ Y^{(i)} =  X_1 + W^{(i)},  \quad   W^{(i)} \sim p_W(\cdot). $$ 
In analogy to the above
derivation, characterize the BNE of the resulting one-shot HMM global  game. (This will require assuming that $\tp$ is TP2.)

\end{compactenum}


\chapter{Stopping Time POMDPs for Quickest Change Detection}  

\begin{compactenum}

\item  For classical detection theory, a ``classic" book is the multi-volume \cite{Van68}.
\item As mentioned in the book, 
there are two approaches to quickest change detection:  Bayesian and minimax.
Chapter \ref{chp:stopapply} of the book deals with Bayesian quickest detection which assumes that the change point distribution is known (e.g.  phase distribution).  The focus of Chapter \ref{chp:stopapply}  was to determine the structure of the optimal policy of the 
Bayesian detector by showing that the problem is a special case of a stopping time POMDP.
\cite{VB13} uses nonlinear renewal theory to analyze the performance of the optimal Bayesian detector.

The minimax formulation for quickest detection assumes that the change point is either deterministic or has an unknown distribution.
For an excellent starting point on performance analysis of change detectors with  minimax formulations  please see \cite{TM10} and  \cite{PH08}.
The papers \cite{Lor71,Mou86}  gives a lucid description of the analysis of change detection in this framework.

\item {\bf  Shiryaev Detection Statistic.} 
In the classical Bayesian formulation of quickest detection described in \secn \ref{sec:convexstop},
a two state  Markov chain is considered to model geometric distributed change times.  Recall  (\ref{eq:tpqdp}), namely,
\beq \tp = \begin{bmatrix} 1 & 0 \\ 1- \tp_{22} & \tp_{22}  \end{bmatrix} , \;  \belief_0 = \begin{bmatrix} 0 \\ 1 \end{bmatrix} , \quad
\tau^0 = \inf\{ k:  \state_k = 1\}. \eeq
where $1-\tp_{22}$ is the parameter of the geometric prior.

In  classical detection theory,  the  belief state $\belief_k$  is written in terms of the  {\em Shiryaev detection
statistic} $r_k$ which is defined as  follows: \index{Shiryaev detection statistic}
\beq
 r_k \ole  \frac{1}{1-\tp_{22}}\times \frac{\belief_k(2) }{1 - \belief_k(2)}  \eeq
Clearly $r_k$ is an increasing function of $\belief_k(2)$ and so all the monotonicity results in the chapter continue to hold. 
In particular Corollary \ref{cor:qdclassical} in the book holds for $r_k$ implying a threshold policy in terms of $r_k$.

In terms of the Shiryaev 
statistic $r_k$, it is straightforward to write
the belief state update (HMM filter for 2  state Markov chain) as a function of the likelihood ratio as follows:
\begin{empheq}[box=\fbox]{equation}
r_k = \frac{1}{1-p}  \left( r_{k-1} + 1 \right)  L(y_k)  \label{eq:shirstat}
\end{empheq}
where
$$ p = 1 - \tp_{22}, \quad   L(y_k ) =  \frac{\oprob_{2y_k}}{ \oprob_{1y_k}}   \text{ (likelihood ratio) }$$
In (\ref{eq:shirstat}) by choosing $p \rightarrow 0$, the Shiryaev detection statistic converges to the so called {\em Shiryaev-Roberts detection statistic}.
Note that as $p \rightarrow 0$, the Markov chain becomes a slow Markov chain. We have analyzed in detail how to track the state of such a slow
Markov chain via a stochastic approximation algorithm in Chapter \ref{chp:markovtracksa} of the book. \index{Shiryaev-Roberts detection statistic}

The Shiryaev-Roberts detector for change detection reads:
\begin{compactenum}
\item Update the Shiryaev-Roberts statistic 
$$ r_k =  \left( r_{k-1} + 1 \right)  L(y_k) $$
\item If $r_k \geq r^*$  then stop and declare a change.
Here $r^*$ is a suitably chosen detection threshold.

\end{compactenum}
Please see \cite{PT12} for a nice survey description of minimax change detection and also the sense in which the above Shiryaev-Roberts detector
is optimal.

\item {\bf Classical Bayesian sequential detection.} 
\index{classical sequential detection}
 This problem shows that classical Bayesian sequential detection is a trivial case of the results
developed in Chapter \ref{ch:pomdpstop}.

Consider a random variable $\state \in \{1,2\}$.
So the transition matrix is $\tp = I$.  Given  noisy observations $\obs_k \sim\oprob_{\state y}$, the aim is to decide if the underlying state is either 1 or 2.  Taking stop action 1 declares that the state is 1 and stops. Taking stop action 2 declares that the state is 2 and stops. 
Taking action 3 at time $k$ simply takes another measurement $\obs_{k+1}$.
The misclassification costs are: $$\cost(\state=2,\action=1) = \cost(\state=1,\action=2) = L.$$ 
The cost of taking an additional measurement is $\cost(\state,\action=3) = C$. What is the optimal policy $\mu^*(\belief)$?

Since $\tp = I$, show that the dynamic programming equation reads
\begin{align*}
 V(\belief) &= \min\{ \belief_2 L, \;\belief_1 L, \; C + \sum_y V(T(\belief,y)) \filterd(\belief,y) \} \\
 T(\belief,y) &= \frac{\oprob_y \belief}{\one^\p \oprob_y \belief} , \quad  \filterd(\belief,\obs) = \one^\p \oprob_y \belief ,
 \end{align*}
where $\belief = [\belief(1), \belief(2)]^\p$ is the belief state. Note that $\obs \in \obspace$ where $\obspace$ can be finite or
 continuum (in which case $\sum$ denotes integration over $\obspace$).

From Theorem \ref{thm:pomdpconvex} we immediately know that the stopping sets $$\region_1 = \{\belief: \mu^*(\belief) = 1\}, \text{ and }
\region_2 =  \{\belief: \mu^*(\belief) = 2\}$$ are convex sets. Since the belief state is two dimensional, in terms of the second component $\belief(2)$,  $\region_1 $ and 
$\region_2$ are intervals in the unit interval $[0,1]$. Clearly $\belief(2) = 0 \in \region_1$ and $\belief(2)=1$ in $\region_2$.
Therefore $\region_1 = [0,\pi_1^*]$ and $\region_2 = [\pi_2^*,1]$ for some $\pi_1^* \leq \pi_2^*$. So the continue region is
$[\pi_1^*,\pi_2^*]$.  

Of course, Theorem \ref{thm:pomdpconvex} is much more general since it does not require $X=2$ states and $\state_k$ can evolve
according to a Markov chain with transition matrix $\tp$ (whereas in the simplistic setting above, $\state$ is a random variable).

\item {\bf Stochastic Ordering of Passage Times for Phase-Distribution.} \index{passage time}
In quickest detection, we formulated the 
change point  $\tau^0$ to have  a   
{\em  phase type (PH) distribution}. 
A systematic investigation of the statistical properties of PH-distributions can be found in \cite{Neu89}.
The family of all PH-distributions forms a dense subset for the set of all distributions
	\cite{Neu89} i.e., for any given distribution function $F$ such that $F(0) = 0$, one can find a sequence of PH-distributions	
$\{F_n , n	\geq	1\}$	 to		approximate	$F$	uniformly over $[0, \infty)$.
Thus PH-distributions  can be used to  approximate   change points with an arbitrary distribution. This is done by
	constructing a multi-state  Markov chain as follows:
Assume    state `1'  (corresponding to belief $e_1$) is an absorbing state
and denotes the state after the jump change.  The states $2,\ldots,X$ (corresponding to beliefs $e_2,\ldots,e_X$) can be viewed as  a single composite state that $x$ resides in before the jump. 
To avoid trivialities, assume that  the change occurs after at least one measurement. So the initial distribution $\belief_0$ satisfies $\belief_0(1) = 0$.
The 
transition probability matrix  is of the form
\beq \label{eq:phmatrix}
\tp = \begin{bmatrix}  1 & 0 \\ \underline{\tp}_{(X-1)\times 1} & \bar{\tp}_{(X-1)\times (X-1)} \end{bmatrix}.
\eeq
The  {\em first passage time} $\tau^0$ to state 1 denotes the time at which $x_k$ enters the absorbing state 1:
 \beq \tau^0 = \min\{k: x_k = 1\}.  \label{eq:tau}  \eeq
 As described in \secn \ref{sec:qdphform} of the book,
the distribution of $\tau^0$ is determined by choosing the transition probabilities $\underline{P}, \bar{P}$ in 
 (\ref{eq:phmatrix}).  
 The 
distribution of the absorption time to state 1 is denoted by  $$\nu_k = \prob(\tau^0 = k)$$ and given by
\beq \label{eq:nu}
 \nu_0 = \belief_0(1), \quad \nu_k = \bar{\belief}_0^\p \bar{P}^{k-1} \underline{P}, \quad k\geq 1, \eeq
 where $\bar{\belief}_0 = [\belief_0(2),\ldots,\belief_0(X)]^\p$.

{\bf Definition. Increasing Hazard  Rate}:  A   pmf  $p$ is said to be increasing hazard rate (IHR) if \index{increasing hazard rate}
$$ \frac{\bar{F}_{i+1}}{{\bar{F}_{i}}}  \downarrow i, \quad \text{ where }  \bar{F}_i = \sum_{j=i}^\infty p_j$$

{\bf Aim}. Show that if the transition matrix $\tp$ in (\ref{eq:phmatrix}) is TP2 and initial condition $\belief_0 = e_X$, then  the passage time distribution $\nu_k$ 
in (\ref{eq:nu})
satisfies the increasing hazard rate (IHR)
property; see \cite{Sha88} for a detailed proof.

\item {\bf Order book high frequency trading and social learning.}  \index{order book high frequency trading}
Agent based models for high frequency trading with an order book have been studied a lot recently \cite{AS08}. Agents trade (buy or sell) stocks by exploiting information about the decisions of previous agents (social learning) via an order book in addition to a private (noisy) signal they receive on the value of the stock. We are interested in the following: (1) Modeling the dynamics of these risk averse agents, (2) Sequential detection of a market shock based on the behavior of these agents.

The agents perform social learning according to the protocol in \secn \ref{sec:herdaa} of the book.  A market maker needs to decide based on the actions of the agents if there
is a sudden change (shock) in the underlying value of an asset. Assume that the shock occurs with a phase distributed change time.
The individual agents perform social learning with a CVaR social learning filter as in \secn \ref{sec:classicalsocial} of the book.  The market maker aims
to determine the shock as soon as possible.

Formulate this decision problem as a quickest detection problem. Simulate the value function and optimal policy. Compare it with the  
market maker's optimal policy
obtained when the agents perform risk
neutral social learning. See \cite{KB16} for details.
\end{compactenum}

\chapter{Myopic Policy Bounds for POMDPs and Sensitivity}  

\begin{compactenum}

\item To obtain upper and lower  bounds to the optimal policy, the key idea was to change the cost vector but still preserve the optimal policy. \cite{KP15} gives a complete description of this idea.
What if a nonlinear cost was subtracted from the costs thereby still keeping the optimal policy the same. Does that allow for larger regions of the belief space
where the upper and lower bounds coincide?
Is it possible to construct different transition matrices that yield the same optimal policy?

\item{\bf First order dominance of Markov chain sample paths.} In \secn \ref{sec:tpcp} of the book we defined the importance concept of copositive dominance
to say that if two transition matrices $\tp_1$ and $\tp_2$ satisfy $\tp_1 \lR \tp_2$   (see Definition  \ref{def:lR}), then the one step ahead predicted belief satisfies the MLR dominance property
$$ \tp_1^\p \belief \lr \tp_2 ^\p \belief. $$  
If we only want first  order stochastic dominance, then the following condition suffices: 
Let $U$ denote the $\statedim \times \statedim $ dimensional triangular  matrix with elements 
$U_{ij} = 0,  i > j$ and $U_{ij} = 1,  i \leq j$.
\begin{compactenum}
\item 
Show the following result:
$$ \tp_1 U \geq \tp_2 U \implies  \tp_1^\p \belief_1 \gs \tp_2^\p \belief_2\;  \text{ if } \belief_1 \gs \belief_2. $$
\item
 Consider the following special case of a POMDP. Suppose the prior belief  $\belief_0 \in \Belief $ is known. From time 1 onwards,
the state is fully observed. How can the structural results in this chapter be used to characterize the optimal policy?
\end{compactenum}

\item In \cite{Lov87} it is  assumed that one can construct a POMDP with observation matrices $\oprob(1),\oprob(2)$ such that 
(i) $\filter(\belief,\obs,2) \gr \filter(\belief,\obs,1)$ for each $\obs$ and (ii)  $\filterd(\belief,2) \gs \filterd(\belief,1) $.  Prove that it is impossible
to construct  an example that satisfies (i) and (ii) apart from the trivial case  where $\oprob(1) = \oprob(2)$. Therefore Theorem \ref{th:theorem1udit} does not apply
when the transition probabilities are the same and only the observation probabilities are action dependent. For such cases, Blackwell dominance is used.

\item {\bf Extensions of Blackwell dominance idea in POMDPs to more general cases.}
Blackwell dominance \index{Blackwell dominance} was used in \secn \ref{sec:blackwelldom}  of the book to construct myopic policies that bound the optimal policy of a POMDP.
But Blackwell dominance
is quite finicky. In  Theorem  \ref{thm:compare2} we assumed that the POMDP has  dependency structure $\state \rightarrow \yii \rightarrow \yi$. That is, the observation distributions are
$\oprob(2) = \pdf(\yii|x)$ and $\oprob(1) = \pdf(\yi|\yii)$.

\begin{compactenum}
\item 
Recall the proof of Theorem  \ref{thm:compare2} which is written element wise below for maximum clarity:
\[
\begin{split}
 \filter_j(\belief,\yi,1) = \frac{\sum_{\yii} \sum_i \belief(i) \tp_{ij}  \pdf(\yii|j,2)  \pdf(\yi | \yii)}
{\sum_m \sum_{\yii} \sum_i \belief(i) \tp_{im}  \pdf(\yii|m,2)\pdf(\yi | \yii) } \\
= 
  \frac{\sum_{\yii}   \dfrac{ \sum_i \belief(i) \tp_{ij}     \pdf(\yii|j,2) }
{ \cancel{ \sum_i \sum_m \belief(i) \tp_{im}     \pdf(\yii|m,2) }}   \cancel{ \sum_i \sum_m \belief(i) \tp_{im}     \pdf(\yii|m,\action)} \,
    \pdf(\yi | \yii)}
{\sum_m \sum_{\yii} \sum_i \belief(i) \tp_{im} \pdf(\yii|m,2) \pdf(\yi | \yii) }  \\ =
 \frac{\sum_{\yii} \filter_j(\belief,\yii,2)   \filterd(\belief,\yii,2) \pdf(\yi| \yii)}{  \sum_{\yii} \filterd(\belief,\yii,2) \pdf(\yi| \yii)}
\end{split}
\]

Then clearly $ \frac{ \filterd(\belief,\yii,2) \pdf(\yi| \yii)}{  \sum_{\yii} \filterd(\belief,\yii,2) \pdf(\yi| \yii)}$
is a probability measure w.r.t $\yii$.

\item Consider now the more general POMDP where $\pdf(\yi| \yii,\state)$ depends on the state $\state$. (In Theorem  \ref{thm:compare2}
this was functionally independent of $\state$.)
Then 
$$  \filter_j(\belief,\yi,1)  = \frac{\sum_{\yii} \filter_j(\belief,\yii,2)   \filterd(\belief,\yii,2) \pdf(\yi| \yii,j)}{  \sum_{\yii} \sum_m\filterd(\belief,\yii,2) \pdf(\yi| \yii,m)} $$
Now 
$$ \frac{   \filterd(\belief,\yii,2) \pdf(\yi| \yii,j)}{  \sum_{\yii} \sum_m\filterd(\belief,\yii,2) \pdf(\yi| \yii,m)} $$ is no longer a probability measure w.r.t. $\yii$.
The proof of Theorem  \ref{thm:compare2} no longer holds.

\item Next consider the case where the observation distribution is $\pdf(\yii_k | \state_k, \state_{k-1})$ and $\pdf(\yi | \yii)$. Then the proof
of  Theorem  \ref{thm:compare2} continues to hold.

\end{compactenum}

\item {\bf Blackwell dominance implies higher channel capacity.}  Show that if $B(1)$ Blackwell dominates $ B(2) $, i.e., $B(2) = B(1) Q$ for some stochastic matrix $Q$, then the capacity of a channel
with likelihood  probabilities given by $B(1)$ is higher than that with likelihood probabilities $B(2)$. 

\item 
Recall that the structural result involving Blackwell dominance deals with action dependent observation probabilities but assumes identical transition matrices for the various actions.
Show that copositive dominance and Blackwell dominance can be combined to deal with a POMDP with action dependent transition and
observation probabilities of the form:
\\  Action $\action = 1$:   $\tp^2, \oprob$.
\\  Action $\action=2$: $\tp, \oprob^2$.
\\
Give numerical examples of POMDPs with the above structure.
\end{compactenum}

\newpage
\chapter{Part IV.  Stochastic Approximation  and Reinforcement Learning}

Here we present three case studies of stochastic approximation algorithms.
The first case study  deals with online HMM parameter estimation and extends the method described in Chapter \ref{chp:markovtracksa}.
The second case study  deals with reinforcement learning of  equilibria in repeated games.  The third case study  deals with discrete stochastic optimization
(recall \secn \ref{sec:dopt} gave two algorithms)
and provides a simple example of such an algorithm.

\section{Case Study. Online HMM parameter estimation}  \index{online HMM estimation} 
Recall from Chapter \ref{chp:markovtracksa} that estimating the parameters of a HMM in real time is motivated by adaptive control of a POMDP.
The parameter estimation algorithm can be used to estimate the parameters of the POMDP
for a fixed policy; then the policy can be updated using dynamic programming (or approximation)  based on the parameters and so on.

This case study outlines several algorithms for recursive estimation of HMM parameters. The reader should implement these 
algorithms in Matlab to get a good feel for how they work.

Consider  the loss function for $N$ data points of a HMM or Gaussian state space model:
\beq  J_N(\th) = \E\{ \sum_{k=1}^N c_\th(x_k,y_k,\belief_k^\model)  \}  \label{eq:lossfunction} \eeq
where $x_k$ denotes the state, $y_k$ denotes the observation,  $\belief_k^\th$ denotes the belief state,  and $\th$ denotes the model variable.

The aim is to determine the model $\th$ that maximizes this loss function.

An offline gradient  algorithm operates iteratively to minimize this loss as follows:
\beq  \th^{(l+1)}  = \th^{(l)}  -  \epsilon  \nabla_\th J_N(\th)\big\vert_{\th = \th^{(l)}} \label{eq:offlinegrad} \eeq
The  notation $ \vert_{\th = \th^{(l)}}$ above means that  the derivatives are evaluated at $\th = \th^{(l)}$.

An offline Newton type algorithm operates iteratively as follows:
\beq  \th^{(l+1)}  = \th^{(l)}  -  \big[ \nabla_\th^2 J_N(\th)  \big]^{-1} \nabla_\th J_N(\th)\big\vert_{\th = \th^{(l)}} 
\label{eq:offlinenewton}
\eeq

\subsection{Recursive Gradient and Gauss-Newton Algorithms}
A 
recursive online version  of the above gradient  algorithm (\ref{eq:offlinegrad})  is 
\begin{empheq}[box=\fbox]{equation}  \begin{split}
 & \th_{k}  = \th_{k-1}  -  \epsilon \,    \nabla_\th c_\th(x_k,y_k,\belief^\th_k)\big\vert_{\th = \th_{k-1}} \, 
 \\ & \belief_k^{\th_{k-1}} =  \filter( \belief_{k-1}^{\th_{k-1}}, y_k;  \th_{k-1})   
 \end{split}  \label{eq:onlinegrad}
\end{empheq}
where $ \filter( \belief_{k-1}^{\th_{k-1}}, y_k;  \th_{k-1})$  is the optimal filtering recursion at time $k$  using  prior $\belief_{k-1}^{\th_{k-1}}$,
model $\th_{k-1}$ and observation $y_k$. The  notation $ \vert_{\th = \th_{k-1}}$ above means that  the derivatives are evaluated at $\th = \th_{k-1}$
Finally,  $\epsilon$ is a small positive step size.

The recursive Gauss Newton algorithm is an online implementation   of (\ref{eq:offlinenewton}) and reads
\begin{empheq}[box=\fbox]{equation} \begin{split}
 & \th_{k}  = \th_{k-1}  -  \mathcal{I}_k^{-1} \,    \nabla_\th c_\th(x_k,y_k,\belief_k^{\th})\big\vert_{\th = \th_{k-1}} \\
& \mathcal{I}_k  = \mathcal{I}_{k-1} + \epsilon \nabla^2  c_\th(x_k,y_k,\belief_k^{\th})\big\vert_{\th = \th_{k-1}}  \\
 &  \belief_k^{\th_{k-1}} =  \filter( \belief_{k-1}^{\th_{k-1}}, y_k;  \th_{k-1}) 
 \end{split} \label{eq:recgaussnewton}
\end{empheq}
Note that the above recursive Gauss Newton is a  stochastic approximation algorithm with a matrix step size $\mathcal{I}_k$.

\subsection{Justification  of (\ref{eq:onlinegrad})}
Before proceeding with examples, we give a 
 heuristic derivation of (\ref{eq:onlinegrad}). Write (\ref{eq:offlinegrad}) as 
$$ \th_{k}^{(k)}  = \th_k^{(k-1)}  -  \epsilon  \nabla_\th J_N(\th)\big\vert_{\th = \th_k^{(k-1)}} $$
Here the subscript $k$ denotes the estimate based on observations $\obs_{1:k}$. The superscript $(k)$ denotes the iteration of the offline optimization
algorithm.

Suppose that at each iteration $k$ we collect one more observation. Then the above algorithm becomes
\beq   \th_{k}^{(k)}  = \th_{k-1}^{(k-1)}  -  \epsilon  \nabla_\th J_k(\th)\big\vert_{\th = \th_{k-1}^{(k-1)}}  \label{eq:onoff}\eeq
Introduce the convenient notation
$$ \th_k = \th_{k}^{(k)} . $$
Next we  use the following two  crucial approximations:
\begin{compactitem}
\item First, make the inductive assumption that $\th_{k-1}$ minimized $J_{k-1}(\th)$ so that  
$$ \nabla_\th J_{k-1}(\th)\big\vert_{\th = \th_{k-1}}   = 0 $$ Then from  (\ref{eq:lossfunction}) it follows that
\beq \nabla_\th J_{k}(\th)\big\vert_{\th = \th_{k-1}} = \nabla_\th \E\{ c_\th(x_k,y_k,\belief^\th_k)  \}\big\vert_{\th = \th_{k-1}} 
\label{eq:approx1}
\eeq
\item  Note that evaluating the right hand side of (\ref{eq:approx1}) requires running a filter and its derivates wrt $\th$  from time 0 to $k$ for fixed model $\th_{k-1}$. We want a recursive
approximation for this.  It is here that the
second approximation is used. We revlaute the  filtering recursion  using
the a sequence of  available model estimates $\th_t$, $t=1,\ldots, k$ at each time $t$. In other words, we make the approximation
\beq   \belief_k^{\th_{k-1}} =  \filter( \belief_{k-1}^{\th_{k-1}}, y_k;  \th_k)   , k = 1,2,,\ldots,  \label{eq:approx2} \eeq

\end{compactitem}
To summarize, introducing  approximations (\ref{eq:approx1}) and (\ref{eq:approx2}) in (\ref{eq:onoff}) yields the online gradient algorithm (\ref{eq:onlinegrad}). The derivation of the Gauss-Newton algorithm is similar.

\subsection{Examples of online HMM estimation algorithm}
With the algorithms (\ref{eq:onlinegrad}) and (\ref{eq:recgaussnewton}) we can obtain several types of online HMM parameter estimators by choosing different loss functions $J$ in (\ref{eq:lossfunction}). Below we outline two popular choices.

\subsubsection*{1.  Recursive EM algorithm\footnote{This name is a misnomer. More accurately the algorithm below is a stochastic approximation algorithm that seeks to approximate the EM algorithm}}     \index{online HMM estimation! recursive EM} 
Recall from the EM algorithm, that the auxiliary likelihood for fixed parameter $\bar{\th}$ is 
$$Q_n(\th,\bar{\th}) =  \E\{ \log (\pdf(x_{0:n},y_{1:n}| \th) | y_{1:n}, \bar{\th})\} =   \E\{  \sum_{k=1}^n \log p_\th(x_k,y_k| x_{k-1}) | y_{1:n}, \bar{\th} \} $$
With $\th^o$ denoting the true model,  $\th$ denoting the model variable, and $\bmodel$ denoting a fixed model value, define
$$J_n(\th,\bmodel) = \E_{y_{1:n}}\{ Q_n(\th, \bar{\th}) | \th^o\}. $$
To be more specific, 
for a  HMM,  from (\ref{eq:emqfn}),
in the notation of (\ref{eq:lossfunction}),
\beq c_\th (y_k,\bmodel) =  
\sum_{i=1}^\statedim  \belief_{k|n}^{\bar{\th}}(i)   \; \log  \oprob^\model_{i\obs_k}
+  \sum_{i=1}^\statedim \sum_{j=1}^\statedim  \belief^{\bmodel}_{k|n} (i,j )
\log \tp^\model_{ij} . \label{eq:receminst}
\eeq
where $\tp^\model$ denotes the transition matrix and $\oprob^\model$ is the observation matrix
and $\bmodel$ is a fixed model for which the smoothed posterior $\belief^{\bmodel}_{k|n}$ is computed.

Note that $c_\th$ is a reward and not a loss; our aim is to maximize $J_n$.
The idea then is to implement a Gauss-Newton stochastic  gradient algorithm for maximizing $J_n(\th,\bmodel)$ for fixed model $\bmodel$,
then update $\bmodel$ and so on. This yields the following {\em recursive EM algorithm}:

\begin{compactenum}
\item  For $ k = n \Delta +1, \ldots (n+1) \Delta$ run
\beq \label{eq:recem}  
\begin{split}  \th_{k} &= \th_{k-1} + \mathcal{I}_k^{-1} 
\sum_i  \nablam c_\th (y_k, \bmodel_{n-1})  \belief_k^{\bmodel_{n}}(i) \\
 \mathcal{I}_k & = \mathcal{I}_{k-1} + \epsilon  \sum_i \nabla^2   c_\th(i,y_k,\belief_k^{\th})\big\vert_{\th = \th_{k-1}} \,   \belief_k^{\bmodel_{n}}(i) \\
   %
      \belief^{\bmodel_{n-1}}_{k} & = \filter(\belief^{\bmodel_{n-1}}_{k-1}, \obs_k; \th_{k-1})  \quad \text{ (HMM filter update) } 
\end{split}
\eeq
Here  $\belief_{k|n}^{\bar{\th}}$ and $\belief^{\bmodel}_{k|n} (i,j )$ in (\ref{eq:receminst}) are replaced by filtered estimates
$\belief_{k}^{\bar{\th}}$ and $\belief^{\bmodel}_{k-1} (i) \tp^{\bmodel}_{ij} \oprob^{\bmodel}_{j }\obs_{k}$.
\item Then update 
$  \bmodel_{n+1} =  \model_{(n+1)\Delta} $, set $n$ to $n+1$ and go to step 1.
\end{compactenum}

To ensure that the transition matrix estimates are a valid stochastic matrix, one can parametrize it in terms of spherical coordinates,
see~(\ref{eq:alfadef}).

As an illustrative  example, suppose we wish to estimate the $X$-dimension vector of state levels $\levels = (\levels(1),\levels(2),\ldots,\levels(X))^\p $
 of a HMM in zero mean Gaussian noise 
with known variance $\sigma^2$. Assume the transition matrix $\tp$ is known.
Then  $\th = g$ and 
$$  c_\th (y_k,\bmodel) 
= - \frac{1}{2 \sigma^2}  \sum_{i} \belief_k^{\bmodel}(i)  \big( y_k - \level(i) \big)^2 + \text{ constant }$$

  \subsubsection*{2. Recursive Prediction Error (RPE)}   \index{online HMM estimation! recursive prediction error} 
 Suppose $g$ is the vector of state levels of the underlying Markov chain and $\tp$ the transition matrix.  Then the model to estimate
 is  $\th = (g,\tp)$. 
 Offline prediction error methods seek to find the model $\th$ that minimizes the loss function
 $$ J_N(\th) = \E\{ \sum_{k=1}^N  (y_k - g^\p \belief_{k|k-1}) ^2\} $$  
 So   squared prediction error at each time $k$ is 
\beq c_\th(\state_k,\model_k,\belief^\model_k) = \big(y_k -  g^\p \tp^\p \belief^\th_{k-1}\big)^2  \label{eq:rpecost} \eeq
Note that unlike (\ref{eq:lossfunction}) there is no conditional expectation in the loss function.
 Note the key difference compared to the recursive EM. In the recursive 
EM $c_\th (x_k,y_k) $ is functionally independent of $\belief^\th$ and hence the recursive EM does not involve derivatives (sensitivity) of the HMM filter. In comparison,  the RPE  cost (\ref{eq:rpecost}) involves derivatives of  $\belief^\th_{k-1}$ with respect to $\th$.
  Then the derivatives with respect to $\th$  can be evaluated as in \secn \ref{sec:rmle}.

  \subsubsection*{3.  Recursive Maximum likelihood}  This was discussed in  \secn \ref{sec:rmle}.  The cost function is
$$c_\th(\state_k,\model_k,\belief^\model_k)  =\log \left[ \one^\p \oprob_{\obs_k}(\th) \belief_{k|k-1}^\th \right]
$$
  \index{online HMM estimation! recursive maximum likelihood}

Recursive versions of the method of moment estimation algorithm for the HMM parameters is presented in \cite{MKW15}.

\section{Case Study. Reinforcement Learning of Correlated Equilibria} \index{reinforcement learning of correlated equilibria}
This case study illustrates the use of stochastic approximation algorithms for learning the correlated equilibrium in a repeated
game. Recall in Chapter \ref{chp:markovtracksa} we used  the ordinary differential equation analysis of a stochastic approximation algorithm to 
characterize where it converges to. For a game, we will show that the stochastic approximation algorithm converges to
a differential inclusion (rather than a differential equation).  Differential inclusions are generalization of ordinary differential equations (ODEs) and arise naturally in game-theoretic learning, since the strategies according to which others play are unknown. Then by a straightforward
Lyapunov function type proof, we show that the differential inclusion converges to the set of correlated equilibria of the game, implying that the stochastic approximation algorithm also converges to the set of correlated equilibria.

\subsection{Finite Game Model}
Consider a finite action static  game\footnote{For notational convenience we assume two players with identical action spaces.} comprising two players $l=1,2$ with
costs
$c_l(\act1,\act2)$ where $\act1, \act2 \in \{1,\ldots,\actiondim\}$. Let $p$ and $q$ denote the randomized policies (strategies) of the two players:
$p(i) = \prob( \act1 = i) $ and $q(i) = \prob(\act2=i)$. So $p,q$ are $\actiondim$ dimensional probability vectors that live  in the $\actiondim-1$ dimensional unit simplex $\Pi$. Then the policies
$(p^*,q^*)$ constitute a Nash equilibrium if   the following inequalities hold:
\begin{empheq}[box=\fbox]{equation} 
\begin{split}
 \sum_{\act1,\act2} c_1(\act1,\act2)\, p^*(\act1) \, q^*(\act2) \leq   \sum_{\act2} c_1(u,\act2) \,  q^*(\act2)  , \quad u =1,\ldots,\actiondim \\
  \sum_{\act1,\act2} c_2(\act1,\act2) \, p^*(\act1) \,q^*(\act2) \leq   \sum_{\act1} c_2(\act1,u) \, p^*(\act1) , \quad u =1,\ldots,\actiondim.
 \end{split}  \label{eq:staticnashone}
 \end{empheq}
Equivalently, $(p^*,q^*)$ constitute a Nash equilibrium if  for all policies $p,q \in \Pi$,
\beq 
\begin{split}
 \sum_{\act1,\act2} c_1(\act1,\act2) \, p^*(\act1) \, q^*(\act2) \leq   \sum_{\act1,\act2} c_1(\act1,\act2)\,  p(\act1) \, q^*(\act2)   \\
  \sum_{\act1,\act2} c_2(\act1,\act2) \, p^*(\act1)\, q^*(\act2) \leq   \sum_{\act1,\act2} c_2(\act1,\act2) \, p^*(\act1) \, q(\act2) 
 \end{split}  \label{eq:staticnash}
 \eeq
The first inequality in (\ref{eq:staticnash}) says that if player 1 cheats and deploys policy $p$ instead of $p^*$, then it is worse off
and incurs an higher cost. The second inequality says that same thing for player 2. So in a non-cooperative game, since collusion is not allowed, there is no rational reason for any of the players to unilaterally  deviate from the Nash equilibrium
$p^*,q^*$.

By a standard application of Kakutani's fixed point theorem, it can be shown that for a finite action game, at least
one  Nash equilibrium always exists. However, computing it can be difficult since the above constraints are bilinear and therefore nonconvex.

\subsection{Correlated Equilibrium} The Nash equilibrium assumes that the player's act independently. The correlated equilibrium is a generalization of the Nash equilibrium. The two players now choose their action from the joint probability distribution
$\belief(\act1,\act2) $ where $$\belief(i,j) = \prob(\act1 = i,\act2=j).$$  \index{correlated equilibrium}
Hence the actions of the players are correlated. Then the policy $\belief^*$ is said to be a correlated equilibrium if
\begin{empheq}[box=\fbox]{equation} 
\begin{split}
 \sum_{\act2} c_1(\act1,\act2)\, \belief^*(\act1,\act2) \leq   \sum_{\act2} c_1(u,\act2)\, \belief^*(\act1,\act2)   \\
  \sum_{\act1} c_2(\act1,\act2) \belief^*(\act1,\act2) \leq   \sum_{\act1} c_2(\act1,u) \,\belief^*(\act1,\act2) 
 \end{split}  \label{eq:staticcoeq}
 \end{empheq}
 Define the set of correlated equilibria as
 \beq \label{eq:coreqset}
 \corr = \bigg\{\belief: \text{(\ref{eq:staticcoeq}) holds and $\belief(\act1,\act2)\geq 0, \sum_{\act1,\act2}
 \belief(\act1,\act2)=1$ } \bigg\} \eeq
{\em Remark}: In the special case where the players act independently, the correlated equilibrium specializes to a Nash equilibrium.
Independence implies the joint distribution $\belief^*(\act1,\act2)$ becomes the product of marginals:
so  $\belief^*(\act1,\act2)= p^*(\act1) q^*(\act2) $. Then clearly (\ref{eq:staticcoeq}) reduces to the definition (\ref{eq:staticnashone})
of a Nash equilibrium. Note that  the set of correlated equilibria specified by (\ref{eq:coreqset}) is a convex polytope in $\belief$.

\subsubsection{Why Correlated Equilibria?}  \index{game theory! correlated equilibrium}

John F. Nash proved in his famous paper~\cite{Nas51} that every game with a finite set of players and actions has at least one mixed strategy Nash equilibrium.
However, as asserted by Robert J. Aumann
\footnote{Robert J. Aumann was awarded the Nobel Memorial Prize in Economics in 2005 for his work on conflict and cooperation through game-theoretic analysis. He is the first to conduct a full-fledged formal analysis of the so-called infinitely repeated games.}
in the following extract from~\cite{Aum87}, ``Nash equilibrium does make sense if one starts by assuming that, for some specified reason, each player knows which strategies the other players are using.'' Evidently, this assumption is rather restrictive and, more importantly, is rarely true in any strategic interactive situation. He adds:
\begin{quote} {\em
``Far from being inconsistent with the Bayesian view of the world, the notion of equilibrium is an unavoidable consequence of that view. It turns out, though, that the appropriate equilibrium notion is not the ordinary mixed strategy equilibrium of Nash (1951), but the more general notion of correlated equilibrium.''} -- Robert J. Aumann
\end{quote}
This, indeed, is the very reason why correlated equilibrium~\cite{Aum87} best suits and is central to the analysis of strategic decision-making.

There is much to be said about correlated equilibrium; see Aumann~\cite{Aum87} for rationality arguments. Some advantages that make it ever more appealing include:
\begin{enumerate}
    \item \emph{Realistic:} Correlated equilibrium is realistic in multi-agent learning. Indeed, Hart and Mas-Colell observe in~\cite{HM01b} that for most simple adaptive procedures, ``\ldots there is a natural coordination device: the common history, observed by all players. It is thus reasonable to expect that, at the end, independence among players will not obtain;''

    \item \emph{Structural Simplicity:} The correlated equilibria set constitutes a compact convex polyhedron, whereas the Nash equilibria are isolated points at the extrema of this set~\cite{NCH04}. Indeed from (\ref{eq:coreqset}), the set of correlated equilibria is a convex polytope in $\belief$.

    \item \emph{Computational Simplicity:} Computing correlated equilibrium only requires solving a linear feasibility problem (linear program with null objective function) that can be done in polynomial time, whereas computing Nash equilibrium requires finding fixed points;

    \item \emph{Payoff Gains:} The coordination among agents in the correlated equilibrium can lead to potentially higher payoffs than if agents take their actions independently (as required by Nash equilibrium)~\cite{Aum87};

    \item \emph{Learning:} There is no natural process that is known to converge to a Nash equilibrium in a general non-cooperative game that is not essentially equivalent to exhaustive search. There are, however, natural processes that do converge to correlated equilibria (the so-called law of conservation of coordination~\cite{HM03}), e.g., \index{regret-matching procedure} regret-matching~\cite{HM00}.
\end{enumerate}

Existence of a centralized coordinating device neglects the distributed essence of social networks. Limited information at each agent about the strategies of others further complicates the process of computing correlated equilibria. In fact, even if agents could compute correlated equilibria, they would need a mechanism that facilitates coordinating on the same equilibrium state in the presence of multiple equilibria---each describing, for instance, a stable coordinated behavior of manufacturers on targeting influential nodes in the competitive diffusion process~\cite{TAM12}. This highlights the significance of adaptive learning algorithms \index{game theory! adaptive learning} that, through repeated interactive play and simple strategy adjustments by agents, ensure reaching correlated equilibrium. The most well-known of such algorithms, fictitious play, was first introduced in 1951~\cite{Rob51}, and is extensively treated in~\cite{FL98}. It, however, requires monitoring the behavior of all other agents that contradicts the information exchange structure in social networks. The focus below  is on the more recent regret-matching learning algorithms~\cite{BHS06,Cah04,HM00,HM01b}. 

Figure~\ref{fig:strategy-sets} illustrates how the various notions of equilibrium are related in terms of the relative size and inclusion in other equilibria sets. As discussed earlier in this subsection, dominant strategies and pure strategy Nash equilibria do not always exist---the game of ``Matching Pennies'' being a simple
example. Every finite game, however, has at least one mixed strategy Nash equilibrium. Therefore, the ``nonexistence critique'' does not apply to any notion that generalizes the mixed strategy Nash equilibrium in Figure~\ref{fig:strategy-sets}. A Hannan consistent strategy (also known as ``universally consistent'' strategies~\cite{FL95}) is one that ensures, no matter what other players do, the player's average payoff is asymptotically no worse than if she were to play any \emph{constant} strategy for in all previous periods. Hannan consistent strategies guarantee no asymptotic external regrets and lead to the so-called ``coarse correlated equilibrium''~\cite{MV78} notion that generalizes the Aumann's correlated equilibrium.

\begin{figure}[h]
\begin{tikzpicture}[decoration=zigzag]
  \filldraw[fill=gray!10, draw=gray!80] (0,0) ellipse (8.0cm and 4.0cm);
  \path[postaction={decorate,
    decoration={
      raise=-1em,
      text along path,
      text={|\tt|Hannan Consistent},
      text align=center,
    },
  }] (0,0) ++(110:8.0cm and 4.0cm) arc(110:70:8.0cm and 4.0cm);
  \filldraw[fill=gray!10, draw=gray!80] (0,0) ellipse (6.4cm and 3.2cm);
  \path[postaction={decorate,
    decoration={
      raise=-1em,
      text along path,
      text={|\tt|Correlated Equilibria},
      text align=center,
    },
  }] (0,0) ++(110:6.4cm and 3.2cm) arc(110:70:6.4cm and 3.2cm);
  \filldraw[fill=gray!20, draw=gray!80] (0,0) ellipse (4.8cm and 2.4cm);
  \path[postaction={decorate,
    decoration={
      raise=-1em,
      text along path,
      text={|\tt|Randomized  Nash},
      text align=center,
    },
  }] (0,0) ++(110:4.8cm and 2.4cm) arc(110:70:4.8cm and 2.4cm);
  \filldraw[fill=gray!30, draw=gray!80] (0,0) ellipse (3.2cm and 1.6cm);
  \path[postaction={decorate,
    decoration={
      raise=-1em,
      text along path,
      text={|\tt|Pure Nash},
      text align=center,
    },
  }] (0,0) ++(110:3.2cm and 1.6cm) arc(110:70:3.2cm and 1.6cm);
  \filldraw[fill=gray!40, draw=gray!80] (0,0) ellipse (1.6cm and 0.8cm) node
  {\texttt{Dominant Strategies}};
\end{tikzpicture}
\caption{Equilibrium notions in non-cooperative games. Enlarging the equilibria set weakens the behavioral sophistication on the player's part to distributively reach equilibrium through repeated plays of the game.}
\label{fig:strategy-sets}
\end{figure}

\subsection{Reinforcement Learning Algorithm}
To describe the learning algorithm and the concept of regret, it is convenient to deal with rewards rather than costs.
Each agent $l$ has utility reward  $\ut(u^{(l)},u^{-l})$ where $u^{(l)}$ denotes the action of agent $l$ and $u^{-l}$ denotes the action of the other agents.
The action space for each agent $l$ is $\{1,2,\ldots, U\}$.
Define the inertia parameter
\beq  \mu \geq   U \big( \max  \ut(\action, \action^{-l}) -  \min \ut(\action, \action^{-l}) \big)  \label{eq:inertia}\eeq

Each agent then runs the  regret matching Algorithm \ref{alg:rm}.  Algorithm \ref{alg:rm} assumes that  once a decision is made by an agent, it is observable by all other agents.  However, agent $l$ does not know the utility function of other agents. Therefore, a learning algorithms
such as Algorithm \ref{alg:rm} is required to learn the correlated equilibria.

The assumption that the actions of each agent are known to all other agents  can be relaxed; see \cite{HM01b} for "blind" algorithms that do not require this.  \index{regret matching algorithm}  \index{reinforcement learning of correlated equilibria! regret matching algorithm}

\begin{algorithm}
Each agent  $l$ with utility reward  $\ut(u^{(l)},u^{-l})$ independently executes the following:
\begin{compactenum}
    \item \bf Initialization: \rm
      Choose action $u_0^{(l)}\in \{1,\ldots,U\}$ arbitrarily.  
     Set  $\btheta^l_1=0.$
    \item
    Repeat for $n=1,2,\ldots$, the following steps:\\
   { \bf Choose Action:}  $u^{(l)}_{n} \in \{1,\ldots,U\}$ with probability
        \begin{align}
        \label{rmtrans} \prob(u^{(l)}_{n}=j| u_{n-1}^{(l)}=i,\btheta^l_n)
        &= \begin{cases}
       \frac{ |\btheta_n^l(i,j)|^+} {\mu } &  j\neq i, \\
        1-\sum_{m\neq i}\frac{ |\btheta_n^l(i,m)|^+}{\mu}  & j=i
\end{cases}
        \end{align}
where inertia parameter $\mu$ is defined in (\ref{eq:inertia}) and $|x|^+ \ole \max\{x,0\}.$

{\bf Regret  Update}: 
              Update the $\actiondim \times \actiondim$ regret matrix $\btheta^l_{n+1}$ as
        \begin{equation}
        \label{sa-d}
        \barray\ad
        \btheta^l_{n+1}(i,j) =\btheta^l_{n}(i,j)+\epsilon\left(I\{u^{(l)}_n=i\}\big(\ut(j,u^{-l}_n)-\ut(i,u^{-l}_n)\big) -\btheta^l_{n}(u,j)\right).
        \earray
        \end{equation}
        Here $\e \ll1$ denotes a constant positive step size.
\end{compactenum}
\caption{Regret Matching Algorithm for Learning Correlated Equilibrium} \label{alg:rm}
\end{algorithm}

\subsubsection{Discussion and Intuition of Algorithm \ref{alg:rm}}
{\em 1. Adaptive Behavior:} In~(\ref{sa-d}), $\stepsize$ serves as a forgetting factor to foster adaptivity to the evolution of the non-cooperative game parameters. That is, as agents repeatedly take actions, the effect of the old underlying parameters on their current decisions vanishes.

\noindent {\em  2. Inertia:} The choice of $\mu$ guarantees that there is always a positive probability of playing the same action as the last period. Therefore, $\mu$ can be viewed as an ``inertia'' parameter: A higher $\mu$
yields switching with lower probabilities. It plays a significant role in breaking away from bad cycles. It is worth emphasizing that the speed of convergence to the correlated equilibria set is closely related to this inertia parameter.

\noindent  {\em 3.  Better-reply vs. Best-reply:} In light of the above discussion, the most distinctive feature of the regret-matching procedure, that differentiates it from other works such as~\cite{FL99a}, is that it implements a better-reply rather than a best-reply strategy\footnote{This has the additional effect of making the behavior
continuous, without need for approximations~\cite{HM00}.}. This inertia assigns positive probabilities to any actions that are just better. Indeed, the behavior of a regret-matching decision maker is very far from that of a rational decision maker that makes optimal decisions given his (more or less well-formed) beliefs about the environment. Instead, it resembles the model of a reflex-oriented individual that reinforces decisions with ``pleasurable'' consequences~\cite{HM01b}.

We also point out the generality of Algorithm \ref{alg:rm}, by noting that it can be easily transformed into the well-known
{\em fictitious play} algorithm  by choosing $u^{(l)}_{n+1}=\arg\max_k \btheta^l_{n+1}(i,j)$ deterministically, where $u_n^{(l)}=i$,
and the extremely simple {\em best response} algorithm by further specifying $\epsilon=1$.


\noindent  {\em 4. Computational Cost:} The computational burden (in terms of calculations per iteration) of the regret-matching algorithm does not grow with the number of agents and is  hence scalable. At each iteration, each agent needs to execute two multiplications, two additions, one comparison and two table lookups (assuming random numbers are stored in a table) to calculate the next decision. Therefore, it is suitable for implementation in sensors with limited local computational capability.

\noindent {\em 5. Global performance metric}
Finally,  we introduce a metric for the global behavior of the system. 
The global behavior $\globbehav_\dtimee$ at time $k$ is defined as the empirical frequency of joint play of all agents up to period $k$. Formally,
\begin{equation}
\label{eq:global-behavior}
\globbehav_\dtimee = \sum_{\tau\leq k} (1-\stepsize)^{k-\tau} \unitvec_{\actprof_\tau}
\end{equation}
where $\unitvec_{\actprof_\tau}$ denotes the unit vector with the element corresponding to the joint play $\actprof_\tau$ being equal to one.  Given $\globbehav_\dtimee$, the average payoff accrued by each agent can be straightforwardly evaluated, hence the name global behavior. It is more convenient to define $\globbehav_\dtimee$ via the stochastic approximation recursion
\begin{equation}
\label{eq:global-behavior-SA}
\globbehav_\dtimee = \globbehav_{k-1} + \stepsize \lb \unitvec_{\actprof_\dtimee} - \globbehav_{k-1}\rb.
\end{equation}

The global behavior $\globbehav_\dtimee$ is a system ``diagnostic'' and is only used for the analysis of the emergent collective behavior of agents. That is, it does not need to be computed by individual agents.
In real-life application such as smart sensor networks, however, a network controller can monitor $\globbehav_\dtimee$ and use it to adjust agents' payoff functions to achieve the desired global behavior.

\subsection{Ordinary Differential Inclusion Analysis of Algorithm \ref{alg:rm}}  \index{ordinary differential inclusion analysis}
Recall from Chapter \ref{chp:markovtracksa} that the dynamics  of a stochastic approximation algorithm can be characterized by an
ordinary differential equation obtained by averaging the equations in the algorithm.  In particular,
using Theorem \ref{thm:odeweak} of Chapter \ref{chp:markovtracksa},  the estimates generated by the stochastic approximation
algorithm converge weakly to the averaged system corresponding to (\ref{sa-d})  and
(\ref{eq:global-behavior-SA}), namely,
\beq
\begin{split}
\frac{d R(i,j)}{dt} & =  \E_{\belief} \left\{ I(u_t= i ) \big(  \ut(j,u^{-l}) -   \ut(i,u^{-l}) \big) - R(i,j) \right\} \\
&=  \sum_{u^{-l}} \biggl[  \belief(i | u^{-l}) \,\bigg(  \ut(j,u^{-l}) -   \ut(i,u^{-l}) \bigg)   \biggr] \belief(u^{-l}) - R(i,j)  \\
\frac{dz}{dt} &= \belief(i|u^{-l})\, \belief(u^{-l}) - z 
\end{split} \label{eq:avgame1}
\eeq
where $\belief(u^{(l)}, u^{-l}) = \belief(u^{-l} |  u^{(l)}) \belief( u^{(l)} ) $ is the stationary distribution of the Markov process $(u^{(l)}, u^{-l})$.

Next note that the transition probabilities  in (\ref{rmtrans})  of $u_n^{(l)}$ given $R_n$  are conditionally independent of $u_n^{-l}$. So given
$R_n$,  $\belief(i | u^{-l}) = \belief(i)$.
So given the transition probabilities  in (\ref{rmtrans}), clearly the stationary distribution $\belief(u^{(l)})$ satisfies the linear algebraic equation
$$
\belief(i) = \belief(i) \biggl[ 1 - \sum_{j\neq i}  \frac{ | R(j,i)|^+ }{ \mu}  \biggr] + \sum_{j\neq i} \belief(j)  \frac{|R(i,j)|^+}{\mu} .
$$
which after cancelling out $\belief(i)$ on both sides yields
\beq
\sum_{i\neq j} \belief(i)  |R(i,j)|^+  = \sum_{i\neq j} \belief(j) |R(j,i)|^+  \label{eq:avgame2}
\eeq
Therefore the stationary distribution $\belief$ is functionally independent of the inertia parameter~$\mu$.

Finally note that as far as player $l$ is concerned, the strategy  $\belief(u^{-l})$ is not known. All is known is that $\belief(u^{-l})$ is a valid pmf.
So we can write
the  averaged dynamics of the regret matching Algorithm \ref{alg:rm} as
  \beq \displaystyle
  \begin{split}
& \begin{rcases} 
 \dfrac{d R(i,j)}{dt}   \in \displaystyle \sum_{u^{-l}} \biggl[  \belief(i ) \,\bigg(  \ut(j,u^{-l}) -   \ut(i,u^{-l}) \bigg)   \biggr] \belief(u^{-l}) - R(i,j)  
\\
 \dfrac{dz}{dt} \in  \belief(i)\, \belief(u^{-l}) - z  \end{rcases}  \; \belief(u^{-l}) \in  \text{ valid pmf} \\
&\sum_{i\neq j} \belief(i)  |R(i,j)|^+  = \sum_{i\neq j} \belief(j) |R(j,i)|^+ \end{split} \label{eq:odi}
\eeq
The above averaged dynamics 
 constitute an algebraically constrained ordinary differential inclusion.\footnote{Differential inclusions \index{differential inclusion} are a generalization of the concept of ordinary differential equations. A generic differential inclusion is of the form $dx/dt \in \mathcal{F}(x,t)$, where $\mathcal{F}(x,t)$ specifies a family of trajectories rather than a single trajectory as in the ordinary differential equations $dx/dt = F(x,t)$.}
We refer the reader to \cite{BHS05,BHS06} for an excellent exposition of the use of differential inclusions for analyzing game theoretical type learning algorithms.

{\em Remark}:
The asymptotics of a stochastic approximation algorithm is typically captured by an ordinary differential equation (ODE). Here, although agents observe $u^{-l}$, they are oblivious to the strategies $\belief(u^{-l})$ from which $u^{-l}$ has been drawn. Different strategies
$\belief(u^{-l})$ result in 
different trajectories of $R_n$. Therefore, $R_t$ and $z_t$ are  specified by a differential inclusions rather than  ODEs \index{ordinary differential equation}.

\subsection{Convergence of Algorithm \ref{alg:rm} to the set of correlated equilibria}
The previous subsection says that the regret matching Algorithm \ref{alg:rm} behaves asymptotically as an 
algebraically constrained differential inclusion (\ref{eq:odi}). So we only need to analyze the behavior of this differential inclusion to characterize the
behavior of the regret matching algorithm. 

\begin{theorem}
Suppose every agent follows the ``regret-matching''Algorithm~\ref{alg:rm}. Then  as  $t\to\infty$:
(i) $R(t) $ converges  to the negative orthant in the sense that
\begin{equation}
\dis\big[ R(t),\mathbb{R}^-\big] = \inf_{\boldsymbol{r} \in \reals^-}\big\| R(t)  - \boldsymbol{r}\big\|  \Rightarrow 0;
\end{equation}

(ii) $z(t) $ converges to the correlated equilibria set $\corr$ in the sense that
\begin{equation}
\label{eq:convergence-1}
\dis [ z(t),\corr ] = \inf_{\boldsymbol{z} \in \corr}\left\| z(t)  - \boldsymbol{z}\right\| \Rightarrow 0.
\end{equation}
\end{theorem}
The proof below shows the simplicity and elegance of the ordinary differential equation (inclusion) approach for analyzing stochastic approximation algorithm. Just a few elementary lines
based on the Lyapunov function yields the proof.

\begin{proof}

Define the Lyapunov function \index{Lyapunov function}
\begin{equation}
\lyap\big(R\big) = \frac{1}{2}\big( \textmd{dist}\big[R,\reals^-\big] \big)^2= \frac{1}{2}\sum_{i,j} \big(\big| R(i,j)\big|^+\big)^2.
\end{equation}
Evaluating the time-derivative  and substituting for $dR(i,j) / dt$ from  (\ref{eq:odi}) we obtain
\begin{align}
{d\over dt}\lyap\big(R\big) &= \sum_{i,j} \big| R(i,j)\big|^+ \cdot {d\over dt} R(i,j) \nonumber\\
& = \sum_{i,j} \big| r(i,j)\big|^+ \Big[ ( \ut\big(j,u^{-l}\big) - \ut\big(i,u^{-l}\big)) \belief(i) - R(i,j) \Big]\nonumber\\
& = \underbrace{\sum_{i,j} \big| R(i,j)\big|^+ \big( \ut \big(j,u^{-l}\big) - \ut\big(i,u^{-l}\big)\big) \belief(i)}_{= 0\;\text{from (\ref{eq:avgame2}) }} - \sum_{i,j} \big| R(i,j)\big|^+ R(i,j)\nonumber\\
& = -2  \lyap\big(R\big).
\end{align}
In the last equality we used
\begin{equation}
\sum_{i,j} \big| R(i,j)\big|^+ R(i,j) = \sum_{i,j} \big( \big| R(i,j)\big|^+ \big)^2 = 2 \lyap\big(R\big).
\end{equation}
This completes the proof of the first assertion, namely that Algorithm \ref{alg:rm} eventually generates regrets that are non-positive.

To prove the second assertion,  from Algorithm \ref{alg:rm}, the elements of the regret matrix are
\begin{align}
R_k(i,j) &= \stepsize \sum_{\tau \leq k} (1-\stepsize)^{k-\tau} \bigg[ \ut \big(j,u_\tau^{-l}\big) - \ut \big(u^{(l)}_\tau, u^{-l}_\tau)\bigg]
I(u_\tau^{(l)} = i) \nonumber\\
& = \sum_{u^{-l}} z  (i,u^{-l}) \big[  \ut(j,u^{-l})  - \ut(i, u^{-l})\big]  
\end{align}
where $z(i,u^{-l})$ denotes the empirical distribution of agent $l$ choosing action $i$ and the rest playing $u^{-l}$. 
On any convergent subsequence $\lbrace z_{\underline{k}}\rbrace_{\underline{k}\geq 0}\rightarrow {\pi}$, then
\begin{equation}
\label{eq:thrm-5-3}
\lim_{k\to\infty} R_k(i,j)  = \sum_{u^{-l}} \pi(i, u^{-l}) \big[ \ut (j,u^{-l}) - \ut(i,u^{-l})\big]
\end{equation}
where $\pi(i,u^{-l})$ denotes the probability of agent $l$ choosing action $i$ and the rest playing $u^{-l}$.
The first assertion of the theorem  proved that the regrets converge to non-positive values (negative orthant). Therefore (\ref {eq:thrm-5-3})
yields that 
$$  \sum_{u^{-l}} \pi(i, u^{-l}) \big[ \ut (j,u^{-l}) - \ut(i,u^{-l})\big] \leq 0 $$
implying that $\belief$ is a correlated equilibrium.
\end{proof}

\subsection{Extension to switched Markov games}   \index{reinforcement learning of correlated equilibria! switched Markov game}
Consider the case now where rewards $\ut(u^{(l)},u^{-l})$ evolve according to an unknown Markov chain $\th_n$. Such a time varying
game can result from utilities in a social network evolving with time or the number of players changing with time. 
The reward for agent $l$ is  now  $\ut(u^{(l)},u^{-l}, \th_n)$.   The aim is to track the set of correlated equilibria $\corr(\th_n)$; that is
use the regret matching algorithm \ref{alg:rm} so that agents eventually deploy strategies from $\corr(\th_n)$.
If $\th_n$ evolves with transition matrix $I + \epsilon^2 Q$ (where $Q$ is a generator), then it is on a slower time scale than the dynamics of the regret
matching Algorithm \ref{alg:rm}.  Then a more general proof in the spirit of  Theorem \ref{ode-lim} yields that the regret matching algorithm can track the time
varying correlated equilibrium set  $\corr(\th_n)$.   Moreover, in analogy to \secn \ref {sec:fastmc}, if the transition matrix for $\th_n$ is $I + \epsilon Q$, then the asymptotic dynamics
are given by a switched Markov differential inclusion, see   \cite{KMY08,NKY13}.

\section{Stochastic Search-Ruler Algorithm} \index{discrete stochastic optimization! stochastic ruler}

We discuss  two simple variants of  Algorithm \ref{alg:RS}  that
require  less restrictive conditions for convergence than condition (O). Assume $c_n(\th)$ are uniformly bounded for $\th \in \Theta$.
Neither of the algorithms given below are particularly novel; but they are useful from a pedagogical point of view.

 It is convenient to normalize the objective (\ref{eq:discobj}) as follows:
 Let $\alpha \leq c_n(\th) \leq \beta$ where $\alpha$ denotes a finite lower
bound and $\beta>0$ denotes a finite upper bound. Define the
normalized costs $m_n(\th)$ as
\begin{equation} \label{eq:mn}
m_n(\th) = \frac{\cost_n(\th) -\alpha}{\beta-\alpha} , \quad
\text{ where }  0 \leq m_n(\th) \leq 1. 
\end{equation}
Then the stochastic optimization problem (\ref{eq:discobj}) is equivalent to 
\beq
\th^* = 
\arg \min_{\th \in \Theta}m(\th) \text{ where }  m(\th) = \E\{m_n(\th)\}  
\label{eq:scale}
\eeq
since scaling the cost function does not affect the minimizing solution.
Recall  $\Theta = \{1,2,\ldots,S\}$.

Define  the loss function
\beq  Y_n(\th,u_n) =  I\left(m_n(\th)  - u_n\right)
\text{ where } I(x) = \begin{cases} 1 & \text{ if } x> 0  \\
                          0 & \text{ otherwise } \end{cases} \label{eq:yn}\eeq
Here
$u_n $ is a independent uniform random number in $[0,1]$. The uniform random number $u_n$ is a stochastic ruler against which the candidate
$m_n(\th)$ is measured.
The result was originally used
in devising stochastic ruler optimization algorithms \cite{AA01} -- although here we 
propose a more efficient algorithm than the stochastic ruler.
Applying Algorithm \ref{alg:RS} to the cost function $ \E\{Y_n(\th,u_n)\}$ defined 
in (\ref{eq:yn}) yields the following stochastic search-ruler algorithm:

\begin{algorithm}
\caption{Stochastic  Search-Ruler}
\label{alg:SR}
Identical to Algorithm \ref{alg:RS} with $c_n(\th_n)$ and
$c_n(\tilde{\theta}_n)$ replaced by $Y_n(\th_n,u_n)$ and
$Y_n(\tilde{\th}_n,\tilde{u}_n)$. Here $u_n$ and $\tilde{u}_n$ are 
independent uniform random numbers in $[0,1]$.
\end{algorithm}

Analogous to Theorem \ref{thm:RS} we have the following result:

\begin{theorem}  \label{thm:SR} Consider the discrete stochastic optimization problem  (\ref{eq:discobj}).
Then the Markov chain $\{\th_n\}$ generated by Algorithm \ref{alg:SR} has 
the following property for its stationary distribution $\belief_\infty$:
\beq \frac{\pi_\infty(\th^*)}{\pi_\infty(\th)} =
 \frac{m(\th)}{m(\th^*)}
\frac{ (1 - m(\th^*))}{ (1-m(\th)) } > 1 .  \label{eq:ratio}\eeq
\end{theorem}

 The theorem says that Algorithm \ref{alg:SR} is attracted to set the global minimizers $\globalopt$. It spends more time in $\globalopt$
than any other candidates. The restrictive condition (O) is not required for  Algorithm \ref{alg:SR} to be attracted to $\globalopt$.
 Theorem \ref{thm:SR}
gives an explicit representation of the discriminative power of the algorithm
between the optimizer $\th^*$ and any other candidate $\th$ in terms of the normalized
expected costs $m(\th)$ and $m(\th^*)$.
 Algorithm \ref{alg:SR} is more efficient than the stochastic ruler algorithm of \cite{And99} when the candidate samples are chosen with equal probability.
The stochastic ruler algorithm of \cite{And99} has asymptotic efficiency
$\pi(\th^*)/\pi(\th) = (1 - m(\th^*))/ (1-m(\th)) $.
So Algorithm \ref{alg:SR} has the additional  improvement in efficiency due to the 
additional multiplicative term $m(\th)/m(\th^*)$ in (\ref{eq:ratio}).

\noindent {\em Variance reduction using common random numbers}:
A more efficient implementation of Algorithm \ref{alg:SR} can be obtained
by using variance reduction based on common random numbers  (discussed in Appendix  \ref{sec:crn} of the book) as follows:
Since $u_n$ is uniformly distributed in $[0,1]$, so is $1-u_n$.
Similar to Theorem \ref{thm:SR}
it can be shown that 
the optimizer
 $\th^*$ is the minimizing solution of
the following stochastic optimization problem
$
\th^* = \arg \min_\th \E\{Z_n(\th,u_n)\} $
where 
\begin{equation}
Z_n(\th,u_n) = \frac{1}{2} \left[ Y_n(\th,u_n) + Y_n(\th,1-u_n)\right]
\label{eq:zn} \end{equation}
where the normalized sample cost $m_n(\th)$ is defined in (\ref{eq:scale}).
Applying Algorithm \ref{alg:SR} with $ Z_n(\th_n,u_n)$ and $Z_n(\tilde{\th}_n,u_n)$ replacing  $ Y_n(\th_n,u_n)$ and $Y_n(\tilde{\th}_n,u_n)$,
respectively,
yields the  variance reduced search-ruler algorithm.

In particular, since the indicator function $I(\cdot)$ in (\ref{eq:yn}) is a monotone function of its argument, it follows that 
$\var\{Z_n(\th,u_n)\} \leq \var\{Y_n(\th,u_n)\} $.
 As a result one would expect that the stochastic
optimization algorithm using $Z_n$ would converge faster.

\begin{proof}
We first show that 
 $\th^*$ defined  in (\ref{eq:scale}) is the minimizing solution of
the stochastic optimization problem
$
\th^* = \arg \min_\th \E\{Y_n(\th,u_n)\}$.
 Using the smoothing property of conditional expectations (\ref{eq:lieeng})
 yields \index{smoothing property of conditional expectation}
\begin{align*}
\E\{ I\left(m_n(\th)  - u_n\right) \} &=  
\E\{\E\{I\left(m_n(\th)  - u_n\right)| m_n(\th) \} \} \\
&\hspace{-1.3cm}= \E\{ \prob(u_n < m_n(\th))\} = \E\{m_n(\th)\} = m(\th)\end{align*}
 The second equality follows since expectation of an indicator function is probability,
 the third equality holds because $u_n$ is a uniform random number in [0,1] so
 that $\prob(u_n<a) = a$ for any $a$ in $[0,1]$.

Next we show that 
the state process $\{\th_n\}$ generated by  Algorithm \ref{alg:SR} is a homogeneous,
 aperiodic,  irreducible, 
Markov chain on the  state space $\Th$
with  transition probabilities
$$ \tp_{ij} = P(\th_n = j|\th_{n-1} = i) = \frac{1}{S-1}m(i)\bigl(1-m(j)\bigr) .
$$
That the process $\{\th_n\}$ is a homogeneous
aperiodic irreducible Markov chain follows from its construction
in Algorithm \ref{alg:SR} -- indeed $\th_n$ only depends probabilistically on $\th_{n-1}$.
From Algorithm \ref{alg:SR}, given candidate $i$ and its associated
cost $Y_n(i,u_n)$, candidate $j$ is accepted if its associated
cost $\tilde{Y}_n(j,\tilde{u}_n)$ is
smaller. So \begin{align*}
 \tp_{ij} &= \frac{1}{S-1} P(\tilde{Y}_n(j,\tilde{u}_n) < Y_n(i,u_n)) 
\\& 
=\frac{1}{S-1} P(m_n(j) < \tilde{u}_n)
P(m_n(i) > u_n)
\end{align*}
Finally,  for this transition matrix, it is easily verified that  
\beq \pi_\infty(\th) = \kappa (1 - m(\th)) \prod_{j\neq \th} m(j) \eeq
is the invariant distribution
 where $\kappa$ denotes a normalization constant.
 Hence 
$$\frac{\pi_\infty(\th^*)}{ \pi_\infty(\th)} =  \frac{m(\th)}{m(\th^*)}\frac{ (1 - m(\th^*))}{ (1-m(\th)) }
= \frac{ 1/m(\th^*)-1}{1/m(\th) - 1} > 1$$
since $m(\th^*)$ is the global minimum
and therefore $ m(\th^*) < m(\th)$ for $\th \in \Theta - \globalopt$.
\end{proof}

\newpage

\bibliographystyle{plain}
\bibliography{$HOME/styles/bib/vkm}

\printindex
\end{document}

%% file: myproblemcommands.tex
\newcommand{\secn}{\S}

\makeatletter
\def\nl#1#2{\begingroup
    #2%
    \def\@currentlabel{#2}%
    \phantomsection\label{#1}\endgroup
}
\makeatother

{\ensuremath{
\begin{empheq}[box=\fbox]{align}
{#1}
\end{empheq}
}}

\newtheorem{theorem}            {Theorem}[section]

\newcommand{\normal}{\mathbfit{N}}  

{
\begin{mdframed}
\par\noindent\textbf{#1:}\begin{rmfamily}\noindent}%
{\end{rmfamily}
\end{mdframed}
}

\newcommand{\dob}{\rho}  

\newcommand{\var}{\operatorname{Var}}

\newcommand{\pdf}{p}
\newcommand{\cdf}{F}
\newcommand{\prob}{\mathbb{P}}
\newcommand{\E}                 {\Bbb{E}}
\newcommand{\bE}{\bar{\E}}

\newcommand{\pref}{q}  

\newcommand{\inp}{u}

\newcommand{\obs}{y}

\newcommand{\snoise}{w}
\newcommand{\onoise}{v}
\newcommand{\statem}{A}

\newcommand{\snoisem}{\Gamma}

\newcommand{\obsm}{C}

\newcommand{\level}{g}

\newcommand{\onoisem}{D}
\newcommand{\inpm}{f}
\newcommand{\oinpm}{g}
\newcommand{\snoisecov}{Q}
\newcommand{\onoisecov}{R}

\newcommand{\state}{x}
\newcommand{\statespace}{\mathcal{X}}
\newcommand{\obspace}{\mathcal{Y}}
\newcommand{\statedim}{X}

\newcommand{\mc}{r}  

\newcommand{\fun}{\phi}

\newcommand{\identity}{I}

\newcommand{\dvar}[2]{\|{#1}-{#2}\|_{\text{\tiny{TV}}}}

\newcommand{\oprob}{B}
\newcommand{\tp}{P}
\newcommand{\utp}{\bar{\tp}}
\newcommand{\ltp}{\underline{\tp}}

\newcommand{\finaltime}{N}

\newcommand{\model}{\theta}

\newcommand{\logl}{\mathcal{L}_\finaltime}  

\newcommand{\nablam}{\nabla_\model}

\newcommand{\belief}{\pi}

\newcommand{\bbelief}{\bar{\pi}}
\newcommand{\ubelief}{{q}}

\newcommand{\lbelief}{\underline{\pi}}

\newcommand{\Belief}{\Pi(\statedim)}

\newcommand{\tbelief}{\pi^0}
\newcommand{\history}{\mathcal{H}}

\newcommand{\kalmancov}{\Sigma}

\newcommand{\ole}{\stackrel{\text{defn}}{=}}

\newcommand{\gr}{\geq_r}
\newcommand{\lr}{\leq_r}
\newcommand{\gs}{\geq_s}
\newcommand{\ls}{\leq_s}

\newcommand{\filterd}{\sigma}

\newcommand{\filter}{T}

\newcommand{\argmin}{\operatornamewithlimits{argmin}}
\newcommand{\argmax}{\operatornamewithlimits{argmax}}

\newcommand{\reals}{{\rm I\hspace{-.07cm}R}}
\newcommand{\R}{{\rm I\hspace{-.07cm}R}}

\newcommand{\beq}{\begin{equation}}
\newcommand{\eeq}{\end{equation}}
\newcommand{\nn}{\nonumber}

\renewcommand{\th}{\theta}
\newcommand{\Th}{\Theta}

\newcommand{\p}{\prime}

\newcommand{\one}{\mathbf{1}}

\newcommand{\samt}{\Delta}

\newcommand{\lR}{\preceq}

\newcommand{\levels}{g}

\newcommand{\bmodel}{\bar{\model}}

\newcommand{\cost}{c}

\newcommand{\yi}{y^{(1)}}
\newcommand{\yii}{y^{(2)}}

\newcommand{\action}{u}

\newcommand{\actionspace}{\,\mathcal{U}}
\newcommand{\actiondim}{U}

\newcommand{\discount}{\rho}

\newcommand{\region}{\mathcal{R}}

\def\Policy{{\boldsymbol{\policy}}}
\newcommand{\policy}{\mu}

\newcommand{\optpolicy}{\policy^*}
\newcommand{\bpolicy}{{\boldsymbol{\mu}}}

\newcommand{\uV}{\underline{V}}

\newcommand{\stopset}{\mathcal{S}}

\newcommand {\policyl} {\underline{\mu}}

\newcommand{\use}{n}
\newcommand{\tar}{z}

\newcommand{\Vb}{\bar{V}}

\newcommand{\stepsize}{\epsilon}

\newcommand{\globalopt}{\mathcal{G}}

\def\lb{\left[}
\def\rb{\right]}

\newcommand{\ad}{&\!\!\!\disp}

\newcommand{\barray}{\begin{array}{ll}}
\newcommand{\earray}{\end{array}}
\newcommand{\disp}{\displaystyle}

\newcommand{\e}{\varepsilon}










%% file: problems.bbl
\begin{thebibliography}{10}

\bibitem{AA01}
M.~Alrefaei and S.~Andradottir.
\newblock A modification of the stochastic ruler method for discrete stochastic
  optimization.
\newblock {\em European Journal of Operational Research}, 133:160--182, 2001.

\bibitem{AM79}
B.~D.~O. Anderson and J.~B. Moore.
\newblock {\em Optimal filtering}.
\newblock Prentice Hall, Englewood Cliffs, New Jersey, 1979.

\bibitem{And99}
S.~Andradottir.
\newblock Accelerating the convergence of random search methods for discrete
  stochastic optimization.
\newblock {\em ACM Transactions on Modelling and Computer Simulation},
  9(4):349--380, Oct. 1999.

\bibitem{AHP07}
G.M. Angeletos, C.~Hellwig, and A.~Pavan.
\newblock Dynamic global games of regime change: Learning, multiplicity, and
  the timing of attacks.
\newblock {\em Econometrica}, 75(3):711--756, 2007.

\bibitem{Aum87}
R.~J. Aumann.
\newblock Correlated equilibrium as an expression of {B}ayesian rationality.
\newblock {\em Econometrica}, 55(1):1--18, 1987.

\bibitem{AS08}
M.~Avellaneda and S.~Stoikov.
\newblock High-frequency trading in a limit order book.
\newblock {\em Quantitative Finance}, 8(3):217--224, Apr 2008.

\bibitem{BOL10}
B.~Bahrami, K.~Olsen, P.~Latham, A.~Roepstorff, G.~Rees, and C.~Frith.
\newblock Optimally interacting minds.
\newblock {\em Science}, 329(5995):1081--1085, 2010.

\bibitem{BO91}
T.~Basar and G.~J. Olsder.
\newblock {\em Dynamic Noncooperative Game Theory}.
\newblock SIAM Series in Classics in Applied Mathematics, 1991.

\bibitem{BHS05}
M.~Benaim, J.~Hofbauer, and S.~Sorin.
\newblock Stochastic approximations and differential inclusions.
\newblock {\em SIAM Journal on Control and Optimization}, 44(1):328--348, 2005.

\bibitem{BHS06}
M.~Benaim, J.~Hofbauer, and S.~Sorin.
\newblock Stochastic approximations and differential inclusions, {P}art {II}:
  Applications.
\newblock {\em Mathematics of Operations Research}, 31(3):673--695, 2006.

\bibitem{Cah04}
A.~Cahn.
\newblock General procedures leading to correlated equilibria.
\newblock {\em International Journal of Game Theory}, 33(1):21--40, Dec. 2004.

\bibitem{CD93}
H.~Carlsson and E.~van Damme.
\newblock Global games and equilibrium selection.
\newblock {\em Econometrica}, 61(5):989--1018, Sept. 1993.

\bibitem{POMDP}
A.~R. Cassandra.
\newblock Tony's {POMDP} page.
\newblock http://www.cs.brown.edu/research/ai/pomdp/index.html.

\bibitem{Cas95}
D.~Casta{\~n}on.
\newblock Optimal search strategies in dynamic hypothesis testing.
\newblock {\em Systems, Man and Cybernetics, IEEE Transactions on},
  25(7):1130--1138, 1995.

\bibitem{CR14a}
K.~Chen and S.~Ross.
\newblock An adaptive stochastic knapsack problem.
\newblock {\em European Journal of Operational Research}, 239(3):625--635,
  2014.

\bibitem{CAM02}
L.~Chen, P.~O. Arambel, and R.~K. Mehra.
\newblock Estimation {U}nder {U}nknown {C}orrelation: {C}ovariance
  {I}ntersection {R}evisited.
\newblock {\em IEEE Transactions on Automatic Control}, 47(11):1879--1882, 11
  2002.

\bibitem{Che59}
Herman Chernoff.
\newblock Sequential design of experiments.
\newblock {\em The Annals of Mathematical Statistics}, 30(3):755--770, 1959.

\bibitem{DLK00}
A.~Doucet, A.~Logothetis, and V.~Krishnamurthy.
\newblock Stochastic sampling algorithms for state estimation of jump {M}arkov
  linear systems.
\newblock {\em IEEE Transactions on Automatic Control}, 45(2):188--202, Feb.
  2000.

\bibitem{EAM95}
R.~J. Elliott, L.~Aggoun, and J.~B. Moore.
\newblock {\em Hidden {M}arkov Models -- Estimation and Control}.
\newblock Springer-Verlag, New York, 1995.

\bibitem{EK99}
R.~J. Elliott and V.~Krishnamurthy.
\newblock New finite dimensional filters for estimation of discrete-time linear
  {G}aussian models.
\newblock {\em IEEE Transactions on Automatic Control}, 44(5):938--951, May
  1999.

\bibitem{EE99}
J.~Evans and R.J. Evans.
\newblock Image-enhanced multiple model tracking.
\newblock {\em Automatica}, 35(11):1769--1786, 1999.

\bibitem{FK14}
M.~Fanaswala and V.~Krishnamurthy.
\newblock Syntactic models for trajectory constrained track-before-detect.
\newblock {\em IEEE Transactions on Signal Processing}, 62(23):6130--6142,
  2014.

\bibitem{FK13}
M.~Fanaswalla and V.~Krishnamurthy.
\newblock Detection of anomalous trajectory patterns in target tracking via
  stochastic context-free grammars and reciprocal process models.
\newblock {\em IEEE Journal on Selected Topics Signal Processing}, 7(1):76--90,
  Feb. 2013.

\bibitem{Fel65}
A.~A. Fel'dbaum.
\newblock {\em Optimal control systems}.
\newblock Academic Press, 1965.

\bibitem{FV12}
J.~Filar and K.~Vrieze.
\newblock {\em Competitive Markov decision processes}.
\newblock Springer Science \& Business Media, 2012.

\bibitem{FL98}
D.~Fudenberg and D.~K. Levine.
\newblock {\em The Theory of Learning in Games}.
\newblock MIT Press, 1998.

\bibitem{FL99a}
D.~Fudenberg and D.~K. Levine.
\newblock Conditional universal consistency.
\newblock {\em Games and Economic Behavior}, 29(1):104--130, Oct. 1999.

\bibitem{FL95}
D.~Fudenberg and D.K. Levine.
\newblock Consistency and cautious fictitious play.
\newblock {\em Journal of Economic Dynamics and Control}, 19(5-7):1065--1089,
  1995.

\bibitem{GSV05}
D.~Guo, S.~Shamai, and S.~Verd{\'u}.
\newblock Mutual information and minimum mean-square error in {G}aussian
  channels.
\newblock {\em IEEE Transactions on Information Theory}, 51(4):1261--1282,
  2005.

\bibitem{HM00}
S.~Hart and A.~Mas-Colell.
\newblock A simple adaptive procedure leading to correlated equilibrium.
\newblock {\em Econometrica}, 68(5):1127--1150, 2000.

\bibitem{HM01b}
S.~Hart and A.~Mas-Colell.
\newblock A reinforcement procedure leading to correlated equilibrium.
\newblock In G.~Debreu, W.~Neuefeind, and W.~Trockel, editors, {\em Economic
  Essays: A Festschrift for Werner Hildenbrand}, pages 181--200. Springer,
  2001.

\bibitem{HM03}
S.~Hart and A.~Mas-Colell.
\newblock Uncoupled dynamics do not lead to nash equilibrium.
\newblock {\em American Economic Review}, 93(5):1830--1836, December 2003.

\bibitem{HS84}
D.~P. Heyman and M.~J. Sobel.
\newblock {\em Stochastic Models in Operations Research}, volume~2.
\newblock McGraw-Hill, 1984.

\bibitem{HL11}
N.~Higham and L.~Lin.
\newblock On pth roots of stochastic matrices.
\newblock {\em Linear Algebra and its Applications}, 435(3):448--463, 2011.

\bibitem{HL04}
D.~Hunter and K.~Lange.
\newblock A tutorial on {MM} algorithms.
\newblock {\em The American Statistician}, 58(1):30--37, 2004.

\bibitem{KT79}
D.~Kahneman and A.~Tversky.
\newblock Prospect theory: An analysis of decision under risk.
\newblock {\em Econometrica}, pages 263--291, 1979.

\bibitem{KLM07}
L.~Karp, I.H. Lee, and R.~Mason.
\newblock A global game with strategic substitutes and complements.
\newblock {\em Games and Economic Behavior}, 60:155--175, 2007.

\bibitem{KK77}
J.~Keilson and A.~Kester.
\newblock {Monotone matrices and monotone {M}arkov processes}.
\newblock {\em Stochastic Processes and their Applications}, 5(3):231--241,
  1977.

\bibitem{Kri08}
V.~Krishnamurthy.
\newblock Decentralalized activation in dense sensor networks via global games.
\newblock {\em IEEE Transactions on Signal Processing}, 56(10):4936--4950,
  2008.

\bibitem{Kri09}
V.~Krishnamurthy.
\newblock Decentralized spectrum access amongst cognitive radios-an interacting
  multivariate global game-theoretic approach.
\newblock {\em IEEE Transactions on Signal Processing}, 57(10):3999--4013, Oct.
  2009.

\bibitem{KV12}
V.~Krishnamurthy and F.~Vazquez Abad.
\newblock Gradient based policy optimization of constrained unichain {M}arkov
  decision processes.
\newblock In S.~Cohen, D.~Madan, and T.~Siu, editors, {\em Stochastic
  Processes, Finance and Control: A Festschrift in Honor of Robert J. Elliott}.
  World Scientific, 2012.
\newblock http://arxiv.org/abs/1110.4946.

\bibitem{KB16}
V.~Krishnamurthy and S.~Bhatt.
\newblock Sequential detection of market shocks with risk-averse cvar social
  sensors.
\newblock {\em IEEE Journal Selected Topics in Signal Processing}, 2016.

\bibitem{KE97a}
V.~Krishnamurthy and R.J. Elliott.
\newblock Filters for estimating {M}arkov modulated poisson processes and image
  based tracking.
\newblock {\em Automatica}, 33(5):821--833, May 1997.

\bibitem{KH15}
V.~Krishnamurthy and W.~Hoiles.
\newblock Online reputation and polling systems: Data incest, social learning
  and revealed preferences.
\newblock {\em IEEE Transactions Computational Social Systems}, 1(3):164--179,
  Jan. 2015.

\bibitem{KMY08}
V.~Krishnamurthy, M.~Maskery, and G.~Yin.
\newblock Decentralized activation in a {Z}ig{B}ee-enabled unattended ground
  sensor network: A correlated equilibrium game theoretic analysis.
\newblock {\em IEEE Transactions on Signal Processing}, 56(12):6086--6101,
  December 2008.

\bibitem{KP15}
V.~Krishnamurthy and U.~Pareek.
\newblock Myopic bounds for optimal policy of {POMDPs}: An extension of
  {L}ovejoy's structural results.
\newblock {\em Operations Research}, 62(2):428--434, 2015.

\bibitem{KR14}
V.~Krishnamurthy and C.~Rojas.
\newblock Reduced complexity {HMM} filtering with stochastic dominance bounds:
  A convex optimization approach.
\newblock {\em IEEE Transactions on Signal Processing}, 62(23):6309--6322,
  2014.

\bibitem{KW09}
V.~Krishnamurthy and B.~Wahlberg.
\newblock {POMDP} multiarmed bandits -- structural results.
\newblock {\em Mathematics of Operations Research}, 34(2):287--302, May 2009.

\bibitem{KV86}
P.~R. Kumar and P.~Varaiya.
\newblock {\em Stochastic systems -- Estimation, Identification and Adaptive
  Control}.
\newblock Prentice-Hall, New Jersey, 1986.

\bibitem{LI99}
A.~Logothetis and A.~Isaksson.
\newblock On sensor scheduling via information theoretic criteria.
\newblock In {\em Proc. American Control Conf.}, pages 2402--2406, San Diego,
  1999.

\bibitem{LK99}
A.~Logothetis and V.~Krishnamurthy.
\newblock Expectation maximization algorithms for {MAP} estimation of jump
  {M}arkov linear systems.
\newblock {\em IEEE Transactions on Signal Processing}, 47(8):2139--2156,
  August 1999.

\bibitem{Lor71}
G.~Lorden.
\newblock Procedures for reacting to a change in distribution.
\newblock {\em The Annals of Mathematical Statistics}, pages 1897--1908, 1971.

\bibitem{Lov87}
W.~S. Lovejoy.
\newblock Some monotonicity results for partially observed {M}arkov decision
  processes.
\newblock {\em Operations Research}, 35(5):736--743, Sept.-Oct. 1987.

\bibitem{Lov93}
W.~S. Lovejoy.
\newblock Suboptimal policies with bounds for parameter adaptive decision
  processes.
\newblock {\em Operations Research}, 41(3):583--599, 1993.

\bibitem{MXJ00}
J.~Ma, L.~Xu, and M.~I. Jordan.
\newblock Asymptotic convergence rate of the em algorithm for gaussian
  mixtures.
\newblock {\em Neural Computation}, 12(12):2881--2907, 2000.

\bibitem{MS99a}
C.~D. Manning and H.~Sch{\"{u}}tze.
\newblock {\em Foundations of Statistical Natural Language Processing}.
\newblock The MIT Press, 1999.

\bibitem{MWG95}
A.~Mas-Colell, M.~Whinston, and J.~Green.
\newblock {\em Microeconomic Theory}.
\newblock Oxford, 1995.

\bibitem{MKW15}
R.~Mattila, V.~Krishnamurthy, and B.~Wahlberg.
\newblock Recursive identification of chain dynamics in hidden markov models
  using non-negative matrix factorization.
\newblock In {\em Proceedings of IEEE CDC 2015}, 2015.

\bibitem{MS00}
S.~Morris and H.S. Shin.
\newblock Global games: Theory and applications.
\newblock In {\em Advances in Economic Theory and Econometrics: Proceedings of
  Eight World Congress of the Econometric Society}, pages 56--114. Cambridge
  University Press, 2000.

\bibitem{MV78}
H.~Moulin and J.-P. Vial.
\newblock Strategically zero-sum games: The class of games whose completely
  mixed equilibria cannot be improved upon.
\newblock {\em International Journal of Game Theory}, 7(3-4):201--221, 1978.

\bibitem{Mou86}
G.~B. Moustakides.
\newblock Optimal stopping times for detecting changes in distributions.
\newblock {\em Annals of Statistics}, 14:1379--1387, 1986.

\bibitem{MS02}
A.~Muller and D.~Stoyan.
\newblock {\em Comparison Methods for Stochastic Models and Risk}.
\newblock Wiley, 2002.

\bibitem{NJ13}
M.~Naghshvar and T.~Javidi.
\newblock Active sequential hypothesis testing.
\newblock {\em The Annals of Statistics}, 41(6):2703--2738, 2013.

\bibitem{Nak95}
T.~Nakai.
\newblock The problem of optimal stopping in a partially observable markov
  chain.
\newblock {\em Journal of Optimization Theory and Applications},
  45(3):425--442, 1985.

\bibitem{NKY13}
O.~Namvar, V.~Krishnamurthy, and G.~Yin.
\newblock Distributed tracking of correlated equilibria in regime switching
  noncooperative games.
\newblock {\em IEEE Transactions on Automatic Control}, 58(10):2435--2450,
  2013.

\bibitem{Nas51}
J.~Nash.
\newblock Non-cooperative games.
\newblock {\em Annals of Mathematics}, 54(2):286--295, Sep. 1951.

\bibitem{NCH04}
R.~Nau, S.~Canovas, and P.~Hansen.
\newblock On the geometry of {N}ash equilibria and correlated equilibria.
\newblock {\em International Journal of Game Theory}, 32(4):443--453, 2004.

\bibitem{Neu89}
M.~F. Neuts.
\newblock {\em Structured stochastic matrices of M/G/1 type and their
  applications}.
\newblock Marcel Dekker, N.Y., 1989.

\bibitem{SM04}
R.~Olfati-Saber and R.~M. Murray.
\newblock Consensus problems in networks of agents with switching topology and
  time-delays.
\newblock {\em IEEE Transactions on Automatic Control}, 49(9):1520--1533, Sept.
  2004.

\bibitem{PC08}
T.~Park and G.~Casella.
\newblock The {B}ayesian lasso.
\newblock {\em Journal of the American Statistical Association},
  103(482):681--686, 2008.

\bibitem{PT12}
A.~Polunchenko and A.~Tartakovsky.
\newblock State-of-the-art in sequential change-point detection.
\newblock {\em Methodology and computing in applied probability},
  14(3):649--684, 2012.

\bibitem{PH08}
H.~V. Poor and O.~Hadjiliadis.
\newblock {\em Quickest Detection}.
\newblock Cambridge University Press, 2008.

\bibitem{QS09}
J.~Quah and B.~Strulovici.
\newblock Comparative statics, informativeness, and the interval dominance
  order.
\newblock {\em Econometrica}, 77(6):1949--1992, 2009.

\bibitem{QS12}
J.~Quah and B.~Strulovici.
\newblock Aggregating the single crossing property.
\newblock {\em Econometrica}, 80(5):2333--2348, 2012.

\bibitem{RO05}
I.~Rapoport and Y.~Oshman.
\newblock A {Cram\'er-Rao}-type estimation lower bound for systems with
  measurement faults.
\newblock {\em IEEE Transactions on Automatic Control}, 50(9):1234--1245, 2005.

\bibitem{Rob51}
J.~Robinson.
\newblock An iterative method of solving a game.
\newblock {\em Annals of Mathematics}, 54(2):296--301, Sep. 1951.

\bibitem{Ros83}
S.~Ross.
\newblock {\em Introduction to Stochastic Dynamic Programming}.
\newblock Academic Press, San Diego, California., 1983.

\bibitem{Sha88}
J.~G. Shanthikumar.
\newblock {DFR} property of first-passage times and its preservation under
  geometric compounding.
\newblock {\em The Annals of Probability}, pages 397--406, 1988.

\bibitem{Sim03}
C.~Sims.
\newblock Implications of rational inattention.
\newblock {\em Journal of Monetary Economics}, 50(3):665--690, 2003.

\bibitem{SS97}
L.~Smith and P.~Sorensen.
\newblock Informational herding and optimal experimentation.
\newblock Economics Papers 139, Economics Group, Nuffield College, University
  of Oxford, 1997.

\bibitem{SSD93}
D.~D. Sworder, P.~F. Singer, D.~Doria, and R.~G. Hutchins.
\newblock Image-enhanced estimation methods.
\newblock {\em Proceedings of the IEEE}, 81(6):797--812, June 1993.

\bibitem{TM10}
A.~Tartakovsky and G.~Moustakides.
\newblock State-of-the-art in bayesian changepoint detection.
\newblock {\em Sequential Analysis}, 29(2):125--145, 2010.

\bibitem{Tib96}
R.~Tibshirani.
\newblock Regression shrinkage and selection via the lasso.
\newblock {\em Journal of the Royal Statistical Society. Series B
  (Methodological)}, pages 267--288, 1996.

\bibitem{TMN98}
P.~Tichavsky, C.~H. Muravchik, and A.~Nehorai.
\newblock Posterior {C}ram\'er-{R}ao bounds for discrete-time nonlinear
  filtering.
\newblock {\em IEEE Transactions on Signal Processing}, 46(5):1386--1396, May
  1998.

\bibitem{Van68}
H.L.~Van Trees.
\newblock {\em {Detection, Estimation and Modulation Theory}}.
\newblock John Wiley \& Sons, 1968.

\bibitem{TAM12}
V.~Tzoumas, C.~Amanatidis, and E.~Markakis.
\newblock A game-theoretic analysis of a competitive diffusion process over
  social networks.
\newblock In {\em Internet and Network Economics}, volume 7695, pages 1--14.
  Springer, 2012.

\bibitem{VB13}
V.~Veeravalli and T.~Banerjee.
\newblock Quickest change detection.
\newblock {\em Academic press library in signal processing: Array and
  statistical signal processing}, 3:209--256, 2013.

\bibitem{Whi82}
W.~Whitt.
\newblock Multivariate monotone likelihood ratio and uniform conditional
  stochastic order.
\newblock {\em Journal Applied Probability}, 19:695--701, 1982.

\bibitem{XJ96}
L.~Xu and M.~I. Jordan.
\newblock On convergence properties of the em algorithm for gaussian mixtures.
\newblock {\em Neural computation}, 8(1):129--151, 1996.

\end{thebibliography}
